\documentclass[12pt]{article}

\usepackage{amsmath, amsthm, amssymb}
\usepackage{setspace}
\usepackage{graphicx}
\def\R{{\mathbb R}}
\def\N{{\mathbb N}}
\def\bpsi{\boldsymbol\psi}
\def\bd{\boldsymbol d}
\def\bV{\boldsymbol V}
\newtheorem {definition} {Definition}

  \newtheorem {theorem}    {Theorem}
  \newtheorem {lemma}      {Lemma}
  \newtheorem {corollary}  {Corollary}
  \newtheorem {proposition}{Proposition}

\numberwithin{equation}{section}
\begin{document}
\title{Existence of a  solution to a  fluid-multi-layered-structure interaction problem}
\author{Boris Muha\\
Department of Mathematics\\ University of Zagreb
\and Sun\v{c}ica \v{C}ani\'{c}\\
Department of Mathematics\\ University of Houston} 
\date{}
\maketitle

\begin{abstract}
We study a nonlinear, unsteady, moving boundary, fluid-structure (FSI) problem in which the structure is composed of two layers:
a thin layer which is in contact with the fluid, and a thick layer which sits on top of the thin structural layer. 
The fluid flow, which is driven by the time-dependent dynamic pressure data, 
is governed by the 2D Navier-Stokes equations for an incompressible, viscous fluid,
defined on a 2D cylinder. 
The elastodynamics of the cylinder wall is governed by the 1D linear wave equation modeling the thin structural layer, 
and by the 2D equations of linear elasticity modeling the thick structural layer.  
The fluid and the structure, as well as the two structural layers, are fully coupled via the kinematic 
and dynamic coupling conditions describing continuity of velocity and balance of contact forces. 
The thin structural layer acts as a fluid-structure interface with mass.
The resulting FSI problem is a nonlinear moving boundary problem of parabolic-hyperbolic type.
This problem is motivated by the flow of blood in elastic arteries whose walls are composed of several layers, each with
different mechanical characteristics and thickness. 
We prove existence
of a weak solution to this nonlinear FSI problem as long as the cylinder radius is greater than zero.
The proof is based on a novel semi-discrete, operator splitting numerical scheme, 
known as the kinematically coupled scheme.
%The backbone of the kinematically coupled scheme is the well-known Marchuk-Yanenko scheme, 
%also known as the Lie splitting scheme.
We effectively prove convergence of that numerical scheme to a solution of the nonlinear fluid-multi-layered-structure interaction problem. 
The spaces of weak solutions presented in this manuscript reveal a striking new feature:
the presence of a thin fluid-structure interface with mass regularizes solutions of the coupled problem. 
\end{abstract}

\section{Introduction}\label{sec:intro}
\subsection{Problem definition}
We consider the flow  
of an incompressible, viscous fluid modeled by 
the Navier-Stokes equations in a 2D,  time-dependent cylindrical fluid domain $\Omega_F(t)$, which is not known {\sl a priori}:
\begin{equation}
{\bf FLUID:} \qquad \left .
\begin{array}{rcl}
\rho_F (\partial_t{\bf u}+{\bf u}\cdot\nabla{\bf u})&=&\nabla\cdot\boldsymbol\sigma,
\\
\nabla\cdot{\bf u}&=&0,
\end{array}
\right \}\ {\rm in}\ \Omega_F(t),\ t\in (0,T),
\label{NS}
\end{equation}
where $\rho_F$ denotes the fluid density; ${\bf u}$ the fluid velocity; 
$\boldsymbol\sigma=-p{\bf I}+2\mu{\bf D}({\bf u})$ is the fluid Cauchy stress tensor; $p$ is the fluid pressure;
$\mu$ is the kinematic viscosity coefficient;  and ${\bf D}({\bf u})=\frac 1 2(\nabla{\bf u}+\nabla^{\tau}{\bf u})$ is the symmetrized gradient of ${\bf u}$.

\begin{figure}[ht]
\centering{
\includegraphics[scale=0.35]{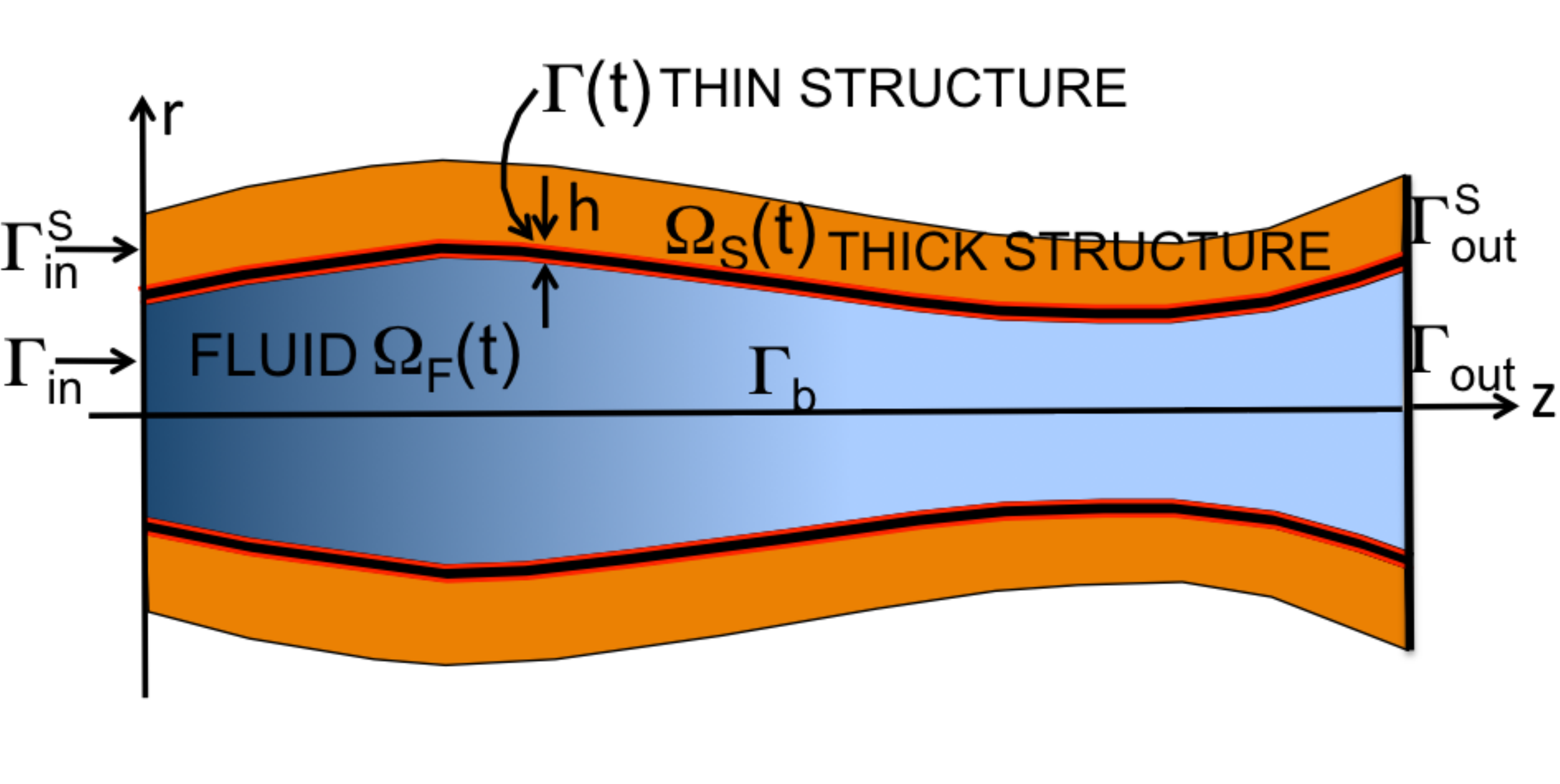}
}
\caption{Domain sketch.}
\label{fig:domain}
\end{figure}

The cylindrical fluid domain is of length $L$, with reference radius $r = R$. 
The radial (vertical) displacement of the cylinder radius at time $t$ and position $z\in(0,L)$ will be denoted by $\eta(t,z)$, giving rise
to a deformed domain with radius $R+\eta(t,z)$.
For simplicity, we will be assuming that longitudinal displacement of the structure is negligible.
This is a common assumption in literature on FSI in blood flow. 
Thus, the fluid domain, sketched in Figure~\ref{fig:domain}, is given by
$$
\Omega_F(t)=\{(z,r)\in\R^2:z\in (0,L),\ r\in (0,R+\eta(t,z)\},
$$
where the lateral boundary of the cylinder corresponds to fluid-structure interface,
denoted by 
$$\Gamma(t) = \{(z,r)\in\R^2:z\in (0,L),\ r = R + \eta(t,z)\}.$$
Without loss of generality we only consider the upper half of the fluid cylinder, with a symmetry boundary condition 
prescribed at the axis of symmetry,
denoted by $\Gamma_b$ in Figure~\ref{fig:domain}. 

The fluid is in contact with a thin elastic structure, which is located between the fluid and the thick structural layer.
The thin structure thereby serves as  a fluid-structure interface with mass.  
In this manuscript we will be assuming that the elastodynamics of the thin elastic structure is governed by the
1D wave equation
\begin{equation}
{\bf THIN \ STRUCTURE: }\qquad  \rho_{s_1} h \ \partial_{tt} \eta = c^2 \partial_{zz} \eta + f, \quad z\in (0,L),\ t\in (0,T),
\label{Koiter}
\end{equation}
where $\eta$ denotes radial (vertical) displacement.
More generally, the wave equation can be viewed as a special case of 
 the linearly (visco)elastic cylindrical Koiter shell model
\begin{equation}\label{Koiter_full}
\rho_{s_1} h {\partial^2_t \eta}+C_0 \eta -C_1{\partial^2_z \eta}+C_2\partial^4_z\eta+D_0{\partial_t \eta}-
D_1{\partial_t\partial^2_z\eta}+D_2\partial_t\partial^4_z\eta  =  f, 
\end{equation}
with  $C_0=C_2=D_0=D_1=D_2=0$.
Here, $\rho_{s_1}$ is structure density, $h$ denotes structure thickness, and $f$ denotes force density in the 
radial (vertical) direction acting on the structure.
The constants $C_i$ and $D_i > 0$ are the material constants describing structural elasticity and viscosity, respectively,
which are given in terms of  four material parameters: the Young's modulus of elasticity $E$,
the Poisson ratio $\sigma$,
and their viscoelastic couter-parts 
 (for a derivation of this model and the exact form of the coefficients, please see \cite{MarSun,SunTam}).
 The results of the manuscript hold in the case when all the coefficients  in the Koiter shell model are different from zero.
 From the analysis point of view, however, the most difficult case is the case of the wave equation,
 and for that reason we consider this case in the present manuscript.

The thick structural layer will be modeled by the equations of linear elasticity
\begin{equation}
{\bf THICK \ STRUCTURE:}\qquad  \rho_{s_2} \ {\partial_{tt}}  {\bf d} = \nabla \cdot {\bf S} \quad \textrm{in}\; {\Omega}_{S},\ t \in (0,T),
 \label{thick_structure}
\end{equation}
where  ${\bf d}(t,z,r) = (d_z(t,z,r), d_r(t,z,r))$ denotes structural displacement of the thick elastic wall at 
point $(z,r)\in\Omega_S$ and time $t$, 
$\boldsymbol S$ is the first Piola-Kirchhoff stress tensor, and $\rho_{s_2}$ is the density of the thick structure. 
Equation \eqref{thick_structure} describes the second Newton's Law of motion for an arbitrary thick structure.
We are interested in linearly elastic structures for which 
\begin{equation}
{\bf S}=\mu\ (\nabla {\bf d} + (\nabla {\bf d})^T)+\lambda(\nabla\cdot{\bf d}){\bf I},
\label{PKST}
\end{equation} 
where $\lambda$ and $\mu$ are Lam\' e constants describing material properties of the structure.
Since structural problems are typically defined in the Lagrangian framework, domain $\Omega_S$ corresponds to a fixed, reference domain 
which is independent of time, and is
given by 
$$\Omega_S=(0,L)\times (R,R+H).$$
A deformation of $\Omega_S$ at time $t$ is denoted by $\Omega_S(t)$ in Figure~\ref{fig:domain}.

{\bf The coupling} between the fluid, the thin structural layer, and the thick structural layer is achieved via two sets of coupling conditions:
the kinematic coupling condition and the dynamic coupling condition. In the present manuscript the kinematic coupling condition
is the   {\bf no-slip} boundary condition  between the fluid and thin structure, as well as between the thin and thick structural layers. 
Depending on the application,
different kinematic coupling conditions can be prescribed between the three 
different physical models.

The dynamic coupling condition describes balance of forces at the 
fluid-structure interface $\Gamma(t)$. Since $\Gamma(t)$ is a fluid-structure interface with mass, the dynamic coupling condition
states that the mass times the acceleration of the interface is balanced by the sum of total forces acting on $\Gamma(t)$. This 
includes the contribution due to the elastic energy of the structure ($\partial_{zz}\eta$), and the balance of contact forces exerted
by the fluid and the thick structure onto $\Gamma(t)$. More precisely, we have the following set of coupling conditions 
written in Lagrangian framework, with $z\in (0,L)$ and  $t\in (0,T)$:
\begin{itemize}
\item {\bf The\  kinematic\  coupling\ condition:}
\begin{equation}\label{Coupling1a}
\begin{array}{cl}
\displaystyle{({\partial_t \eta}(t,z),0)} = \boldsymbol{u}(t,z,R+\eta(t,z)),&\ {\rm (continuity \ of \ velocity)}
\\
\displaystyle{(\eta(t,z),0)={\bf d}(t,z,R)},&\ {\rm (continuity \ of \ displacement)}
\end{array}
\end{equation}
\item {\bf The\  dynamic\  coupling\ condition:}
\begin{equation}
\label{Coupling1b}
\rho_{s_1} h \partial_{tt}\eta =  c^2 \partial_{zz}\eta -J(t,z)(\sigma{\bf n})|_{(t,z,R+\eta(t,z))}\cdot{\bf e}_r+{\bf S}(t,z,R){\bf e}_r\cdot{\bf e}_r.
\end{equation}
Here 
$J(t,z)=\displaystyle{\sqrt{1+({\partial_z \eta}(t,z))^2}}$
denotes the Jacobian of the transformation from  Eulerian to Lagrangian coordinates, and ${\bf e}_r$ is the 
unit vector  associated with the vertical, $r$-direction.  
\end{itemize}

Problem \eqref{NS}-\eqref{Coupling1b} is supplemented with initial and boundary conditions. 
At the {\bf inlet and outlet boundaries} to the fluid domain we prescribe zero tangential velocity and a given dynamic pressure 
$p+\frac{\rho_f}{2}|u|^2$ (see e.g. \cite{CMP}):
\begin{equation}\label{PD}
\left. \begin{array}{rcl}
\displaystyle{p+\frac{\rho_f}{2}|u|^2}&=&P_{in/out}(t), \\
u_r &=& 0,\\
\end{array}
\right\} \quad {\rm on}\ \Gamma_{in/out},
\end{equation}
where $P_{in/out}\in L^{2}_{loc}(0,\infty)$ are given. 
Therefore, the fluid flow is driven by a prescribed dynamic pressure drop, 
and the flow enters and leaves the fluid domain orthogonally to the inlet and outlet boundary. 

At the {\bf bottom boundary} we prescribe the symmetry boundary condition:
\begin{equation}
u_r=\partial_r u_z=0,\quad {\rm on}\ \Gamma_b.
\label{SymBC}
\end{equation}

At the {\bf end points of the thin structure} we prescribe zero displacement:
\begin{equation}\label{BCthin}
\eta(t,0)=\eta(t,L)=0.
\end{equation}

For the {\bf thick structure}, we assume that the {\bf external (top) boundary} $r=H$ 
is exposed to an external ambient pressure $P_e$:
\begin{equation}
{\bf S}{\bf e}_r= -P_e {\bf e}_r,\quad {\rm on}\ \Gamma_{ext},
\label{DSBC}
\end{equation}
while at the {\bf end points} of the annular sections of the thick structure, $\Gamma^s_{in/out}$, we assume that the displacement is zero
$$
{\bf d}(t,0,r)={\bf d}(t,L,r) = 0, \quad {\rm for}\  r \in (R,H).
$$

The {\bf initial} fluid and structural velocities, and the initial displacements are given by
\begin{equation}
{\bf u}(0,.)={\bf u}_0,\ \eta(0,.)=\eta_0,\ \partial_t\eta(0,.)=v_0,\; {\bf d}(0,.)={\bf d}_0,\; \partial_t{\bf d}(0,.)={\bf V}_0, 
\label{IC}
\end{equation}
and are assumed to belong to the following spaces: 
${\bf u}_0\in L^2(\Omega_F(0))$, $\eta_0\in H^1_0(0,1)$, 
$v_0\in L^2(0,1)$,  $\bV_0\in L^2(\Omega_S)$, $\bd_0\in H^1(\Omega_S)$,
satisfying the following compatibility conditions:
\begin{equation}
\begin{array}{c}
(\eta_0(z),0)={\bf d}_0(z,R),
\\
\eta_0(0)=\eta_0(L)=v_0(0)=v_0(L)=0={\bf d}_0(0,.)={\bf d}_0(L,.)={\bf V}_0(0,.)={\bf V}_0(L,.),
\\
R+\eta_0(z) > 0,\quad z\in [0,L].
\end{array}
\label{CC}
\end{equation}

We study the existence of a weak solution to the nonlinear FSI problem 
\eqref{NS}-\eqref{CC}, in which the flow is driven by the time-dependent inlet and outlet dynamic pressure data.

For simplicity, in the rest of the manuscript, we will be setting all the parameters in the problem to be equal to 1. 
This includes the domain parameters $R$ and $L$, the Lam\'{e} constants $\lambda$ and $\mu$, and the 
structure parameters $\rho_{s_1}, \rho_{s_2}$ and $h$. 
Furthermore, we will be assuming that the external pressure, given in \eqref{DSBC}, is equal to zero. 
Correspondingly, we  subtract the constant external pressure data from the inlet and outlet dynamic pressure data
to obtain an equivalent problem. 

\subsection{Motivation}
\begin{figure}[ht]
\centering{
\includegraphics[scale=0.25]{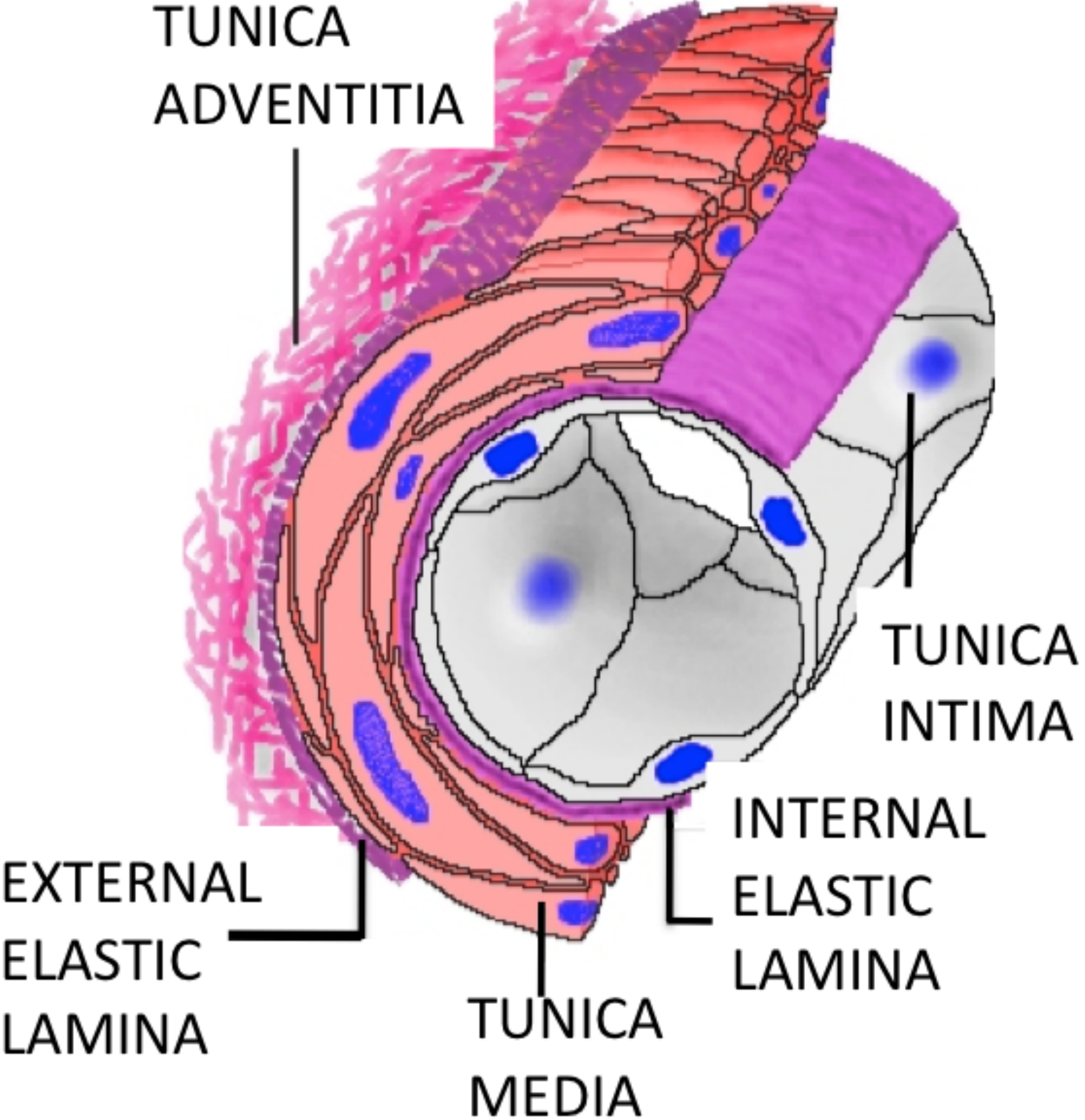}
}
\caption{Arterial wall layers.}
\label{fig:aretria}
\end{figure}

This work was motivated by blood flow in major human arteries. 
In medium-to-large human arteries, such as the aorta or coronary arteries, blood can be modeled as an incompressible,
viscous, Newtonian fluid. Arterial walls of major arteries are composed of 
several layers, each with 
different mechanical characteristics. The main layers are the tunica intima, the tunica media, and the tunica adventitia. 
They are separated by the thin elastic laminae,
see Figure~\ref{fig:aretria}.
To this date, there have been no fluid-structure interaction models or computational solvers 
of arterial flow that take into account the multi-layered structure of arterial walls.
In this manuscript we take a first step in this direction by proposing to study a benchmark problem in fluid-multi-layered-structure interaction
in which the structure consists of two layers, a thin and a thick layer, where the thin layer serves as a fluid-structure interface with mass.
The proposed problem is a nonlinear moving-boundary problem of parabolic-hyperbolic type for which the questions of
well-posedness  and numerical simulation
are wide open.

\subsection{Literature review}
Fluid-structure interaction problems have been extensively studied for the past 20 years by many authors. The focus has been
exclusively  on FSI problems with structures consisting of a single material. The field has evolved
from first studying FSI between an incompressible, viscous fluid and a rigid structure immersed in a fluid, to
considering compliant (elastic/viscoelastic) structures interacting with a fluid. Concerning compliant structures, 
the coupling between the structure and the fluid was first assumed to take place along a fixed fluid domain boundary
(linear coupling). 
This was then extended to FSI problems in which the coupling was evaluated at a deformed fluid-structure interface,
giving rise to an additional nonlinearity in the problem (nonlinear coupling).

Well-posedness results in which the structure was assumed to be a  rigid body immersed in a fluid, 
or described by a finite number of modal functions,
were studied in \cite{Boulakia,CSJ,CumTak,DE2,DEGLT,Feireisl,Gal,Starovoitov}.
FSI problems coupling the Navier-Stokes equations with linear elasticity  where the coupling was calculated at a fixed fluid domain boundary,
were considered in \cite{Gunzburger}, and  in \cite{BarGruLasTuff2,BarGruLasTuff,KukavicaTuffahaZiane}
where an additional nonlinear coupling term was added at the interface.
A study of well-posedness for FSI problems between an incompressible, viscous fluid and an elastic/viscoelastic structure
with nonlinear coupling evaluated at a moving interface started with the result by daVeiga  \cite{BdV1},
where existence of a strong solution was obtained locally in time 
for an interaction between a $2D$ fluid and a $1D$ viscoelastic string,
assuming periodic boundary conditions.
This result was extended by
Lequeurre in \cite{Leq11,Leq13}, where the existence of a unique, local in time, strong solution
for any data, and the existence of a global strong solution for small data, was proved in the case when the structure was modeled as a clamped viscoelastic beam.
D.~Coutand and S.~Shkoller proved existence, locally in time, of a unique, regular solution for
an interaction between a viscous, incompressible fluid in $3D$ and a $3D$ structure, immersed in the fluid,
where the structure was modeled by the equations of linear \cite{CSS1}, or quasi-linear \cite{CSS2} elasticity. 
In the case when the structure (solid) is modeled by a linear wave equation, 
I. Kukavica and A. Tufahha proved the existence, locally in time, of a strong solution, 
assuming lower regularity  for the initial data \cite{Kuk}. 
A similar result for compressible flows can be found in \cite{KukavicaNSLame}.
 A fluid-structure interaction between a viscous, incompressible fluid in $3D$, and $2D$ elastic shells was considered in 
\cite{ChenShkoller,ChengShkollerCoutand} where existence, locally in time,
of a unique regular solution was proved. 
All the above mentioned existence results for strong solutions are local in time.
We also mention that the works of Shkoller et al., and Kukavica at al. were obtained in the context 
of Lagrangian coordinates, which were used for both the structure and fluid problems.

In the context of weak solutions, the following results have been obtained.
Continuous dependence of  weak solutions on initial data for a
fluid structure interaction problem with a free boundary type coupling condition was studied in \cite{GioPad}.
Existence of a weak solution for a FSI problem between a $3D$ incompressible, viscous fluid 
and a $2D$ viscoelastic plate was considered by Chambolle et al. in \cite{CDEM}, while
Grandmont improved this result in \cite{CG}  to hold for
a $2D$ elastic plate. These results were extended to a more general geometry in \cite{LenRuz},
and then to the case of generalized Newtonian fluids in \cite{Len12}, and to a non-Newtonian shear dependent fluid in \cite{NecasovaLukacova}.
In these works existence of a weak solution was proved for as long as the elastic boundary does not touch "the bottom" (rigid)
portion of the fluid domain boundary. 

Muha and \v{C}ani\'{c}  recently proved existence of weak solutions to a class of FSI problems modeling the flow
of an incompressible, viscous, Newtonian fluid flowing through a cylinder whose lateral wall was modeled by 
either the linearly viscoelastic, or by the linearly elastic Koiter shell equations \cite{BorSun}, assuming nonlinear coupling at
the deformed fluid-structure interface. The fluid flow boundary conditions  were not periodic, but rather, the flow
was driven by the dynamic pressure drop data.
The methodology of proof in \cite{BorSun} was based on a semi-discrete, operator splitting Lie scheme, 
which was used in \cite{GioSun} to design a 
stable, loosely coupled partitioned numerical scheme, called the kinematically coupled scheme (see also \cite{MarSun}).
Ideas based on the Lie operator splitting scheme were also used by Temam
in \cite{TemCar} to prove the existence of a solution to the nonlinear Carleman equation.

Since the kinematically-coupled scheme is modular, it is particularly suitable for dealing with problems in which the structure consists
of several layers, since modeling each additional layer can be accomplished by adding a new module to the partitioned scheme.
Indeed, in the present manuscript we use the kinematically coupled scheme to prove the existence of a weak solution 
to a fluid-multi-layered structure interaction problem described in \eqref{NS}-\eqref{CC}.
The method of proof, first introduced by the authors in \cite{BorSun}, is robust in the sense that it can be extended to 
the multi-layered structural case considered in this manuscript. 
In particular, the method presented in \cite{BorSun} can be adopted to prove the existence of a 
FSI solution when the thin structure is modeled as a linearly elastic Koiter shell, and the elastodynamics of the thick structure is described
by the equations of linear elasticity. 
In this manuscript we make  further progress in this direction by considering  the linear wave equation,
and not the full Koiter shell model for our thin structural model. This is a more difficult case since
the fourth-order flexural term $\partial^4_z$, which provides higher regularity of weak solutions in the Koiter shell problem, is not present in 
the pure membrane model described by the linear wave equation. As a result, the analysis is more involved,
and a non-standard version of the Trace Theorem (see \cite{BorisTrag} and Theorem \ref{TraceTm}) needs to be
used to obtain the existence result.

The existence proof  presented in this manuscript is constructive.
 To deal with the motion of the fluid domain we adopt the Arbitrary Lagrangian Eulerian (ALE) approach.
 We construct a sequence of approximate solutions to the problem written in 
{ALE  weak formulation}  by performing the {time-discretization via Lie operator splitting}.
At each time step,
the full FSI problem is split into a fluid and a structure sub-problem.
To achieve stability and convergence of the corresponding splitting scheme, the splitting
is performed in a special way in which the fluid sub-problem 
includes structure inertia  via a  "Robin-type" boundary condition.  
The fact that structure inertia are included implicitly in the fluid sub-problem, enabled us, in the present work, to get
appropriate energy estimates for the approximate solutions, independently of the size of the time discretization.
Passing to the limit, as the size of the time step converges to zero, is achieved by the use of compactness
arguments alla Simon, and by a careful construction of the appropriate test functions associated with moving domains.

Our analysis revealed a striking new result  that concerns solutions of coupled, multi-physics problems.
We found that the presence of 
a thin structure with mass at the fluid-structure interface, regularizes the FSI solution. 
In \cite{BorSun1D} it is shown
that this is not just a consequence of our mathematical approach, but a physical property of the problem. 

The partitioned, loosely coupled scheme used in the proof in this manuscript
has already been implemented in the design
of several stable computational FSI schemes for 
simulation of blood flow in human arteries with thin structural models \cite{MarSun,Fernandez,GioSun,Lukacova},
with a thick structural model \cite{ThickPaper},
and with a multi-layered structural model \cite{multi-layered}.
We effectively prove in this manuscript that this numerical scheme converges to a weak solution of the 
nonlinear fluid-multi-layered structure interaction problem.

\if 1 = 0
This manuscript is organized as follows. In Section~\ref{sec:problem}
we describe the underlying FSI problem and present the energy of the coupled problem.
In Section~\ref{sec:ALE_and Lie} we introduce the 
corresponding ALE formulation and  perform the Lie operator splitting, thereby introducing two sub-problems, one
for the fluid and one for the structure.
In Section~\ref{sec:weak_solutions} we define solution spaces and weak solutions.
The time discretization and the definition of approximate solutions for each $\Delta t > 0$, are introduced in
Section~\ref{sec:approximate_solutions}. In this section the main uniform energy estimates, independent of $\Delta t$, are obtained. 
Convergence of approximate solutions, as $\Delta t \to 0$, is discussed in Section~\ref{sec:convergence}.
Finally, in Section~\ref{sec:limit} we show that the limiting functions satisfy the weak form of the coupled FSI problem.
\fi

\section{The energy of the coupled problem}
We begin by first showing that the coupled FSI problem \eqref{NS}-\eqref{CC} is well-formulated in the sense that the 
total energy of the problem is bounded in terms of the prescribed data. More precisely, we now show that the following energy estimate holds:
\begin{equation}\label{EE}
\frac{d}{dt}\left (E_{kin}(t)+E_{el}(t)\right )+D(t)\leq C(P_{in}(t),P_{out}(t)),
\end{equation}
where 
\begin{equation}\label{energy}
\begin{array}{lcr}
E_{kin}(t) &:=&\frac 1 2\left ( \|{\bf u}\|^2_{L^2(\Omega_F(t))}+\|\partial_t \eta\|^2_{L^2(\Gamma)}
+\|{\partial_t}{\bf d}\|^2_{L^2(\Omega_S)}\right ),\\ \\
E_{el}(t) &:=&\frac 1 2 \left ( \|\partial_z \eta\|^2_{L^2(\Gamma)}
+2 \|{\bf D}({\bf d})\|^2_{L^2(\Omega_S)}+\|\nabla\cdot{\bf d}\|^2_{L^2(\Omega_S)}\right ),
\end{array}
\end{equation}
denote the kinetic and elastic energy of the coupled problem, respectively,
and the term $D(t)$ captures viscous dissipation in the fluid:
\begin{equation}\label{dissipation}
D(t) :=  \|{\bf D}({\bf u})\|^2_{L^2(\Omega_F(t))}.
\end{equation}
The constant $C(P_{in}(t),P_{out}(t)))$ depends only on the inlet and outlet pressure data, which are both functions of time.
Notice that, due to the presence of an elastic fluid-structure interface with mass, the kinetic energy term $E_{kin}(t)$ 
contains a contribution from the kinetic energy of the fluid-structure interface $\|\partial_t\eta\|^2_{L^2(\Gamma)}$ incorporating the 
interface inertia, and the
elastic energy $E_{el}(t)$ of the FSI problem accounts for the elastic energy $\|\partial_z\eta\|^2_{L^2(\Gamma)}$ of the interface.
If a FSI problem between the fluid and a thick structure was considered without the thin FSI interface with mass,
these terms would not have been present. In fact, the traces of the displacement and velocity at the fluid-structure interface
of that FSI problem would not have been even defined for weak solutions.

To show that \eqref{EE} holds, we first multiply equation (\ref{NS}) by ${\bf u}$, integrate over $\Omega_{F}(t)$, and formally integrate by parts to obtain:
$$
\int_{\Omega_{F}(t)}\big (\partial_t{\bf u}\cdot{\bf u}+({\bf u}\cdot\nabla){\bf u}\cdot{\bf u}\big )+2\int_{\Omega_{F}(t)}|{\bf D}{\bf u}|^2
-\int_{\partial\Omega_{F}(t)}(-p{\bf I}+2{\bf D}({\bf u})){\bf n}(t)\cdot{\bf u}=0.
$$
To deal with the inertia term we first recall that  $\Omega_F(t)$ is moving in time and that the velocity of the 
lateral boundary is given by ${\bf u}|_{\Gamma(t)}$.
The transport theorem applied to the first term on the left hand-side of the above equation then gives:
$$
\int_{\Omega_{F}(t)}\partial_t{\bf u}\cdot{\bf u} = \frac 1 2\frac{d}{dt}\int_{\Omega_{F}(t)}|{\bf u}|^2-
\frac 1 2\int_{\Gamma(t)}|{\bf u}|^2{\bf u}\cdot{\bf n}(t).
$$
The second term on the left hand side can can be rewritten by using integration by parts, and the divergence-free condition, to obtain:
$$
\int_{\Omega_{F}(t)}({\bf u}\cdot\nabla){\bf u}\cdot{\bf u}=\frac 1 2\int_{\partial\Omega_{F}(t)}|{\bf u}|^2{\bf u}\cdot{\bf n}(t)=
\frac 1 2\big (\int_{\Gamma(t)}|{\bf u}|^2{\bf u}\cdot{\bf n}(t) 
$$
$$
-\int_{\Gamma_{in}}|{\bf u}|^2u_z+\int_{\Gamma_{out}}|{\bf u}|^2u_z.
\big )
$$
These two terms added together give
\begin{equation}\label{inertia}
\int_{\Omega_{F}(t)}\partial_t{\bf u}\cdot{\bf u}+\int_{\Omega_{F}(t)}({\bf u}\cdot\nabla){\bf u}\cdot{\bf u}=
\frac 1 2\frac{d}{dt}\int_{\Omega_{F}(t)}|{\bf u}|^2 -\frac 1 2\int_{\Gamma_{in}}|{\bf u}|^2u_z+\frac 1 2\int_{\Gamma_{out}}|{\bf u}|^2u_z.
\end{equation}
Notice the importance of nonlinear advection  in canceling the cubic term $\int_{\Gamma(t)}|{\bf u}|^2{\bf u}\cdot{\bf n}(t)$!

To deal with the boundary integral over $\partial\Omega_F(t)$, we first notice that on $\Gamma_{in/out}$
the boundary condition (\ref{PD}) implies $u_r=0$. Combined with the divergence-free condition we obtain
$\partial_z u_z=-\partial_r u_r=0$.
Now, using the fact that the normal to $\Gamma_{in/out}$ is ${\bf n}=(\mp 1,0)$ we get:
\begin{equation}\label{in_out}
\int_{\Gamma_{in/out}}(-p{\bf I}+2{\bf D}({\bf u})){\bf n}\cdot{\bf u}=\int_{\Gamma_{in}}P_{in}u_z-\int_{\Gamma_{out}}P_{out}u_z.
\end{equation}
In a similar way, using the symmetry boundary conditions (\ref{SymBC}), we get:
$$
\int_{\Gamma_b}(-p{\bf I}+2{\bf D}({\bf u})){\bf n}\cdot{\bf u}=0.
$$

What is left is to calculate the remaining boundary integral over $\Gamma(t)$. 
To do this, consider  the wave equation \eqref{Koiter}, multiply it by $\partial_t \eta$,
and integrate by parts to obtain
\begin{eqnarray}\label{Koiter_energy}
&\displaystyle{\int_0^1f\partial_t\eta = \frac{1}{2}\frac{d}{dt}\|\partial_t\eta\|^2_{L^2(0,1)}}+\frac 1 2\|\partial_z\eta\|^2_{L^2(0,1)}.
\nonumber 
\end{eqnarray}
Next, consider the elasticity equation \eqref{thick_structure}, multiply it by $\partial_t{\bf d}$ and integrate by parts to obtain:
\begin{eqnarray}\label{ElasticEnergy}
&\displaystyle{ \frac{1}{2}\frac{d}{dt}\big (\|\partial_t{\bf d}\|^2_{L^2(\Omega_S)}}
+2\|{\bf D}({\bf d})\|^2_{L^2(\Omega_S)}+\|\nabla\cdot{\bf d}\|^2_{L^2(\Omega_S)}\big )=
-\int_0^1{\bf S}{\bf e}_r\cdot\partial_t{\bf d}.\\
\nonumber 
\end{eqnarray}
By enforcing the dynamic and kinematic coupling conditions 
\eqref{Coupling1a}, \eqref{Coupling1b}, we obtain 
\begin{equation}\label{dynamic_energy}
-\int_{\Gamma(t)}\sigma{\bf n}(t)\cdot{\bf u}=-\int_0^1J\sigma{\bf n}\cdot{\bf u}=\int_0^1(f-{\bf S}{\bf e}_r\cdot{\bf e}_r)\partial_t\eta.
\end{equation}
Finally, by combining  \eqref{dynamic_energy} with \eqref{Koiter_energy}, and by adding the remaining
contributions to the energy of the FSI problem calculated in equations \eqref{inertia} and \eqref{PD}, one obtains the following energy equality: 
\begin{eqnarray}
\nonumber
\displaystyle{\frac{1}{2}\frac{d}{dt}\int_{\Omega_{F}(t)}|{\bf u}|^2
+\frac{1}{2}\frac{d}{dt}\|\partial_t\eta\|^2_{L^2(0,1)}
+2\int_{\Omega_{F}(t)}|{\bf D}{\bf u}|^2
+\frac{1}{2}
\frac{d}{dt}\|\partial_z\eta\|^2_{L^2(0,1)}}
 \\ 
 \nonumber
 \\
 \label{EnergyEq}
+ \displaystyle{\frac 1 2\frac{d}{dt}\big (\|\partial_t{\bf d}\|^2_{L^2(\Omega_S)}
+2\|{\bf D}({\bf d})\|^2_{L^2(\Omega_S)}+\|\nabla\cdot{\bf d}\|^2_{L^2(\Omega_S)}  \big)=\pm P_{in/out}(t)\int_{\Gamma_{in/out}}u_z}
\end{eqnarray}
By using the trace inequality and Korn inequality one can estimate:
$$
|P_{in/out}(t)\int_{\Sigma_{in/out}}u_z|\leq C |P_{in/out}|\|{\bf u}\|_{H^1(\Omega_{F}(t))}\leq
\frac{C}{2\epsilon}|P_{in/out}|^2+\frac{\epsilon C}{2}\|{\bf D}({\bf u})\|^2_{L^2(\Omega_{F}(t))}.
$$
By choosing $\epsilon$ such that $\frac{\epsilon C}{2}\leq 1$ we get the energy inequality \eqref{EE}.

\section{The ALE formulation and Lie splitting}\label{sec:ALE_and Lie}
\subsection{First order ALE formulation}
Since we consider nonlinear coupling between the fluid and structure, the fluid domain changes in time.
To prove the existence of a weak solution to \eqref{NS}-\eqref{CC}
 it is convenient to map the fluid domain onto a fixed domain $\Omega_F$. 
 The structural problems are already defined on fixed domains since they are formulated in the Lagrangian framework.
 We follow the approach typical of numerical methods for FSI problems
and map our fluid  domain $\Omega_{F}(t)$ onto $\Omega_F$ by using an Arbitrary Lagrangian-Eulerian (ALE) mapping
\cite{MarSun,GioSun,donea1983arbitrary,quaini2007semi,QTV00}. 
We remark here that in our problem it is not convenient to use Lagrangian formulation for the fluid sub-problem, 
as is done in e.g., \cite{CSS2,ChenShkoller,Kuk}, since, in our problem,
the fluid domain consists of a fixed, control volume of a cylinder, with prescribed inlet and outlet pressure data, 
which does not follow Largangian flow.

We begin by defining a family of ALE mappings $A_{\eta}$ parameterized by $\eta$:
\begin{equation}
A_{\eta}(t):\Omega_F\rightarrow\Omega_{F}(t),\quad
A_{\eta}(t)(\tilde{z},\tilde{r}):=\left (\begin{array}{c}\tilde{z}\\(1+\eta(t,\tilde z))\tilde{r}\end{array}\right ),\quad (\tilde{z},\tilde{r})\in\Omega_F,
\label{RefTrans}
\end{equation} 
where $(\tilde z,\tilde r)$ denote the coordinates in the reference domain $\Omega_F=(0,1)\times (0,1)$.
The mapping $A_\eta(t)$ is a bijection, 
and its Jacobian is given by
\begin{equation}\label{ALE_Jacobian}
|{\rm det} \nabla A_\eta(t)| = |1 + \eta(t,\tilde z)|.
\end{equation}
Composite functions with the ALE mapping will be denoted by
\begin{equation}\label{u_eta}
{\bf u}^{\eta}(t,.)={\bf u}(t,.)\circ A_{\eta}(t) \quad {\rm and} \quad p^{\eta}(t,.)=p(t,.)\circ A_{\eta}(t).
\end{equation}
The derivatives of composite functions satisfy:
\begin{equation}
\label{dercomp}
\partial_t {\bf u}=\partial_t{\bf u}^{\eta}-({\bf w}^{\eta}\cdot\nabla^{\eta}){\bf u}^{\eta},\quad \nabla{\bf u}=\nabla^{\eta}{\bf u}^{\eta},
\end{equation}
where the ALE domain velocity, ${\bf w}^{\eta}$, and the transformed gradient, $\nabla^{\eta}$, are given by:
\begin{equation}
\displaystyle{{\bf w}^{\eta}=\partial_t\eta{\tilde r}{\bf e}_r,\quad \nabla^{\eta}=}
\left (\begin{array}{c}
\displaystyle{\partial_{\tilde{z}}-\tilde{r}\frac{\partial_z\eta}{1+\eta}\partial_{\tilde{r}}}
\\
\displaystyle{
\frac{1}{1+\eta}\partial_{\tilde{r}}}
\end{array}\right ).
\label{nablaeta}
\end{equation}
One can see that $
\nabla^{\eta}{\bf v}=\nabla{\bf v}(\nabla A_{\eta})^{-1}.
$
For the purposes of the existence proof we also introduce the following notation:
$$
 \sigma^{\eta}=-p^{\eta}{\bf I}+2{\bf D}^{\eta}({\bf u}^{\eta}),\quad  
{\bf D}^{\eta}({\bf u}^{\eta})=\frac 1 2(\nabla^{\eta}{\bf u}^{\eta}+(\nabla^{\eta})^{\tau}{\bf u}^{\eta}).
$$
We are now ready to rewrite problem \eqref{NS}-\eqref{CC} in the ALE formulation. 
However, before we do that, we will make one more important step in our strategy to prove the existence of a weak solution to \eqref{NS}-\eqref{CC}.
Namely, as mentioned earlier, we would like to ``solve''  the coupled FSI problem by approximating the problem using the
time-discretization via Lie operator splitting.
Since Lie operator splitting is defined for systems that are first-order in time,
see Section~\ref{sec:LieSplitting}, 
we have to replace the second-order time-derivatives of $\eta$ and ${\bf d}$, with the first-order time-derivatives 
of the thin and thick structure velocities, respectively. 
Furthermore, 
we will use the kinematic coupling condition \eqref{Coupling1a} which implies that
the fluid-structure interface velocity is equal to the normal trace of the fluid velocity on $\Gamma_\eta(t)$.
Thus, we will introduce
a new variable, $v$, to denote this trace, and will replace $\partial_t \eta$ by $v$ {\bf everywhere} in the structure equation.
This has deep consequences both for the existence proof presented in this manuscript, as well as for the 
proof of  stability of the underlying numerical scheme, 
as it enforces the kinematic coupling condition implicitly in all the steps of the scheme. 
We also introduce another new variable ${\bf V}=\partial_t{\bf d}$ 
which denotes the  thick structure velocity.
This enables us to rewrite problem \eqref{NS}-\eqref{CC} as a first-order system in time.

Thus, the ALE formulation of problem \eqref{NS}-\eqref{CC}, defined on the reference domain $\Omega_F$,
and written as a first-order system in time, is given by the following:

{\sl{
Find ${\bf u}(t,\tilde{z},\tilde{r}),p(t,\tilde{z},\tilde{r}),\eta(t,\tilde{z})$, $v(t,\tilde{z})$, ${\bf d}(t,\tilde{z})$ and ${\bf V}(t,\tilde{z})$ such that
}}
 
\begin{equation}\label{FSIeqRef}
\left.
\begin{array}{rcl}
{\partial_t{\bf u}}+(({\bf u}-{\bf w}^{\eta})\cdot\nabla^{\eta}){\bf u}&=&\nabla^{\eta}\cdot\sigma^{\eta}, \\
\nabla^{\eta}\cdot{\bf u}&=&0,
\end{array}
\right\} \textrm{in}\; (0,T)\times\Omega_F,  
\end{equation}

\begin{equation}
\left.
\begin{array}{rcl}
u_r&=&0,\\
{\partial_r u_z}&=&0
\end{array}
\right\} \textrm{on}\; (0,T)\times\Gamma_b,
\end{equation}

\begin{equation}\label{in_out}
\left.
\begin{array}{rcl}
p+\frac{1}{2}|u|^2&=&P_{in/out}(t),\\
u_r&=&0,
\end{array}
\right\} \textrm{on}\; (0,T)\times\Gamma_{in/out},
\end{equation}

\begin{equation}\label{structure_ref}
\left.
\begin{array}{rcl}
{\bf u}&=&v{\bf e}_r,  \\
{\bf d}&=&\eta{\bf e}_r, \\
\partial_t \eta&=&v,    \\
\partial_t v-\partial^2_z \eta&=&-J\sigma{\bf n}\cdot {\bf e}_r+{\bf S}{\bf e}_r\cdot{\bf e}_r
\end{array}
\right\} \textrm{on}\; (0,T)\times(0,1) ,
\end{equation}

\begin{equation}\label{thick_structure_ALE}
\left .
\begin{array}{rcl}
\partial_t{\bf d}&=&{\bf V},\\
\partial_t{\bf V}&=&\nabla\cdot{\bf S},
\end{array}
\right \}
\quad \textrm{on}\;  \Omega_S,
\end{equation}

\begin{equation}
\begin{array}{rcl}
\eta&=&0\; {\rm on}\; (0,T)\times \partial\Gamma,
\\
{\bf d}&=&0\;{\rm on}\; (0,T)\times \Gamma^s_{in/out}
\end{array}
\end{equation}

\begin{equation}\label{thickBC_ALE}
\begin{array}{rcl}
{\bf S}{\bf e}_r&=&0\; {\rm on}\; (0,T)\times \Gamma_{ext}.
\end{array}
\end{equation}
\begin{equation}\label{IC_ALE}
{\bf u}(0,.)={\bf u}_0,\eta(0,.)=\eta_0,  v(0,.)=v_0, {\bf d}(0,.)={\bf d}_0, {\bf V}(0,.)={\bf V}_0 \quad {\rm at}\quad t = 0.
\end{equation}
Here, we have dropped the superscript $\eta$ in ${\bf u}^\eta$ to simplify notation.
This defines a parabolic-hyperbolic-hyperbolic nonlinear moving boundary problem. 
The nonlinearity appears in the equations \eqref{FSIeqRef}, and in 
the coupling conditions \eqref{structure_ref} where the fluid quantities are evaluated at the deformed
fluid-structure interface $\eta(t,z)$.
Parabolic features are associated with the fluid problem \eqref{FSIeqRef}-\eqref{in_out},
 while hyperbolic features  come from the 2D equations of elasticity,
 and from the 1D wave equation modeling the fluid-structure interface, described by
 the last equation in \eqref{structure_ref}.

\subsection{The operator splitting scheme}\label{sec:LieSplitting}
To prove the existence of a weak solution to \eqref{FSIeqRef}-\eqref{IC_ALE} we use the time-discretization via operator splitting,
known as the Lie splitting or the Marchuk-Yanenko splitting scheme.
The underlying multi-physics problem will be split into the fluid and structure sub-problems,
following the different ``physics'' in the problem, but the splitting will be performed
in a particularly clever manner so that the resulting problem defines a scheme 
that converges to a weak solution of the continuous problem.
The basic ideas behind the Lie splitting can be summarized as follows.

Let $N\in\N$, $\Delta t=T/N$ and $t_n=n\Delta t$.
Consider the following
initial-value problem:
$$
\frac{d\phi}{dt}+A\phi=0\quad {\rm in}\ (0,T),\quad \phi(0)=\phi_0,
$$
where $A$ is an operator defined on a Hilbert space, and $A$ can be written as $A=A_1+A_2$. 
Set $\phi^0=\phi_0$, and, for $n=0,\dots, N-1$ and $i=1,2$, compute $\phi^{n+\frac i 2}$ by solving
\begin{equation*}
\left . 
\begin{array}{rcl}
\displaystyle{\frac{d}{dt}\phi_i+A_i\phi_i}&=&0\\
\phi_i(t_n)&=&\phi^{n+\frac{i-1}{2}}
\end{array}
\right\}\quad {\rm in}\ (t_n,t_{n+1}),
\end{equation*}
and then set
$
\phi^{n+\frac i 2}=\phi_i(t_{n+1}), \ {\rm for}\ i = 1,2.
$
It can be shown that this method is first-order accurate in time, see e.g.,  \cite{glowinski2003finite}.

We apply this approach to split problem \eqref{FSIeqRef}-\eqref{IC_ALE} 
 into two sub-problems: a structure and a fluid sub-problem defining operators $A_1$ and $A_2$. 

\vskip 0.1in
{\bf Problem A1: The structure elastodynamics problem.} In this step we solve an elastodynamics problem for the location of the 
multi-layered cylinder wall. The problem is driven only by the initial data, i.e., the initial boundary velocity, 
taken from the previous time step as the trace of the fluid velocity at the fluid-structure interface. 
%Thus, the structure sub-problem communicates with the fluid sub-problem only through the initial data,
%i.e., through the kinematic coupling condition. 
The fluid velocity $\bf u$ remains unchanged in this step.
More precisely, the problem reads: 
\\
{\sl Given $({\bf u}^n,\eta^n,v^n,\bd^n,\bV^n)$ from the previous time step, find 
$({\bf u},v,\eta,\bV,\bd)$ such that:
\begin{equation}\label{Step1DF}
  \begin{array}{ll} 
\quad {\partial_t\mathbf{u}} = 0,   &\quad \textrm{in} \; (t_n, t_{n+1})\times\Omega_F,  \\  \\
 \left.
 \begin{array}{rcl}
  \partial_t{\bf V} &=& \nabla\cdot{\bf S},\; \\
   \partial_t{\bf d}&=&{\bf V} 
   \end{array}\right\} 
                                     &\quad \textrm{in} \; (t_n, t_{n+1})\times\Omega_S, \\  
\quad  {\bf d}=0  &\quad {\rm on}\; \Gamma_{in/out}^s,  \\ 
\quad  {\bf S}{\bf e}_r=0  &\quad {\rm on}\; (t_n,t_{n+1})\times \Gamma_{ext}, 
\end{array}
\end{equation}
\begin{equation}
  \begin{array}{ll} 
\quad  \quad  \quad \quad \  \; {\bf d}=\eta{\bf e}_r  &\quad \textrm{on} \; (t_n, t_{n+1})\times (0,1),\\ 
  \left.
  \begin{array}{rcl}
     \partial_t v-\partial^2_z\eta&=&0 ,\\ 
                          \partial_t{\eta}  &=& v
  \end{array}
   \right\}
                    &\quad \textrm{on} \; (t_n, t_{n+1})\times (0,1),  \\ 
\quad \eta(0)=\eta(1)=0, & 
\end{array}
\end{equation}
with ${\bf u}(t_n)={\bf u}^{n},\; \eta(t_n)=\eta^n,\; v(t_n)=v^n,\; {\bf d}(t_n)={\bf d}^n,\; {\bf V}(t_n)={\bf V}^n.$\\
Then set ${\bf u}^{n+\frac 1 2}={\bf u}(t_{n+1})$, $\eta^{n+\frac 1 2}=\eta(t_{n+1})$, $v^{n+\frac 1 2}=v(t_{n+1})$,
${\bf d}^{n+\frac 1 2}={\bf d}(t_{n+1})$, ${\bf V}^{n+\frac 1 2}={\bf V}(t_{n+1})$.}

\vskip 0.1in
{\bf Problem A2: The fluid problem.} In this step we solve the Navier-Stokes equations coupled with structure inertia 
through a ``Robin-type'' boundary condition  on $\Gamma$ (lines 5 and 6 in \eqref{A2} below).
The kinematic coupling condition is implicitly satisfied.
The structure displacement remains unchanged.
With a slight abuse of notation, the problem can be written as follows:
\\
{\sl Find 
$({\bf u},v,\eta,\bV,\bd)$ such that:
\begin{eqnarray}
\nonumber
 \partial_t\eta  = 0 \hskip 0.85in  &\textrm{on} \; (t_n, t_{n+1})\times(0,1),\\
 \nonumber
 \partial_t{\bf d}=0 \hskip 0.85in  &\textrm{on} \; (t_n, t_{n+1})\times\Omega_S,\\
 \nonumber
\left.
\begin{array}{rcl}
 \partial_t \mathbf{u}+(({\bf u}^n-{\bf w}^{\eta^{n+\frac 1 2}})\cdot\nabla^{\eta^n}){\bf u} &=&\nabla^{\eta^n} \cdot{\sigma}^{\eta^n}   \\
 \nabla^{\eta^n} \cdot \mathbf{u}&=&0
 \end{array}
\right\} 
&\textrm{in} \; (t_n, t_{n+1})\times\Omega_F, 
\\
\label{A2}
\left.
\begin{array}{rcl}
 \partial_t v&=&- J\sigma{\bf n}\cdot {\bf e}_r\\
 {\bf u}&=&v{\bf e}_r
 \end{array}
\right\} 
&\textrm{on} \; (t_n, t_{n+1})\times (0,1), \\
\nonumber
 \left.
 \begin{array}{rcl}
 u_r&=&0\\
 {\partial_r u_z}&=&0
 \end{array}
 \right\} 
 & \textrm{on}\; (t_n, t_{n+1})\times\Gamma_b,\\ 
 \nonumber
%  \end{eqnarray}
 %\begin{eqnarray*}
\left.
 \begin{array}{rcl}
 p+\frac{\rho_f}{2}|u|^2&=&P_{in/out}(t)\\
 u_r&=&0
 \end{array}
\right\}
& \textrm{on}\; (t_n, t_{n+1})\times\Gamma_{in/out},
\end{eqnarray}
\begin{equation*}
 {\rm with}\ 
{\bf u}(t_n,.)={\bf u}^{n+\frac 1 2},\; \eta(t_n,.)=\eta^{n+\frac 1 2}, \;  v(t_n,.)=v^{n+\frac 1 2},
\; {\bf d}(t_n,.)={\bf d}^{n+\frac 1 2}, \;  {\bf V}(t_n,.)={\bf V}^{n+\frac 1 2}. 
\end{equation*}
Then set $\mathbf{u}^{n+1}=\mathbf{u}(t_{n+1}), \; {\eta}^{n+1}={\eta}(t_{n+1}),\; v^{n+1}=v(t_{n+1}),\;
{\bf d}^{n+1}={\eta}(t_{n+1}),\; {\bf V}^{n+1}={\bf V}(t_{n+1}).$
}
\vskip 0.1in
Notice that, since  in this step $\eta$ does not change,  this problem is linear. 
Furthermore, it can be viewed as a stationary Navier-Stokes-like
problem on a fixed domain with a Robin-type boundary condition. 
In numerical simulations, one can use the ALE transformation $A_{\eta^n}$ to 
``transform'' the problem back to domain $\Omega_{\eta^n}$
and solve it there, thereby avoiding the un-necessary calculation
of the transformed gradient $\nabla^{\eta^n}$. The ALE velocity is the only extra term
that needs to be included with that approach. See, e.g., \cite{MarSun} for
more details. For the purposes of our proof, we will, however, remain in the fixed, reference domain $\Omega_F$.

It is important to notice that in Problem A2, the problem is ``linearized'' around the previous location of the boundary, i.e., 
we work with the domain determined by $\eta^n$,
and not by $\eta^{n+1/2}$. This is in direct relation with the implementation of the numerical scheme studied in \cite{MarSun,SunBorMar}.
However, we also notice that  ALE velocity, $w^{n+\frac 1 2}$, is taken from the just calculated Problem A1! 
This choice is {\sl crucial} for obtaining a semi-discrete version of an energy inequality, discussed in Section~\ref{sec:approximate_solutions}.

In the remainder of this paper we use the splitting scheme described above to define approximate solutions of 
\eqref{FSIeqRef}-\eqref{IC_ALE} (or equivalently of problem \eqref{NS}-\eqref{CC} ) and
show that the approximate solutions converge to a weak solution,  as $\Delta t\rightarrow 0$.

\section{Weak solutions}\label{sec:weak_solutions}

\subsection{Notation and function spaces}
To define weak solutions of the moving-bounday problem \eqref{NS}-\eqref{CC} and of the moving-boundary problem \eqref{FSIeqRef}-\eqref{IC_ALE}
defined on a fixed domain, we introduce the following notation.
We use $a_S$ to denote the bilinear form
associated with the elastic energy of the thick structure:
\begin{equation}\label{Elastic}
a_S({\bf d},\bpsi)=\int_{\Omega_S}\big (2{\bf D}({\bf d}):{\bf D}({\bpsi})+(\nabla\cdot{\bf d})\cdot(\nabla\cdot{\bpsi})\big ).
\end{equation}
Here $A : B :={\rm tr} \left[A B^T\right]$.
Furthermore, we will be using $b$ to denote the following trilinear form  corresponding to the 
(symmetrized) nonlinear advection term in the Navier-Stokes equations:
\begin{equation}\label{transport}
b(t,{\bf u},{\bf v},{\bf w})=\frac 1 2\int_{\Omega_{F}(t)}({\bf u}\cdot\nabla){\bf v}\cdot{\bf w}-
\frac 1 2\int_{\Omega_{F}(t)}({\bf u}\cdot\nabla){\bf w}\cdot{\bf v}.
\end{equation}
Finally, we define a linear functional which associates the inlet and outlet dynamic pressure boundary data to
a test function $\bf v$ in the following way:
$$
\langle F(t),{\bf v}\rangle_{\Gamma_{in/out}}=P_{in}(t)\int_{\Gamma_{in}}v_z-P_{out}(t)\int_{\Gamma_{out}}v_z.
$$

The following functions spaces define our weak solutions.
For the fluid velocity we would like to work with the classical function space associated with weak solutions of the Navier-Stokes equations. 
This, however, requires some additional consideration. Namely, since our thin structure is governed by the linear wave equation, 
lacking the bending rigidity terms, weak solutions cannot be expected to be Lipschitz-continuous. 
Indeed,  from the energy inequality \eqref{EE} we only have $\eta\in H^1(0,1)$, and from  
Sobolev embedding we get that $\eta\in C^{0,1/2}(0,1)$, which means that $\Omega_F(t)$ is not necessarily a Lipshitz domain.
However, $\Omega_F(t)$ is locally a sub-graph of a H\" older continuous function.
In that case one can define a``Lagrangian" trace
\begin{equation}\label{trace}
\begin{split}
\gamma_{\Gamma(t)} &:C^1(\overline{\Omega_F(t)}) \to  C(\Gamma),\\
\gamma_{\Gamma(t)} &:v\mapsto v(t,z,r+\eta(t,z)).
\end{split}
\end{equation} 
Furthermore, it was shown in \cite{CDEM,CG,BorisTrag} that the trace operator $\gamma_{\Gamma(t)}$ can be extended 
by continuity to a linear operator from $H^1(\Omega_F(t))$ to $H^s(\Gamma)$, $0\leq s<\frac 1 4$. For a precise statement of the results about
``Lagrangian'' trace, we refer the reader to 
 Theorem \ref{TraceTm} below \cite{BorisTrag}. Now, we define the velocity solution space in the following way:
\begin{equation}\label{V_Fglatke}
\begin{array}{rcl}
\displaystyle{ V_F(t)}&=&\displaystyle{\{{\mathbf u}=(u_z,u_r)\in C^1(\overline{\Omega_{F}(t)})^2:\nabla\cdot{\bf u}=0,} \\ 
 && \displaystyle{u_z=0\ {\rm on}\ \Gamma(t),
\ u_r=0\ {\rm on}\ \partial\Omega_{F}(t)\setminus\Gamma(t)\},} \\
\displaystyle{{\cal V}_F(t)}&=&\displaystyle{\overline{V_F(t)}^{H^1(\Omega_{F}(t))}.}
\end{array}
\end{equation}
Using the fact that $\Omega_{F}(t)$ is locally a sub-graph of a H\"older continuous function we can get the following characterization of 
the velocity solution space ${\cal V}_F(t)$
(see \cite{CDEM,CG}):
\begin{equation}\label{V_F}
\begin{array}{rcl}
\displaystyle{{\cal V}_F(t)}&=&\displaystyle{\{{\mathbf u}=(u_z,u_r)\in H^1({\Omega_{\eta}(t)})^2:\nabla\cdot{\bf u}=0,} \\ 
 && \displaystyle{u_z=0\ {\rm on}\ \Gamma(t),
\ u_r=0\ {\rm on}\ \partial\Omega_{\eta}(t)\setminus\Gamma(t)\}.} \\
\end{array}
\end{equation}
The function space associated with weak solutions of the 1D linear wave equation and the thick wall are given, respectively, by
\begin{equation}\label{V_K}
{\cal V}_W=H^1_0(0,1), 
\end{equation}
\begin{equation}\label{V_S}
{\cal V}_S=\{\bd=(d_z,d_r)\in H^1(\Omega_S)^2:d_z=0\ {\rm on}\ \Gamma,\ \bd=0\ {\rm on}\ \Gamma_{in/out}^s\}.
\end{equation}
Motivated by the energy inequality we also define the corresponding evolution spaces
for the fluid and structure sub-problems, respectively:
\begin{equation}\label{vel_test}
{\cal W}_F(0,T)=L^{\infty}(0,T;L^2(\Omega_{F}(t)))\cap L^2(0,T;{\cal V}_F(t)),
\end{equation}
\begin{equation}\label{struc_test}
{\cal W}_W(0,T)=W^{1,\infty}(0,T;L^2(0,1))\cap L^2(0,T;{\cal V}_W),
\end{equation}
\begin{equation}\label{thick_test}
{\cal W}_S(0,T)=W^{1,\infty}(0,T;L^2(\Omega_S))\cap L^2(0,T;{\cal V}_S).
\end{equation}
Finally, we are in a position to define the solution space for the coupled fluid-multi-layered-structure interaction problem.
This space must involve
the kinematic coupling condition. The dynamic coupling condition will be enforced in a weak sense, through the integration by parts 
in the weak formulation of the problem. Thus, we define
\begin{equation}\label{W}
\begin{array}{c}
{\cal W}(0,T)=\{({\bf u},\eta,\bd)\in {\cal W}_F(0,T)\times{\cal W}_W(0,T)\times{\cal W}_S(0,T):
\\ 
{\bf u}(t,z,1+\eta(t,z))=\partial_t\eta(t,z){\bf e}_r,\; \bd(t,z,1)=\eta(t,z){\bf e}_r\}.
\end{array}
\end{equation}
Equality ${\bf u}(t,z,1+\eta(t,z))=\partial_t\eta(t,z){\bf e}_r$ is taken in the sense defined in \cite{CDEM,BorisTrag}.
The corresponding test space will be denoted by
\begin{equation}\label{Q}
\begin{array}{c}
{\cal Q}(0,T)=\{({\bf q},\psi,\bpsi)\in C^1_c([0,T);{\cal V}_F\times{\cal V}_W\times{\cal V}_S):{\bf q}(t,z,1+\eta(t,z))=\psi(t,z){\bf e}_r=\bpsi(t,z,1)\}.
\end{array}
\end{equation}

\subsection{Weak solutions for the problem defined on the moving domain}
We are now in a position to define weak solutions of fluid-multi-layered structure interaction problem,
defined on the moving domain $\Omega_F(t)$.
\begin{definition}\label{DefWS}
We say that $({\bf u},\eta,\bd)\in{\cal W}(0,T)$ is a weak solution of problem \eqref{NS}-\eqref{CC}  if
for every $({\bf q},\psi,\bpsi)\in{\cal Q}(0,T)$  the  following equality holds:
\begin{equation}
\begin{array}{c}
\displaystyle{-\int_0^T\int_{\Omega_{F}(t)}{\bf u}\cdot\partial_t{\bf q}+\int_0^T b(t,{\bf u},{\bf u},{\bf q})+
2\int_0^T\int_{\Omega_{F}(t)}{\bf D}({\bf u}):{\bf D}({\bf q})-\frac{1}{2}\int_0^T\int_0^1(\partial_t\eta)^2\psi}
\\ \\
\displaystyle{
-\int_0^T\int_0^1 \partial_t\eta\partial_t\psi+\int_0^T\int_0^1\partial_z\eta\partial_z\psi
-\int_0^T\int_{\Omega_S}\bd\cdot\partial_t\bpsi
+\int_0^Ta_S(\bd,\bpsi)}
\\ \\
\displaystyle{=\int_0^T\langle F(t),{\bf q}\rangle_{\Gamma_{in/out}}+\int_{\Omega_{\eta_0}}{\bf u}_0\cdot{\bf q}(0)+\int_0^1v_0\psi(0)+\int_{\Omega_S}{\bf V}_0\cdot\bpsi(0).}
\end{array}
\label{VF}
\end{equation}
\end{definition}

In deriving the weak formulation we used integration by parts in a classical way, and the following equalities which hold for smooth functions:
\begin{eqnarray*}
\int_{\Omega_{F}(t)}({\bf u}\cdot\nabla){\bf u}\cdot{\bf q}=&
\displaystyle{\frac 1 2\int_{\Omega_{F}(t)}({\bf u}\cdot\nabla){\bf u}\cdot{\bf q}-
\frac 1 2\int_{\Omega_{F}(t)}({\bf u}\cdot\nabla){\bf q}\cdot{\bf u}}\\
&\displaystyle{+\frac 1 2\int_0^1(\partial_t\eta)^2\psi\pm\frac 1 2\int_{\Gamma_{out/in}}|u_r|^2v_r},
\end{eqnarray*}
$$
\int_0^T\int_{\Omega_{F}(t)}\partial_t{\bf u}\cdot{\bf q}=-\int_0^T\int_{\Omega_{F}(t)}{\bf u}\cdot\partial_t{\bf q}-\int_{\Omega_{\eta_0}}{\bf u}_0\cdot{\bf q}(0)
-\int_0^T\int_0^1(\partial_t\eta)^2\psi.
$$

\subsection{Weak solutions for the problem defined on a fixed, reference domain}
Since most of the analysis will be performed on the problem mapped to $\Omega_F$, 
we rewrite the above definition in terms of  $\Omega_F$ using the ALE mapping
$A_\eta(t)$ defined in \eqref{RefTrans}.
For this purpose, the following notation will be useful.
We define the transformed trilinear functional $b^\eta$:
\begin{equation}\label{b_eta}
b^{\eta}({\bf u},{\bf u},{\bf q}):=\frac 1 2\int_{\Omega_F}(1+\eta)(({\bf u}-{\bf w}^{\eta})\cdot\nabla^{\eta}){\bf u}\cdot{\bf q}
- \frac 1 2\int_{\Omega_F} (1+\eta) (({\bf u}-{\bf w}^{\eta})\cdot\nabla^{\eta}){\bf q}\cdot{\bf u},
\end{equation}
where $1+\eta$ is the Jacobian of the ALE mapping, calculated in \eqref{ALE_Jacobian}. 
Notice that we have included the ALE domain velocity ${\bf w}^{\eta}$ into $b^{\eta}$.

It is important to point out that
the transformed fluid velocity ${\bf u}^\eta$ is not divergence-free anymore.
Rather,  it satisfies the transformed divergence-free condition $\nabla^{\eta}\cdot{\bf u}^\eta=0$.
Furthermore, since $\eta$ is not a Lipschitz function, the ALE mapping is not necessarily a Lipschitz function,
and, as a result, ${\bf u}^{\eta}$ is not necessarily an $H^1$ function on $\Omega_F$.
Therefore we need to redefine the function spaces for the fluid velocity by introducing
$$
{\cal V}_F^{\eta}=\{{\bf u}^{\eta}:{\bf u}\in {\cal V}_F(t)\},
$$
where ${\bf u}^{\eta}$ is defined in \eqref{u_eta}. Under the assumption $1+\eta(z)>0$, $z\in [0,1]$, 
we can define a scalar product on ${\cal V}_F^{\eta}$ 
in the following way:
$$
({\bf u}^{\eta},{\bf v}^{\eta})_{{\cal V}_F^{\eta}}=\int_{\Omega_F}(1+\eta)\big ({\bf u}^{\eta}\cdot{\bf v}^{\eta}+\nabla^{\eta}{\bf u}^{\eta}:\nabla^{\eta}{\bf v}^{\eta}\big)
=({\bf u},{\bf v})_{H^1(\Omega_{F}(t))}.
$$
Therefore, ${\bf u}\mapsto {\bf u}^{\eta}$ is an isometric isomorphism between ${\cal V}_F(t)$ and ${\cal V}_F^{\eta}$, so ${\cal V}_F^{\eta}$ is also a Hilbert space.
The function spaces ${\cal W}_F^{\eta}(0,T)$ and ${\cal W}^{\eta}(0,T)$ are defined the same as before, but with ${\cal V}_F^{\eta}$ instead ${\cal V}_F(t)$.
More precisely:
\begin{equation}\label{W_F_eta}
{\cal W}_F^\eta(0,T)=L^{\infty}(0,T;L^2(\Omega_F))\cap L^2(0,T;{\cal V}_F^\eta(t)),
\end{equation}
\begin{equation}\label{W_eta}
{\cal W^\eta}(0,T)=\{({\bf u},\eta,\bd)\in {\cal W}_F^\eta(0,T)\times{\cal W}_W(0,T)\times{\cal W}_S(0,T):{\bf u}(t,z,1)=\partial_t\eta(t,z){\bf e}_r,\; \eta(t,z)=\bd(t,z,1)\}.
\end{equation}
The corresponding test space is defined by
\begin{equation}\label{Q_eta}
{\cal Q^\eta}(0,T)=\{({\bf q},\psi,\bd)\in C^1_c([0,T);{\cal V}_F^\eta\times{\cal V}_W\times{\cal V}_S):{\bf q}(t,z,1)=\psi(t,z){\bf e}_r=\bd(t,z,1)\}.
\end{equation}

\begin{definition}\label{DefWSRef}
We say that $({\bf u},\eta,\bd)\in{\cal W}^{\eta}(0,T)$ is a weak solution of problem \eqref{FSIeqRef}-\eqref{IC_ALE} defined on the reference domain $\Omega_F$,
if for every $({\bf q},\psi,\bpsi)\in {\cal Q^\eta}(0,T)$ the following equality holds:
\begin{equation}
\begin{array}{c}
\displaystyle{-\int_0^T\int_{\Omega_F}(1+\eta){\bf u}\cdot\partial_t{\bf q}+\int_0^T b^{\eta}({\bf u},{\bf u},{\bf q})}
\displaystyle{+2\int_0^T\int_{\Omega_F}(1+\eta){\bf D}^{\eta}({\bf u}):{\bf D}^{\eta}({\bf q})}
\\ \\
\displaystyle{-\frac{1}{2}\int_0^T\int_{\Omega_F}(\partial_t\eta){\bf u}\cdot{\bf q}}
\displaystyle{-\int_0^T\int_0^1\partial_t\eta\partial_t\psi+\int_0^T\int_0^1\partial_z\eta\partial_z\psi}
\\ \\
\displaystyle{-\int_0^T\int_{\Omega_S}\bd\cdot\partial_t\bpsi
+\int_0^Ta_S(\bd,\bpsi)}
\\ \\
\displaystyle{=\int_0^T\langle F(t),{\bf q}\rangle_{\Gamma_{in/out}}+\int_{\Omega_{\eta_0}}{\bf u}_0\cdot{\bf q}(0)+\int_0^1v_0\psi(0)+\int_{\Omega_S}{\bf V}_0\cdot\bpsi(0).}
\end{array}
\label{VFRef}
\end{equation}
\end{definition}

To see that this is consistent with the weak solution defined in Definition~\ref{DefWS}, we present the main steps in the 
transformation of the first integral on the left hand-side in \eqref{VF}, responsible for the fluid kinetic energy. 
Namely, we formally calculate:
$$
-\int_{\Omega_{F}(t)}{\bf u}\cdot\partial_t{\bf q}=
-\int_{\Omega_F}(1+\eta){\bf u}^{\eta}\cdot(\partial_t{\bf q}-({\bf w}^{\eta}\cdot{\nabla}^{\eta}){\bf q})
=-\int_{\Omega_F}(1+\eta){\bf u}^{\eta}\cdot\partial_t{\bf q}
$$
$$
+ \frac 1 2\int_{\Omega_F}(1+\eta)({\bf w}^{\eta}\cdot{\nabla}^{\eta}){\bf q}\cdot{\bf u}^{\eta}
+ \frac 1 2\int_{\Omega_F}(1+\eta)({\bf w}^{\eta}\cdot{\nabla}^{\eta}){\bf q}\cdot{\bf u}^{\eta}.$$
In the last integral on the right hand-side we use the definition of ${\bf w}^\eta$ and of $\nabla^\eta$,
given in \eqref{nablaeta}, to obtain
$$
\int_{\Omega_F}(1+\eta)({\bf w}^{\eta}\cdot{\nabla}^{\eta}){\bf q}\cdot{\bf u}^{\eta}=\int_{\Omega_F}\partial_t\eta \ {\tilde r}\ \partial_{\tilde r}{\bf q}\cdot{\bf u}^{\eta}.
$$
Using integration by parts with respect to $r$, keeping in mind that $\eta$ does not depend on $r$,
we obtain
$$
-\int_{\Omega_{F}(t)}{\bf u}\cdot\partial_t{\bf q}=
-\int_{\Omega_F}(1+\eta){\bf u}^{\eta}\cdot(\partial_t{\bf q}-({\bf w}^{\eta}\cdot{\nabla}^{\eta}){\bf q})
=-\int_{\Omega_F}(1+\eta){\bf u}^{\eta}\cdot\partial_t{\bf q}
$$
$$
+ \frac 1 2\int_{\Omega_F}(1+\eta)({\bf w}^{\eta}\cdot{\nabla}^{\eta}){\bf q}\cdot{\bf u}^{\eta}-
\frac 1 2\int_{\Omega_F}(1+\eta)({\bf w}^{\eta}\cdot{\nabla}^{\eta}){\bf u}^{\eta}\cdot{\bf q}-\frac 1 2\int_{\Omega_F}\partial_t\eta{\bf u}^{\eta}\cdot{\bf q}
+\frac 1 2\int_0^1(\partial_t\eta)^2\psi,
$$
By using this identity in \eqref{VF}, and by recalling the definitions for $b$ and $b^\eta$, we obtain exactly the weak form \eqref{VFRef}.

In the remainder of this manuscript we will be working on the fluid-multi-layered structure interaction problem defined on the fixed domain $\Omega_F$,
satisfying the weak formulation presented in Definition~\ref{DefWSRef}
For brevity of notation, since no confusion is possible, we omit the superscript ``tilde'' which is used to denote the coordinates of points in $\Omega_F$.

\section{Approximate solutions}\label{sec:approximate_solutions}
In this section we use the Lie operator splitting scheme and semi-discretization to define a sequence of approximate solutions of
the FSI problem  \eqref{FSIeqRef}-\eqref{IC_ALE}. 
Each of the sub-problems defined by the Lie splitting in Section~\ref{sec:LieSplitting} as Problem A1 and Problem A2, 
will be discretized in time using the Backward Euler scheme. 
This approach defines a time step, which will be denoted by $\Delta t$, and a number of time sub-intervals $N\in\N$, so that
$$
(0,T) = \cup_{n=0}^{N-1} (t^n,t^{n+1}), \quad t^n = n\Delta t, \ n=0,...,N-1.
$$
For every subdivision containing $N\in\N$ sub-intervals,  the vector of unknown approximate solutions will be denoted by
%\begin{equation}
%{\bf X}_N^{n+\frac i 2}=\left (\begin{array}{c} {\bf u}_N^{n+\frac i 2} \\ v_N^{n+\frac i 2} \\ \eta_N^{n+\frac i 2}
%\\ \bV_N^{n+\frac i 2} \\ \bd_N^{n+\frac i 2} \end{array}\right ), n=0,1,\dots,N-1,\,\ i=1,2,
%\label{X}
%\end{equation}
\begin{equation}
{\bf X}_N^{n+\frac i 2}=\left ({\bf u}_N^{n+\frac i 2}, v_N^{n+\frac i 2},  \eta_N^{n+\frac i 2}, \bV_N^{n+\frac i 2}, \; \bd_N^{n+\frac i 2} \right )^T, n=0,1,\dots,N-1,\,\ i=1,2,
\label{X}
\end{equation}
where $i = 1,2$ denotes the solution of Problem A1 or A2, respectively.
The initial condition will be denoted by
$
{\bf X}^0=\left ( {\bf u}_0, v_0, \eta_0,  \bV_0, \bd_0  \right )^T.
$

The semi-discretization and the splitting of the problem will be performed in 
such a way that the semi-discrete version of the energy inequality \eqref{EE} is preserved 
at every time step.
This is a crucial ingredient for the existence proof.

The semi-discrete versions of the kinetic and elastic energy  \eqref{energy}, and of dissipation  \eqref{dissipation} are defined  by
the following:
 \begin{equation}
\begin{array}{c}
\displaystyle{E_{kin,N}^{n+\frac i 2}=\frac 1 2\Big (\int_{\Omega_F}(1+\eta^{n-1+i})|{\bf u}^{n+\frac i 2}_N|^2+
\|v^{n+\frac i 2}_N\|^2_{L^2(0,1)}+\|\bV^{n+\frac i 2}_N\|^2_{L^2(\Omega_S)}\Big ),}\\
\displaystyle{E_{el,N}^{n+1}=\frac 1 2\Big (
\|\partial_z\eta^{n+\frac 1 2}_N\|^2_{L^2(0,1)}+2\|{\bf D}(\bd^{n+\frac 1 2}_N)\|^2_{L^2(\Omega_S)}+\|\nabla\cdot\bd^{n+\frac 1 2}_N\|^2_{L^2(\Omega_S)}
\Big ),}
\\
\displaystyle{E_N^{n+\frac i 2}=E_{kin,N}^{n+\frac i 2}+E_{el,N}^{n+1},}
\end{array}
\label{kenergija}
\end{equation}
\begin{equation}
\begin{array}{c}
\displaystyle{D_N^{n+1}=\Delta t\int_{\Omega_F}(1+\eta^n)|D^{\eta^n}({\bf u}_N^{n+1})|^2,\; n=0,\dots,N-1,\; i=0,1.}
\end{array}
\label{kdisipacija}
\end{equation}
Throughout the rest of this section we fix the time step $\Delta t$, i.e., we keep $N\in\N$ fixed, and study 
the semi-discretized sub-problems defined by the Lie splitting.
To simplify notation, we  omit the subscript $N$ and write
$({\bf u}^{n+\frac i 2},v^{n+\frac i 2},\eta^{n+\frac i 2},\bV^{n+\frac i 2},\bd^{n+\frac i 2})$ instead of 
$({\bf u}^{n+\frac i 2}_N,v^{n+\frac i 2}_N,\eta^{n+\frac i 2}_N,\bV^{n+\frac i 2}_N,\bd^{n+\frac i 2}_N)$.

%We will show bellow that each sub-problem satisfies a discrete version of the energy inequality \eqref{EE},
%involving the discrete energy and dissipation, defined in \eqref{kenergija} and \eqref{kdisipacija}.
%This will be used later in the proof of convergence of approximate solutions to a weak solution of Problem~\ref{FSIref}.

\subsection{Semi-discretization of Problem A1}
A semi-discrete version of Problem A1 (Structure Elastodynamics), defined by the Lie splitting in \eqref{Step1DF} can be written as follows.
First, in this step ${\bf u}$ does not change, and so
$${\bf u}^{n+\frac 1 2}={\bf u}^n.$$
We define  $(v^{n+\frac 1 2},\eta^{n+ \frac 1 2},{\bf V}^{n+\frac 1 2},{\bf U}^{n+\frac 1 2})\in{\cal V}_W^2\times{\cal V}_S^2$ as a solution of the following problem, written in weak form:
\begin{equation}\label{SProb1}
\begin{array}{c}
{\bd}^{n+\frac 1 2}(z,1)=\eta^{n+\frac 1 2}(z,1){\bf e}_r,\ z\in (0,1),\\
\displaystyle{\frac{{\bd}^{n+\frac 1 2}-{\bd}^{n}}{\Delta t}={\bV}^{n+\frac 1 2}},\ 
\displaystyle{\frac{\eta^{n+\frac 1 2}-\eta^{n}}{\Delta t}=v^{n+\frac 1 2}}, \\
\displaystyle{\int_{\Omega_S}\frac{{\bf V}^{n+\frac 1 2}-{\bf V}^{n}}{\Delta t}}\cdot\boldsymbol\Psi+
\displaystyle{\int_0^1\frac{v^{n+\frac 1 2}-v^{n}}{\Delta t}}\psi+a_S({\bd}^{n+\frac 1 2},\boldsymbol\Psi)+\int_0^1\partial_z \eta^{n+\frac 1 2}\partial_z\psi=0,  
\end{array}
\end{equation}
which holds for all $(\psi,\boldsymbol\Psi)\in {\cal V}_W\times{\cal V}_S$ such that $\boldsymbol\Psi(t,z,1)=\psi(t,z)$.
The first equation enforces the kinematic coupling condition, the second row in \eqref{SProb1} introduces 
the structure velocities, while the third equation
corresponds to a weak form of the semi-discretized elastodynamics problem. Notice that we solve
the  thin and thick structure problems as one coupled problem. The thin structure enters as a boundary condition for the thick structure problem.
\begin{proposition}
For each fixed $\Delta t > 0$, problem \eqref{SProb1} has a unique solution $(v^{n+\frac 1 2},\eta^{n+\frac 1 2},\bV^{n+\frac 1 2},\bd^{n+\frac 1 2})\in {\cal V}_W^2\times{\cal V}_S^2$.
\end{proposition}

\proof 
First notice that Korn's inequality implies that the bilinear form $a_S$ is coercive on ${\cal V}_S$. From here,
the proof is a direct consequence of the Lax-Milgram Lemma applied to the weak form
\begin{equation}\nonumber
\begin{array}{c}
\displaystyle{\int_0^1\eta^{n+\frac 1 2}\psi+\int_{\Omega_S}{\bd}^{n+1}\cdot\boldsymbol\Psi+
(\Delta t)^2\big (\int_0^1\partial_z\eta\partial_z\psi+a_S({\bd}^{n+\frac 1 2},\boldsymbol\Psi)\big )}\\ 
\displaystyle{=\int_0^L\big (\Delta tv^n+\eta^n\big )\psi+\int_{\Omega_S}(\Delta t{\bV}^n+{\bd}^n)\cdot\boldsymbol\Psi,\ 
\forall (\psi,\boldsymbol\Psi)\in \{{\cal V}_W\times{\cal V}_S| \boldsymbol\Psi(t,z,1)=\psi(z,1)\},}
\end{array}
\end{equation}
which is obtained after a substitution of $v^{n+\frac 1 2}$ and $\bV^{n+\frac 1 2}$ in the third equation in \eqref{SProb1},
by using the equations (\ref{SProb1})$_2$.
\qed

\begin{proposition}\label{prop:Energy1}
For each fixed $\Delta t > 0$, a solution of problem \eqref{SProb1} satisfies the following discrete energy equality:
\begin{equation}
\begin{array}{c}
\displaystyle{E_{kin,N}^{n+\frac 1 2}+E_{el,N}^{n+1}+\frac 1 2\big (\|v^{n+\frac 1 2}-v^{n}\|^2_{L^2(0,1)}+
\|\bV^{n+\frac 1 2}-\bV^n\|^2_{L^2(\Omega_S)}}\\
[0.3cm]
\displaystyle{+\|\partial_z(\eta^{n+\frac 1 2}-\eta^{n})\|^2_{L^2(0,1)}+a_S(\bd^{n+\frac 1 2}-\bd^{n},\bd^{n+\frac 1 2}-\bd^{n}) \big)=E_{kin,N}^n+E_{el,N}^n,}
\end{array}
\label{DEO31}
\end{equation}
where the kinetic and elastic energy, $E_{kin,N}^n$, $E_{el,N}^n$, are  defined in \eqref{kenergija}.
\end{proposition}

\proof
From the second  row in (\ref{SProb1}) we immediately get 
$$v^{n+\frac 1 2}=\frac{\eta^{n+\frac 1 2}-\eta^{n}}{\Delta t}\in {\cal V}_W,\;
\bV^{n+\frac 1 2}=\frac{\bd^{n+\frac 1 2}-\bd^{n}}{\Delta t}\in {\cal V}_S.
$$
Therefore, we can proceed as usual, by substituting the test functions in  \eqref{SProb1} with structure velocities.
More precisely, we replace the test function $(\psi,\bpsi)$ by $(v^{n+\frac 1 2},\bV^{n+\frac 1 2})$ in the first term on the left hand-side, 
and then replace $(\psi,\bpsi)$ by $(({\eta^{n+\frac 1 2}-\eta^{n}})/{\Delta t},({\bd^{n+\frac 1 2}-\bd^{n}})/{\Delta t})$ in the bilinear forms 
that correspond to the elastic energy. 
To deal with the terms 
$(v^{n+1/2}-v^n)v^{n+1/2}$, $(\eta^{n+1/2}-\eta^n)\eta^{n+1/2}$, $(\bV^{n+1/2}-\bV^n)\cdot\bV^{n+1/2}$, and $(\bd^{n+1/2}-\bd^n)\cdot\bd^{n+1/2}$,
we use the algebraic identity
$(a-b)\cdot a=\frac 1 2(|a|^2+|a-b|^2-|b|^2)$.
After multiplying the entire equation by $\Delta t$, 
the third equation in  \eqref{SProb1} can be written as:
$$
(\|v^{n+\frac 1 2}\|^2_{L^2(0,1)}+\|v^{n+\frac 1 2}-v^{n}\|^2_{L^2(0,1)})+
(\|\bV^{n+\frac 1 2}\|^2_{L^2(\Omega_S)}+\|\bV^{n+\frac 1 2}-\bV^{n}\|^2_{L^2(\Omega_S)})
$$
$$
\|\partial_z\eta^{n+\frac 1 2}\|^2_{L^2(0,1)}+
\|\partial_z(\eta^{n+\frac 1 2}-\partial_z\eta^{n})\|^2_{L^2(0,1)}
+a_S(\bd^{n+\frac 1  2},\bd^{n+\frac 1  2})
$$
$$
+a_S(\bd^{n+\frac 1  2}-\bd^n,\bd^{n+\frac 1  2}-\bd^n)=\|v^{n}\|^2_{L^2(0,1)}+\|\bV^{n}\|^2_{L^2(\Omega_S)}+
\|\partial_z\eta^{n}\|^2_{L^2(0,1)}+a_S(\bd^n,\bd^n).
$$
Since in this sub-problem ${\bf u}^{n+\frac 1 2}={\bf u}^n$, we can add 
$\rho_f \int_{\Omega_F} (1+\eta^n){\bf u}^{n+1/2}$ on the left hand-side, and
$\rho_f \int_{\Omega _F}(1+\eta^n){\bf u}^{n}$ on the right hand-side of the equation. Furthermore,
displacements ${\bf d}^{n+\frac 1 2}$ and $\eta^{n+\frac 1 2}$ do not change in Problem A2 (see (\ref{A2dis})), and
so we can replace $\bd^n$ and $\eta^n$ on the right hand-side of the equation with $\bd^{n-\frac 1 2}$ and $\eta^{n-\frac 1 2}$,
respectively,
to obtain exactly the energy equality \eqref{DEO31}.
\qed

\subsection{Semi-discretization of Problem A2}
In this step $\eta$, $\bd$ and $\bV$ do not change, and so 
\begin{equation}\label{A2dis}
\eta^{n+1}=\eta^{n+\frac 1 2},\; \bd^{n+1}=\bd^{n+\frac 1 2},\; \bV^{n+1}=\bV^{n+\frac 1 2}.
\end{equation}
Then, define
$({\bf u}^{n+1},v^{n+1})\in {\cal V}_F^{\eta^n}\times L^2 (0,1)$ so that the weak formulation of problem
\eqref{A2} is satisfied. Namely,
for each $({\bf q},\psi)\in{\cal V}_F^{\eta^n}\times L^2 (0,1)$ such that
${\bf q}_{|\Gamma}=\psi{\bf e}_r$, velocities $({\bf u}^{n+1},v^{n+1})$ must satisfy:
%\end{problem}
\begin{equation}
\begin{array}{c}
\displaystyle{\int_{\Omega}(1+\eta^{n})   \left( \frac{{\bf u}^{n+1}-{\bf u}^{n+\frac 1 2}}{\Delta t}\cdot{\bf q}+
\frac 1 2\left[({\bf u}^n-v^{n+\frac 1 2}r{\bf e}_r)\cdot\nabla^{\eta^n}\right]{\bf u}^{n+1}\cdot{\bf q}\right.}\\ 
[0.3cm]
\displaystyle{ \left. -\frac 1 2 \left[({\bf u}^n-v^{n+\frac 1 2}r{\bf e}_r)\cdot\nabla^{\eta^n} \right] {\bf q}\cdot{\bf u}^{n+1}\right)
+\frac{1}{2} \int_{\Omega}{v^{n+\frac 1 2}}{\bf u}^{n+1}\cdot{\bf q}}\\ 
[0.4cm]
+2\int_{\Omega}(1+\eta^n){\bf D}^{\eta^n}({\bf u}):{\bf D}^{\eta^n}({\bf q})\\ 
[0.4cm]
\displaystyle{+\rho_sh\int_0^1\frac{v^{n+1}-v^{n+\frac 1 2}}{\Delta t}\psi
=\big (
P^n_{in}\int_0^1(q_z)_{|z=0}-P^n_{out}\int_0^1(q_z)_{|z=L}\big ),}
\\
\\
{\rm with}\ \nabla^{\eta^n}\cdot{\bf u}^{n+1}=0,\quad {\bf u}^{n+1}_{|\Gamma}=v^{n+1}{\bf e}_r,
\end{array}
\label{D1Prob1}
\end{equation}
{\rm where}  $\displaystyle{P_{in/out}^n=\frac 1 {\Delta t}\int_{n\Delta t}^{(n+1)\Delta t}P_{in/out}(t)dt}$.

The existence of a unique weak solution and energy estimate are given by the following proposition.
\begin{proposition}\label{existenceA2}
Let $\Delta t > 0$, and assume that $\eta^n$ are such that $1+\eta^n  \ge R_{\rm min} > 0, n=0,...,N$. Then: 
\begin{enumerate}
\item The fluid sub-problem
defined by (\ref{D1Prob1}) has a  
unique weak solution $({\bf u}^{n+1},v^{n+1})\in {\cal V}_F^{\eta^n}\times L^2 (0,1)$;
\item Solution of problem (\ref{D1Prob1}) satisfies the following discrete energy inequality:
\begin{equation}
\begin{array}{c}
E_{kin,N}^{n+1}+\displaystyle{\frac{1}{2}\int_{\Omega_F}(1+\eta^n)|{\bf u}^{n+1}-{\bf u}^n|^2+\frac{1}{2}\|v^{n+1}-v^{n+\frac 1 2}\|^2_{L^2(0,1)}}\\
+D^{n+1}_N\leq E_{kin,N}^{n+\frac 1 2}+C\Delta t((P_{in}^n)^2+(P_{out}^n)^2),
\end{array}
\label{DEE}
\end{equation}
where the kinetic energy $E_N^n$ and dissipation $D^{n}_N$ are defined in \eqref{kenergija} and \eqref{kdisipacija},
and the constant $C$ depends only on the parameters in the problem, and not on $\Delta t$ (or $N$).
\end{enumerate}
\end{proposition}

The proof of this proposition is identical to the proof presented in \cite{BorSun} which concerns a FSI problem between an incompressible, 
viscous fluid and a thin elastic structure modeled by a linearly elastic Koiter shell model. 
The fluid sub-problems presented in \cite{BorSun} and in the present manuscript (Problem A2) are the same, except 
for the fact that $\eta$ in this manuscript satisfies the linear wave equation.
As a consequence,  the fluid domain boundary in the full, continuous problem, is not necessarily Lipschitz.
This is, however, not the case in 
the semi-discrete approximations of the fluid multi-layered structure interaction problem,  
since the regularity of the approximation $\eta^{n+1/2}$ obtained from the previous step (Problem A1) is $H^2(0,1)$,
and so the fluid domain in the semi-discretized Problem A2 is, in fact, Lipschitz. This is because $\eta^{n+1/2}$ satisfies an 
elliptic problem for the Laplace operator with the right hand-side given in terms of 
approximate velocities $v^n, v^{n+1/2} \in L^2(0,1)$ (see equation \eqref{SProb1}).
Therefore, the proof of Proposition~\ref{existenceA2} is the same as the proof of Proposition 3\cite{BorSun} (for statement 1) and 
the proof of Proposition 4\cite{BorSun} (for statement 2).

\if 1 = 0
\proof
The proof is again a consequence of the Lax-Milgram Lemma.
More precisely, donote by ${\cal U}$ the Hilbert space
\begin{equation}
{\cal U}=\{({\bf u},v)\in {\cal V}_F^{\eta^n}\times L^2 (0,1):{\bf u}_{|\Gamma}=v{\bf e}_z\},
\label{FSProb2}
\end{equation}
and define the bilinear form associated with problem \eqref{D1Prob1}:
\begin{equation*}
\begin{array}{rcl}
\displaystyle{a(({\bf u},v),({\bf q},\psi)) } &:=& 
\displaystyle{\int_{\Omega}(1+\eta^n)\left( {\bf u}\cdot{\bf q}
+\frac{\Delta t}{2}\left[ ({\bf u}^n-v^{n+\frac 1 2}r{\bf e}_r)\cdot\nabla^{\eta^n}\right]  {\bf u}\cdot{\bf q}\right.}
\\
&-& \displaystyle{ \left. \frac{\Delta t}{2} \left[({\bf u}^n-v^{n+\frac 1 2}r{\bf e}_r)\cdot\nabla^{\eta^n} \right] {\bf q}\cdot{\bf u}\right)}
\\
&+& \displaystyle{\Delta t\frac{1}{2} \int_{\Omega}{v^{n+\frac 1 2}}{\bf u}\cdot{\bf q}+\Delta t2\int_{\Omega}(1+\eta^n){\bf D}^{\eta^n}({\bf u}):{\bf D}^{\eta^n}({\bf q})}
\\
&+&\displaystyle{ \int_0^1v\psi,\qquad ({\bf q},\psi)\in{\cal U}.}
\end{array}
\end{equation*}
We need to prove that this bilinear form $a$ is coercive and continuous on ${\cal U}$.
To see that $a$ is coercive, we write
$$
a(({\bf u},v),({\bf u},v))=\int_{\Omega}(1+\eta^n+\frac{\Delta t}{2}v^{n+\frac 1 2})|{\bf u}|^2+\rho_s h\int_0^1v^2
+2 \Delta t\int_{\Omega}(1+\eta^n)|{\bf D}^{\eta^n}({\bf u})|^2.
$$
Coercivity follows immediately after recalling that $\eta^n$ are such that $1+\eta^n \ge R_{\rm min} > 0$,
which implies that $1+\eta^n+\frac{\Delta t}{2}v^{n+\frac 1 2}=1+\frac 1 2(\eta^n+\eta^{n+\frac 1 2}) \ge R_{\rm min} >0$. 

Before we prove continuity notice that for each fixed $\eta^n$ 
we have ${\cal V}_F^{\eta^n}\hookrightarrow L^4(\Omega)$. 
Indeed, for every function ${\bf u}\in{\cal V}_F^{\eta^n}$
there exists a ${\bf v}\in H^1(\Omega^{\eta^n})$ such that ${\bf u}={\bf v}^{\eta^n}$. Now, the statement follows directly
from the Sobolev embedding $H^1(\Omega^{\eta^n})\hookrightarrow L^4(\Omega^{\eta^n})$.

Therefore, by applying the generalized H\"{o}lder inequality, we obtain
$$
a(({\bf u},v),({\bf q},\psi))\leq C\left(\|{\bf u}\|_{L^2(\Omega)}\|{\bf q}\|_{L^2(\Omega)}+\|v\|_{L^2(0,1)}\|\psi\|_{L^2(0,1)}\right.
$$
$$
\left .
+
\Delta t\big ((\|{\bf u}^n\|_{{\cal V}_F^{\eta^n}}+\|v^{n+\frac 1 2}\|_{H^1(0,1)})\big )\|{\nabla^{\eta^n}{\bf u}}\|_{L^2(\Omega)}\|{\bf q}\|_{{\cal V}_F^{\eta^n}}
+\Delta t\|{\bf u}\|_{{\cal V}_F^{\eta^n}}\|{\bf q}\|_{{\cal V}_F^{\eta^n}}\right).
$$
This shows that $a$ is continuous.
The Lax-Milgram lemma now implies the existence of a unique solution $({\bf u}^{n+1},v^{n+1})$ of problem (\ref{D1Prob1}).
\qed

\begin{proposition}
For each fixed $\Delta t > 0$, solution of problem (\ref{D1Prob1}) satisfies the following discrete energy inequality:
\begin{equation}
\begin{array}{c}
E_{kin,N}^{n+1}+\displaystyle{\frac{1}{2}\int_{\Omega_F}(1+\eta^n)|{\bf u}^{n+1}-{\bf u}^n|^2+\frac{1}{2}\|v^{n+1}-v^{n+\frac 1 2}\|^2_{L^2(0,1)}}\\
+D^{n+1}_N\leq E_{kin,N}^{n+\frac 1 2}+C\Delta t((P_{in}^n)^2+(P_{out}^n)^2),
\end{array}
\label{DEE}
\end{equation}
where the kinetic energy $E_N^n$ and dissipation $D^{n}_N$ are defined in \eqref{kenergija} and \eqref{kdisipacija},
and the constant $C$ depends only on the parameters in the problem, and not on $\Delta t$ (or $N$).
\end{proposition}

TODO Ove dvije propozicije vezane uz problem A2 (egzistencija i energija) su potpuno iste kao u clanku sa Koiterom. Mozda ih ne treba dokazivati,
vec samo izreci i pozvati se na ARMA clanak?

\proof
We begin by focusing on the weak formulation \eqref{D1Prob1} in which we replace the test functions ${\bf q}$ by ${\bf u}^{n+1}$ and $\psi$ by $v^{n+1}$.
We multiply the resulting equation by $\Delta t$, and notice that the first term on the right hand-side is given by
$$
\frac{1}{2} \int_{\Omega_F}(1+\eta^n)|{\bf u}^{n+1}|^2.
$$
This is the term that contributes to the discrete kinetic energy at the time step $n+1$, but it does not have the correct form, since the discrete kinetic
energy at $n+1$ is given in terms of the structure location at $n+1$, and not at $n$, namely, the discrete kinetic energy at $n+1$ involves
$$
\frac{\rho_f}{2}\int_{\Omega_F}(1+\eta^{n+1})|{\bf u}^{n+1}|^2.
$$
To get around this difficulty it is crucial that the advection term is present in the fluid sub-problem. 
The advection term is responsible for the presence of the integral
$$
\frac{\rho_f}{2}\int_{\Omega_F} \Delta tv^{n+\frac 1 2}|{\bf u}^{n+1}|^2
$$
which can be re-written by noticing that  $\Delta tv^{n+\frac 1 2}:= (\eta^{n+1/2}-\eta^n)$ which is equal to $(\eta^{n+1}-\eta^n)$ since,
in this sub-problem $\eta^{n+1} =  \eta^{n+1/2}$. This implies
$$
\frac{\rho_f}{2}\Big (\int_{\Omega_F}(1+\eta^n)|{\bf u}^{n+1}|^2+\Delta tv^{n+\frac 1 2}|{\bf u}^{n+1}|^2\Big )=\frac{\rho_f}{2}\int_{\Omega_F}(1+\eta^{n+1})|{\bf u}^{n+1}|^2.
$$
Thus, these two terms combined provide the discrete kinetic energy at the time step $n+1$.
It is interesting to notice how the nonlinearity of the coupling at the deformed  boundary requires the presence of nonlinear advection
in order for the discrete kinetic energy of the fluid sub-problem to be decreasing in time,
and to thus satisfy the desired energy estimate.

To complete the proof one simply uses the algebraic identity $(a-b)\cdot a=\frac 1 2(|a|^2+|a-b|^2-|b|^2)$
in the same way as in the proof of Proposition~\ref{prop:Energy1}.
\qed
\fi

\vskip 0.1in
We pause for a second, and summarize what we have accomplished so far.
For a given $\Delta t >0$, the time interval $(0,T)$ was divided into $N=T/\Delta t$ sub-intervals $(t^n,t^{n+1}), n = 0,...,N-1$.
On each sub-interval  $(t^n,t^{n+1})$ we ``solved'' the coupled FSI problem by applying the Lie splitting scheme. 
First, Problem A1 was solved for the structure position and velocity, both thick and thin, and then Problem A2 was solved to update
fluid velocity and fluid-structure interface velocity.
We showed that each sub-problem has a unique solution, provided that $1+\eta^n  \ge R_{\rm min} > 0, n=0,...,N$, and that 
each sub-problem solution satisfies
an energy estimate. When combined, the two energy estimates provide a discrete version of the energy estimate \eqref{EE}. 
Thus, for each $\Delta t$ we have designed a time-marching, splitting scheme, which defines an approximate solution on $(0,T)$ of 
our main FSI problem~\eqref{FSIeqRef}-\eqref{IC_ALE}. Furthermore, the scheme is designed in
such a way that  for each $\Delta t >0$
the approximate FSI solution satisfies a discrete 
version of an energy estimate for the continuous problem.

We would like to ultimately show that, as $\Delta t \to 0$, the sequence of solutions parameterized by $N$ (or $\Delta t$),
converges to a weak solution of~\eqref{FSIeqRef}-\eqref{IC_ALE}.
Furthermore, we also need to show that $1+\eta^n  \ge R_{\rm min} > 0$ is satisfied for each $n= 0,...,N-1$.
In order to obtain this result, it is crucial to show that the discrete energy of the approximate FSI solutions defined for each $\Delta t$,
is {\sl uniformly bounded}, independently of $\Delta t$ (or $N$).
This result is obtained by the following Lemma.

\begin{lemma}\label{stabilnost}{\bf(The uniform energy estimates)}
Let $\Delta t > 0$ and $N=T/\Delta t > 0$. 
Furthermore, let $E_{N}^{n+\frac 1 2}$, $E_{N}^{n+1}$, and $D_N^j$ be the total energy and dissipation 
given by \eqref{kenergija} and \eqref{kdisipacija}, respectively.

There exists a constant $C>0$ independent of $\Delta t$ (and $N$) 
such that the following estimates hold:
\begin{enumerate}
\item
$E_{N}^{n+\frac 1 2}\leq C,\; E_{N}^{n+1}\leq C$, for all $ n = 0,...,N-1, $
\item $ \sum_{j=1}^ND_N^j\leq C,$
\item
$\displaystyle{\sum_{n=0}^{N-1}\left(\int_{\Omega_F}(1+\eta^n)|{\bf u}^{n+1}-{\bf u}^n|^2+\|v^{n+1}-v^{n+\frac 1 2}\|^2_{L^2(0,1)}\right.}$\\
\phantom {} \hskip .5in $\left. +\|v^{n+\frac 1 2}-v^{n}\|^2_{L^2(0,1)}+\|\bV^{n+1}-\bV^n\|^2_{L^2(\Omega_S)}\right)\leq C,$
\item
$
\displaystyle{\sum_{n=0}^{N-1}\Big ((\|\partial_z(\eta^{n+1}-\eta^{n})\|^2_{L^2(0,1)}
+a_S\big (\bd^{n+1}-\bd^{n},\bd^{n+1}-\bd^{n}\big )\Big )\le C.}
$
\end{enumerate}
In fact, $C = E_0 + \tilde{C} \left(\|P_{in}\|_{L^2(0,T)}^2 + \|P_{out}\|_{L^2(0,T)}^2\right)$, where $\tilde{C}$ is the constant from 
\eqref{DEE},
which depends only on the parameters in the problem.
\end{lemma}
\proof
We begin by adding the energy estimates (\ref{DEO31}) and (\ref{DEE}) to obtain
$$
E_N^{n+1}+D_N^{n+1}+\frac 1 2\Big (\int_{\Omega_F}(1+\eta^n)|{\bf u}^{n+1}-{\bf u}^n|^2+
\|v^{n+1}-v^{n+\frac 1 2}\|^2_{L^2(0,1)}
$$
$$
+\|v^{n+\frac 1 2}-v^{n}\|^2_{L^2(0,1)}+\|\bV^{n+1}-\bV^n\|^2_{L^2(\Omega_S)}+
\|\partial_z(\eta^{n+\frac 1 2}-\eta^{n})\|^2_{L^2(0,1)}
$$
$$
+a_S\big (\bd^{n+1}-\bd^{n},\bd^{n+1}-\bd^{n}\big ) \Big )\leq 
E_N^{n}+ \tilde{C} \Delta t((P_{in}^n)^2+(P_{out}^n)^2),\quad n=0,\dots,N-1.
$$
Then, we calculate the sum, on both sides, and cancel out like terms in the kinetic energy that appear on both sides of the inequality
to obtain
$$
E_N^{N}+\sum_{n=0}^{N-1} D_N^{n+1}
+\frac 1 2\sum_{n=0}^{N-1} \Big (\int_{\Omega_F}(1+\eta^n)|{\bf u}^{n+1}-{\bf u}^n|^2+
\|v^{n+1}-v^{n+\frac 1 2}\|^2_{L^2(0,1)}
$$
$$
+\|v^{n+\frac 1 2}-v^{n}\|^2_{L^2(0,1)}+\|\bV^{n+1}-\bV^n\|^2_{L^2(\Omega_S)}+
\|\partial_z(\eta^{n+\frac 1 2}-\eta^{n})\|^2_{L^2(0,1)}
$$
$$
+a_S\big (\bd^{n+1}-\bd^{n},\bd^{n+1}-\bd^{n}\big ) \Big )\leq 
E_0+ \tilde{C} \Delta t  \sum_{n=0}^{N-1} ((P_{in}^n)^2+(P_{out}^n)^2).
$$
To estimate the term involving the inlet and outlet pressure,
we recall that on every sub-interval $(t^n,t^{n+1})$ the pressure data is approximated 
by a constant which is equal to the average value of pressure over that time interval.
Therefore, we have, after using H\"{o}lder's inequality: 
$$
\Delta t\sum_{n=0}^{N-1}(P_{in}^n)^2 = \Delta t\sum_{n=0}^{N-1}\left( \frac{1}{\Delta t}\int_{n\Delta t}^{(n+1)\Delta t}P_{in}(t)dt\right)^2
\le \|P_{in}\|_{L^2(0,T)}^2.
$$
By using the pressure estimate to bound the 
right hand-side in the above energy estimate, we have obtained all the statements in the Lemma,
with the constant $C$ given by
$C = E_0 + \tilde{C}  \|P_{in/out}\|_{L^2(0,T)}^2 $.

Notice that Statement 1 can be obtained in the same way by summing from $0$ to $n-1$, for each $n$,
 instead of from $0$ to $N-1$.
\qed
\vskip 0.1in
We will use this Lemma in the next section
to show convergence of approximate solutions.

\section{Convergence of approximate solutions}\label{sec:convergence}
We define approximate solutions of problem \eqref{FSIeqRef}-\eqref{IC_ALE} on $(0,T)$ to be the functions which are piece-wise constant 
on each sub-interval $((n-1)\Delta t,n\Delta t],\ n=1\dots N$ of $(0,T)$, such that for 
$t\in ((n-1)\Delta t,n\Delta t],\ n=1\dots N,$
\begin{equation}
{\bf u}_N(t,.)={\bf u}_N^n,\ \eta_N(t,.)=\eta_N^n,\ v_N(t,.)=v_N^n,\ v^*_N(t,.)=v^{n-\frac 1 2}_N,\; \bd_N(t,.)=\bd_N^n,\; \bV_N(t,.)=\bV_N^n.
\label{aproxNS}
\end{equation}
\begin{figure}[ht]
\centering{
\includegraphics[scale=0.4]{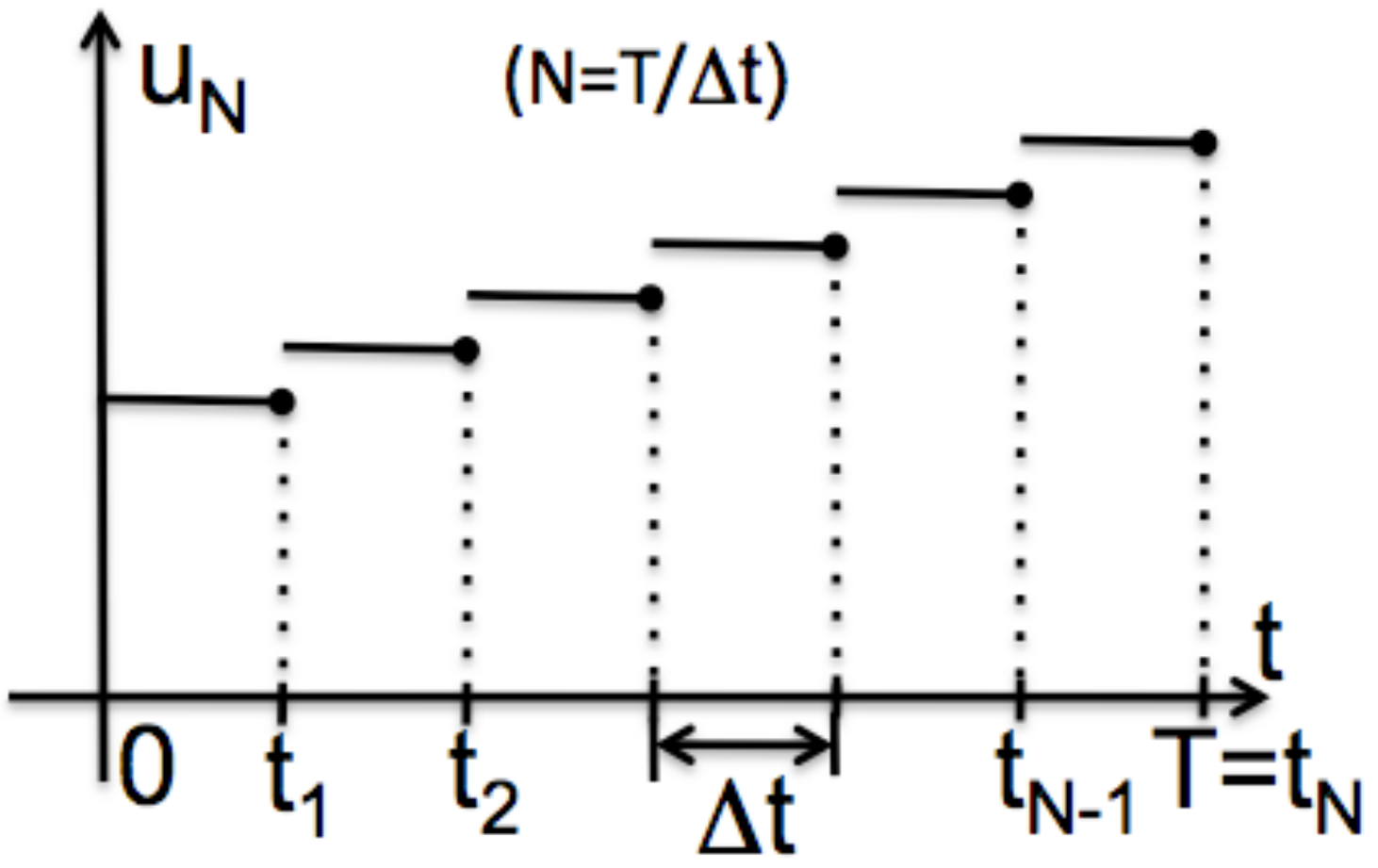}
}
\caption{A sketch of $u_N$.}
\label{fig:u_N}
\end{figure}
See Figure~\ref{fig:u_N}.
Notice that functions $v^*_N=v_N^{n-1/2}$ are determined by Problem A1 (the elastodynamics sub-problem), while functions $v_N=v^{n}_N$ are determined by 
Problem A2 (the fluid sub-problem). As a consequence, functions $v_N$ are equal to the normal trace of the fluid velocity on $\Gamma$,
i.e., ${\bf{u_N}} = v_N \bf e_r$, which may be different from  $v^*_N$. However, we will show later that $\|v_N-v_N^*\|_{L^2(0,1)}\to0$, as $N\to\infty$.

Using Lemma~\ref{stabilnost} we now show that these sequences are uniformly bounded in the appropriate solution spaces.

We begin by showing that $(\eta_N)_{N\in\N}$ is uniformly bounded in $L^\infty(0,T;H_0^1(0,1))$, and that there exists a $T > 0$
for which 
$1+\eta_N^n>0$ holds independently of $N$ and $n$. 

\begin{proposition}\label{eta_bound}
The sequence $(\eta_N)_{N\in\N}$ is uniformly bounded in 
$$
L^\infty(0,T;H_0^1(0,1)).
$$
Moreover, for $T$ small enough, we have
\begin{equation}\label{eta_bounds}
0 < R_{\rm min} \le 1 + \eta_N(t,z) \le R_{\rm max}, \ \forall N\in \N, z\in(0,1), t\in(0,T).
\end{equation}
\end{proposition}

\proof 
%The proof relies on Lemma~\ref{stabilnost} whose proof does not require the assumption  $R+\eta_N^n>0$.
The proof is similar to the corresponding proof in \cite{BorSun}, except that our structure displacement is 
bounded uniformly in $H^1_0$-norm, and not in $H^2$-norm, as in \cite{BorSun}.
This is, however, still sufficient to obtain the desired result.
More precisely, from the energy estimate in  Lemma~\ref{stabilnost}  we have
$$
\|\eta_N(t)\|^2_{L^2(0,1)}+ 
\|\partial_z \eta_N(t)\|^2_{L^2(0,1)}, 
\le C, \ \forall t\in [0,T],
$$
which implies
$$
\|\eta_N\|_{L^\infty(0,T;H_0^1(0,1))} \le C.
$$
To show that the radius $1+\eta_N$ is uniformly bounded away from zero for $T$ small enough,
we first  notice that the above inequality implies
$$
\|\eta^n_N-\eta_0\|_{H_0^1(0,1)}\leq 2C,\; n=1,\dots, N,\; N\in\N.
$$
Furthermore, we calculate
$$
\|\eta^n_N-\eta_0\|_{L^2(0,1)}\leq \sum_{i=0}^{n-1}\|\eta^{i+1}_N-\eta^{i}_N\|_{L^2(0,1)}
=\Delta t\sum_{i=0}^{n-1}\|v^{i+\frac 1 2}_N\|_{L^2(0,1)},
$$
where we recall that $\eta^{0}_N=\eta_0$.
Lemma~\ref{stabilnost} implies  that $E_N^{n+\frac 1 2}\leq C$, where $C$ is independent of $N$. This 
combined with the above inequality implies
$$
\|\eta^n_N-\eta_0\|_{L^2(0,1)}\leq C n \Delta t\leq CT,\; n=1,\dots, N,\; N\in\N.
$$
Now, since $\|\eta^n_N-\eta_0\|_{L^2(0,1)}$ and $\|\eta^n_N-\eta_0\|_{H_0^1(0,1)}$ are uniformly bounded,
we can use the interpolation inequality for Sobolev spaces, Thm. 4.17, p. 79 in  \cite{ADA}, to get
$$
\|\eta^n_N-\eta_0\|_{H^s(0,1)}\leq 2CT^{1-s},\; n=1,\dots, N,\; N\in\N,\; {\rm for}\  0<s<1.
$$ 
From Lemma~\ref{stabilnost} we see that $C$ depends on $T$ through the norms of the inlet and outlet data in such a way that $C$ is an increasing function of $T$.
Therefore, by choosing $T$ small, we can make $\|\eta^n_N-\eta_0\|_{H^s(0,1)}$ arbitrarily small for $n=1,.\dots,N$, $N\in\N$. Because of 
the Sobolev embedding of $H^s(0,1)$ into $C[0,1]$, for $s>1/2$, we can also make $\|\eta^n_N-\eta_0\|_{C[0,1]}$ arbitrarily small. 
Since the initial data $\eta_0$ is such that $1+\eta_0(z) > 0$ (due to the conditions listed in \eqref{CC}),
we see that for $T>0$ small enough, there exist $R_{\min}, R_{\rm max} > 0$, such that 
$$0 < R_{\rm min} \le 1 + \eta_N(t,z) \le R_{\rm max}, \ \forall N\in \N, z\in(0,1), t\in(0,T).$$
\qed

We will show in the end that our existence result holds not only locally in time, i.e., for small $T > 0$, 
but rather, it can be extended all the way until either $T = \infty$, or until the lateral walls of the channel touch each other.

Proposition \ref{eta_bound} implies, among other things, that the standard $L^2$-norm, and the following weighted $L^2$-norm
 are equivalent: for every $f\in L^2(\Omega_F)$, there exist
constants $C_1, C_2 > 0$, which depend only on $R_{\rm min}, R_{\rm max}$, and not on $f$ or $N$, such that
\begin{equation}\label{equivalence}
C_1\int_{\Omega_F}(1+\eta_N)f^2\leq\|f\|^2_{L^2(\Omega_F)}\leq C_2\int_{\Omega_F}(1+\eta_N)f^2.
\end{equation}
We will be using this property in the next section to prove strong convergence of approximate solutions.

Next we show that the sequences of approximate solutions for the velocity and its trace on the lateral boundary,
as well as the displacement of the thick structure and the thick structure velocity,
are uniformly bounded in the appropriate norms.
To do that, we introduce the following notation which will be useful in 
the remainder of this manuscript to prove compactness:
denote by $\tau_h$ the translation in time by $h$ of a function $f$
\begin{equation}\label{shift}
\tau_h f(t,.)=f(t-h,.),\  h\in\R. 
\end{equation}
\begin{proposition}\label{velocity_bounds}
The following statements hold:
\begin{enumerate}
\item $(v_N)_{N\in\N}$, $(v_N^*)_{N\in\N}$ are uniformly bounded in $L^\infty(0,T;L^2(0,1))$.
\item $({\bf u}_N)_{N\in\N}$ is uniformly bounded in $L^\infty(0,T;L^2(\Omega_F))$.
\item $(\nabla^{\tau_{\Delta t}\eta_N}{\bf u}_N)_{N\in\N}$ is uniformly bounded in $L^2((0,T)\times\Omega_F)$.
\item $(\bd_N)_{N\in\N}$ is uniformly bounded in $L^{\infty}(0,T;H^1(\Omega_S))$.
\item $(\bV_N)_{N\in\N}$ is uniformly bounded in $L^{\infty}(0,T;L^2(\Omega_S))$.
\end{enumerate}

\end{proposition}
\proof
The uniform boundedness of $(v_N)_{N\in\N}, (v_N^*)_{N\in\N}$, $(\bd_N)_{N\in\N}, (\bV_N)_{N\in\N}$, and the uniform boundedness of $({\bf u}_N)_{N\in\N}$ in $L^\infty(0,T;L^2(\Omega_F))$
follow directly from Statements 1 and 2 of Lemma~\ref{stabilnost}, and from the definition of 
$(v_N)_{n\in\N}, (v_N^*)_{N\in\N}$, $(\bd_N)_{N\in\N}, (\bV_N)_{N\in\N}$ and $({\bf u}_N)_{N\in\N}$ as step-functions in $t$ so that 
$$
\int_0^T \|v_N\|^2_{L^2(0,1)} dt = \sum_{n=0}^{N-1} \|v_N^n\|^2_{L^2(0,1)} \Delta t.
$$
It remains to show uniform boundedness of $(\nabla^{\tau_{\Delta t}\eta_N}{\bf u}_N)_{N\in\N}$
in $L^2((0,T)\times\Omega_F)$.
From Lemma~\ref{stabilnost} we only know that the {\sl symmetrized gradient} is bounded 
in the following way:
\begin{equation}\label{estKorn}
\displaystyle{\sum_{n=1}^N\int_{\Omega_F}(1+\eta^{n-1}_N)|{\bf D}^{\eta^{n-1}}_N({\bf u}^n_N)|^2 \Delta t\leq C.}
\end{equation}
We cannot immediately apply Korn's inequality since estimate \eqref{estKorn}
is given in terms of the transformed symmetrized gradient. Thus, there are some technical difficulties that need to be
overcome due to the fact that our problem involves moving domains.
To get around this difficulty we take the following approach.
We first transform the problem back to the physical fluid domain 
$\Omega^{\eta_N^{n-1}}_F$ which is defined by the lateral boundary $\eta_N^{n-1}$,
on which $u_N$ is defined. There, instead of the transformed gradient, we have the standard gradient,
and we can apply Korn's inequality in the usual way.
However, since the Korn constant depends on the domain, we will need a result which provides a 
universal Korn constant, independent of the family of domains under consideration. 
Indeed, a result of this kind was obtained in \cite{CDEM,Velcic,BorSun}, assuming certain domain regularity.
In particular, in \cite{Velcic} as in our previous work \cite{BorSun},
the family of domains $\Omega^{\eta_N^{n-1}}_F$ had a uniform 
Lipschitz constant, which is not the case in the present paper, since $\eta_N$ are uniformly bounded in $H_0^1(0,1)$ and not in $H^2(0,1)$.
This is why we take the approach similar to \cite{CDEM}, where the universal Korn constant was calculated explicitly, 
by utilizing the precise form of boundary data.
We have the following.

For each fixed $N \in \N$, and for all $n = 1,\dots,N$, transform ${\bf u}^n_N$ back to the physical domain 
which is determined by the location of $\eta_N^{n-1}$. We will be using super-script $N$ to denote functions defined
on physical domains:
$${\bf u}^{N,n}:={\bf u}^n_N\circ A_{\eta^{n-1}_N}, \  n=1,\dots,N,\; N\in\N.$$
By using formula (\ref{nablaeta}) we get
$$
\int_{\Omega_F}(1+\eta^{n-1}_N)|{\bf D}^{\eta^{n-1}}_N({\bf u}^n_N)|^2=\int_{\Omega_F^{\eta^{n-1}_N}}|{\bf D}({\bf u}^{N,n})|^2=
\|{\bf D}({\bf u}^{N,n})\|^2_{L^2(\Omega_F^{\eta^{n-1}_N})}.
$$
We now show that the following Korn's equality holds for the space ${\cal V}_F(t)$:
\begin{equation}\label{Korn}
\|\nabla{\bf u}^{N,n}\|^2_{L^2(\Omega_F^{\eta^{n-1}_N})}=
2\|{\bf D}({\bf u}^{N,n})\|^2_{L^2(\Omega_F^{\eta^{n-1}_N})}.
\end{equation}
Notice that the Korn constant (the number 2) is domain independent.
The proof of this Korn equality is similar to the proof in Chambolle et al. \cite{CDEM}, Lemma 6, pg. 377.
However, since our assumptions are a somewhat different from those in \cite{CDEM}, we present a sketch of the proof here.
By writing the symmetrized gradient on the right hand-side of \eqref{Korn} in terms of the gradient, and by calculating the square of the norms on both sides,
one can see that it is enough to show that
$$
\int_{\Omega_F^{\eta^{n-1}_N}}\nabla{\bf u}^n_N:\nabla^{\tau}{\bf u}^n_N=0.
$$
To simplify notation, in the proof of this equality
we omit the subscripts and superscripts, i.e. we write $\eta$ and ${\bf u}$ instead
of $\eta^{n-1}_N$ and ${\bf u}^n_N$, respectively. First, we prove the above equality for smooth functions ${\bf u}$ and then
the conclusion follows by a density argument. By using integration by parts and $\nabla\cdot{\bf u}=0$ we get
$$
\int_{\Omega_F^{\eta}}\nabla{\bf u}:\nabla^{\tau}{\bf u}=-\int_{\Omega_F^{\eta}}{\bf u}\cdot\nabla(\nabla\cdot{\bf u})+
\int_{\partial\Omega_{F}^\eta}(\nabla^{\tau}{\bf u}){\bf n}\cdot{\bf u}=\int_{\partial\Omega_{F}^\eta}(\nabla^{\tau}{\bf u}){\bf n}\cdot{\bf u},
$$
where ${\bf n}=(-\eta',1)^{\tau}$. We now show that $(\nabla^{\tau}{\bf u}){\bf n}\cdot{\bf u}=0$ on $\partial\Omega_F$. 
Since $\partial\Omega_F=\Gamma^\eta\cup\Gamma_{in/out}\cup\Gamma_b$
we consider each part of the boundary separately:
\begin{enumerate}
\item On $\Gamma^{\eta}$ we have ${\bf u}=(0,u_r)$, i.e., we have
$
u_z(z,1+\eta(z))=0.
$
Since ${\bf u}$ is smooth we can differentiate this equality w.r.t. $z$ to get
$
\partial_z u_z+\partial_r u_z\eta'=0\; {\rm on}\; \Gamma,
$
i.e., for $z\in(0,L)$.
By using $\nabla\cdot{\bf u}=0$, we get:
$
-\partial_r u_z\eta'=\partial_z u_z=-\partial_r u_r.
$
By using ${\bf n}=(-\eta',1)^{\tau}$ we get
$$
(\nabla^{\tau}{\bf u}){\bf n}\cdot{\bf u}=((\nabla^{\tau}{\bf u}){\bf n})_r u_r=(-\partial_r u_z\eta'+\partial_r u_r)u_r=0.
$$
\item On $\Gamma_{in/out}$ we have ${\bf u}=(u_z,0)$ and ${\bf n}=(\pm 1,0)$. Hence,
$$
(\nabla^{\tau}{\bf u}){\bf n}\cdot{\bf u}=u_z((\nabla^{\tau}{\bf u}){\bf n})_z=u_z(\partial_z u_z)=-u_z\partial_r u_r=0.
$$
\item On $\Gamma_b$ we have ${\bf u}=(u_z,0)$, $\partial_r u_z=0$ and ${\bf n}=(0,-1)$. Hence,
$$
(\nabla^{\tau}{\bf u}){\bf n}\cdot{\bf u}=u_z((\nabla^{\tau}{\bf u}){\bf n})_z=u_z(-\partial_z u_r)=0.
$$
\end{enumerate}
This concludes the proof of Korn's equality \eqref{Korn}.

Now, by using \eqref{Korn} and by mapping everything back to the fixed domain $\Omega_F$, 
we obtain the following Korn's equality on $\Omega_F$:
\begin{equation}\label{KornFixed}
2\int_{\Omega_F}(1+\eta^{n-1}_N)|{\bf D}^{\eta^{n-1}}_N({\bf u}^n_N)|^2=\int_{\Omega_F}(1+\eta^{n-1}_N)|\nabla^{\eta^{n-1}}_N({\bf u}^n_N)|^2.
\end{equation}
By summing equalities \eqref{KornFixed} for $n=,1\dots, N$, and by using \eqref{equivalence}, we get uniform boundedness of 
$(\nabla^{\tau_{\Delta t}\eta_N}{\bf u}_N)_{N\in\N}$
in $L^2((0,T)\times\Omega_F)$.
\qed

From the uniform boundedness of approximate sequences, the following weak and weak* convergence results follow.
\begin{lemma}{\bf (Weak and weak* convergence results)} \label{weak_convergence}
There exist subsequences $(\eta_N)_{N\in\N}$, $(v_N)_{N\in\N}$, $(v^*_N)_{N\in\N}$,
$(\bd_N)_{N\in\N}$, $(\bV_N)_{N\in\N}$ and $({\bf u}_N)_{N\in\N}$, 
and the functions $\eta \in L^{\infty}(0,T;H^1_0(0,1))$, 
$v, v^*\in L^{\infty}(0,T;L^2(0,1))$, 
$\bd\in L^{\infty}(0,T;{\cal V}_S)$, 
$\bV\in L^{\infty}(0,T;L^2(\Omega_S))$, ${\bf u}\in L^{\infty}(0,T;L^2(\Omega_F))$
and ${\bf G}\in L^2((0,T)\times\Omega_F)$
such that
\begin{equation}\label{weakconv}
\begin{array}{rcl}
\eta_N &\rightharpoonup & \eta \; {\rm weakly*}\; {\rm in}\; L^{\infty}(0,T;H^1_0(0,1)),
\\
v_N &\rightharpoonup & v\; {\rm weakly*}\; {\rm in}\; L^{\infty}(0,T;L^2(0,1)),
\\
v^*_N &\rightharpoonup & v^*\; {\rm weakly*}\; {\rm in}\; L^{\infty}(0,T;L^2(0,1)),
\\
\\
\bd_N &\rightharpoonup & \bd \; {\rm weakly*}\; {\rm in}\; L^{\infty}(0,T;H^1(\Omega_S)),
\\
\bV_N &\rightharpoonup & \bV\; {\rm weakly*}\; {\rm in}\; L^{\infty}(0,T;L^2(\Omega_S)),
\\
{\bf u}_N &\rightharpoonup & {\bf u}\; {\rm weakly*}\; {\rm in}\; L^{\infty}(0,T;L^2(\Omega_F)),
\\
\nabla^{\tau_{\Delta t}\eta_N}{\bf u}_N &\rightharpoonup & {\bf G}\; {\rm weakly}\; {\rm in}\; L^{2}((0,T)\times\Omega_F).
\end{array}
\end{equation}
Furthermore,
\begin{equation}\label{v_star}
v = v^*.
\end{equation}
\end{lemma}
\proof
The only thing left to show is that $v = v^*$. For this purpose, 
we multiply the second statement in Lemma~\ref{stabilnost} by $\Delta t$,
and 
notice again that 
$\|v_N\|_{L^2((0,T)\times(0,1))}^2=\Delta t\sum_{n=1}^N\|v^{n}_N\|^2_{L^2(0,1)}$.
This implies
$\|v_N-v^*_N\|_{L^2((0,T)\times (0,1))}\leq C\sqrt{\Delta t}$, and we have that in the limit, as $\Delta t \to 0$, $v=v^*$.
\qed

Naturally, our goal is to prove that ${\bf G}=\nabla^{\eta}{\bf u}$. However, to achieve this goal we will need some stronger
convergence properties of approximate solutions. Therefore, we postpone the proof until Section \ref{sec:limit}.

\subsection{Strong convergence of approximate sequences}

To show that the limits obtained in the previous Lemma satisfy the weak form of problem~\eqref{FSIeqRef}-\eqref{IC_ALE},
we will need to show that our sequences converge strongly in the appropriate function spaces.
The strong convergence results will be achieved by using 
the following compactness result by Simon~\cite{Simon}:
\begin{theorem}{\rm \cite{Simon}}\label{SimonLemma}
Let $X$ be a Banach space and $F\hookrightarrow L^q(0,T;X)$ with $1\leq q<\infty$. Then $F$ is a relatively compact set
in $L^q(0,T;X)$ if and only if
\begin{enumerate}
\item[(i)] $\displaystyle{\Big \{\int_{t_1}^{t_2}f(t)dt:f\in F\Big \}}$ is relatively compact in $X$, $0<t_1<t_2<T$,
\item[(ii)] $\displaystyle{\|\tau_h f-f\|_{L^q(h,T;X)}\to 0}$ as $h$ goes to zero, uniformly with respect to $f\in F$.
\end{enumerate}
\end{theorem}
We used this result in \cite{BorSun} to show compactness, but the proof was simpler because of the higher regularity of 
the lateral boundary of the fluid domain, namely, of the fluid-structure interface.
In the present paper we need to obtain some additional regularity for the fluid velocity ${\bf u}_N$ on $\Omega_F$
and its trace ${\bf v}_N$ on the lateral boundary, before we can use Theorem~\ref{SimonLemma} to show strong convergence
of our approximate sequences.
Notice, we only have that our fluid velocity on $\Omega_F$ is uniformly bounded in $L^2(\Omega_F)$, plus a condition that the transformed 
gradient $\nabla^{\tau_{\Delta t}\eta_N} {\bf u}_N$ is uniformly bounded in $L^2$. Since $\eta$ is not Lipschitz, we cannot get that
the gradient $\nabla {\bf u}_N$ is uniformly bounded in $L^2$ on $\Omega_F$. 
This lower regularity of $\eta_N$ will give us some trouble when showing regularity of ${\bf u}_N$ on $\Omega_F$, namely it will imply
lower regularity of ${\bf u}_N$ in the sense that ${\bf u} \in H^s(\Omega_F)$, for $0<s<1/2$, and not $H^1(\Omega_F)$. Luckily, according to
the trace theorem in \cite{BorisTrag}, this will still allow us to make sense of the trace of ${\bf u}_N$ 
on $\Gamma$. 
More precisely,  we prove the following Lemma.
\begin{lemma}\label{HsBrzina}
The following statements hold:
\begin{enumerate}
\item $({\bf u}_N)_{N\in\N}$ is uniformly bounded in $L^2(0,T;H^s(\Omega_F))$, $0<s<1/2$;
\item $(v_N)_{N\in\N}$ is uniformly bounded in $L^2(0,T;H^{s/2}(0,1))$, $0<s<1/2$.
\end{enumerate}
\end{lemma}

\proof
We start by  mapping  the fluid velocity ${\bf u}_N$ defined on $\Omega_F$, back to the physical fluid domain 
with the lateral boundary 
$\tau_{\Delta t}\eta_N(t,z) = \eta_N(t-\Delta t,z)$.
%(Recall, $\eta_N$ is a piece-wise constant function in $t$, incorporating a sequence of functions $\eta_N^n$.)
We denote by ${\bf u}^N(t,.)$ the fluid velocity on the physical domain $\Omega_{\tau_{\Delta t}\eta_N}$:
$$
{\bf u}^N(t,.)={\bf u}_N(t,.)\circ A^{-1}_{\tau_{\Delta t}\eta_N}(t),\; N\in\N.
$$
As before, we use sub-script $N$ to denote fluid velocity defined on the physical space.
From  \eqref{dercomp} we see that 
$$
\nabla {\bf u}^N = \nabla^{\tau_{\Delta t}\eta_N}{\bf u}_N.
$$
Proposition \ref{velocity_bounds}, statement 3, implies that the sequence $(\nabla {\bf u}^N)_{N\in\N}$ 
is uniformly bounded in $L^2$, and so
we have that 
$\|{\bf u}^N\|_{L^2(0,T;H^1(\Omega_{\tau_{\Delta t}\eta}))}$ is uniformly bounded.

Now, from the fact that the fluid velocities ${\bf u}^N$ defined on the physical domains are uniformly bounded in $H^1$,
we would like to obtain a similar result for the velocities ${\bf u}_N$ defined on the reference domain $\Omega_F$.
For this purpose, we recall that the functions $\eta_N, N\in \N$ that are involved in the ALE mappings
$A_{\tau_{\Delta t}\eta_N}(t)$, $N \in \N$, 
are uniformly bounded in $H^1(0,1)$.
This is, unfortunately,  not sufficient to obtain uniform boundedness of the gradients $(\nabla u_N)_{N\in \N}$ in $L^2(\Omega_F)$.
However, from the Sobolev embedding $H^1(0,1)\hookrightarrow C^{0,1/2}(0,1)$ we have that 
the sequence $(\eta_N)_{N\in\N}$ is uniformly bounded in $L^{\infty}(0,T;C^{0,1/2}(0,1))$.
This will help us obtain uniform boundedness of $({\bf u}_N)_{n\in\N}$ in a slightly lower-regularity space,
namely in the space $L^2(0,T;H^s(\Omega_F))$, $0 < s < 1/2$.
To see this,  we first notice
that ${\bf u}_N$ on $\Omega_F$ can be expressed in terms of function ${\bf u}^N$ defined on $\Omega_{\tau_{\Delta t}\eta_N}$ as
\begin{equation}\label{vel}
{\bf u}_N(t,\tilde z,\tilde r) = {\bf u}^N (t,\tilde z, (1+ \tau_{\Delta t}\eta_N)(t,\tilde z))\tilde r), \ (\tilde z, \tilde r) \in \Omega_F.
\end{equation}
Therefore, ${\bf u}_N$ can be written as an $H^1$-function ${\bf u}^N$ composed with a $C^{0,1/2}$-function $\eta_N$,
in the way described in \eqref{vel}.
The following Lemma, proved in \cite{BorisTrag}, implies that ${\bf u}_N$ belongs to 
a space with asymmetric regularity (more regular in $\tilde r$ than
in $\tilde z$) in the sense that ${\bf u}_N \in L^2(0,1;H^s(0,1)), 0 < s < 1/2$, and $\partial_{\tilde{r}}  {\bf u}_N \in L^2(0,1;L^2(0,1))$.
We use notation from Lions and Magenes \cite{LionsMagenes}, pg.~10, to denote the corresponding function space by
$$
W(0,1;s)=\{f:f\in L^2(0,1;H^s(0,1))),\; \partial_{\tilde{r}} f\in L^2(0,1;L^2(0,1))\}.
$$
More precisely, Lemma 3.3 from \cite{BorisTrag} states the following:
\begin{lemma}{\rm \cite{BorisTrag}}\label{TraceLemma}
Let $\eta\in C^{0,\alpha}$, $0<\alpha<1$, and let $u\in H^1(\Omega_{\eta})$.
Define 
\begin{equation}\label{velALE}
{\tilde u}(\tilde{r},\tilde{z})=u(\tilde{z},(1+\eta(\tilde{z}))\tilde{r}),\quad (\tilde{z},\tilde{r})\in\Omega_F.
\end{equation}
Then ${\tilde u}\in W(0,1;s)$ for $0<s<\alpha$. 
\end{lemma}
Thus, Lemma~\ref{TraceLemma} implies that ${\bf u}_N(t,.)\in W(0,1;s)$ for $0<s<1/2$.
Now, using the fact $W(0,1;s)\hookrightarrow H^s(\Omega_F)$
we get 
$$
\|{\bf u}_N(t,.)\|^2_{H^s(\Omega_F)}\leq C\|{\bf u}^N(t,.)\|^2_{H^1(\Omega_{\eta(t-\Delta t)})},\; a.a.\; t\in (0,T),\; 0< s< 1/2.
$$
By integrating the above inequality w.r.t. $t$ we get the first statement of Lemma~\ref{HsBrzina}.

To prove the second statement of Lemma~\ref{HsBrzina} we use Theorem~3.1 of \cite{BorisTrag}, 
which states that the notion of the 
trace for the functions of the form \eqref{vel} for which ${\bf u}^N \in H^1$ and $\eta_N \in C^{0,1/2}$,
can be defined in the sense of $H^{s/2}$, $0<s<1/2$.
For completeness, we state Theorem 3.2 of \cite{BorisTrag} here. 

\begin{theorem}{\rm \cite{BorisTrag}}\label{TraceTm}
Let $\alpha<1$ and let $\eta$ be such that 
$$
\eta\in C^{0,\alpha}(0,1),\; \eta({z})\geq \eta_{min}>-1,\; {z}\in [0,1],\;  \eta(0)=\eta(1)=1.
$$
Then, the  trace operator 
$$
\gamma_{\eta}:C^1(\overline{\Omega_{\eta}}) \to  C(\Gamma)
$$
that associates to each function $u \in C^1(\overline{\Omega_{\eta}})$ 
its ``Lagrangian trace'' $u(\tilde{z},1+\eta(\tilde{z})) \in C(\Gamma)$,
defined via \eqref{velALE} for $\tilde{r}=1$,
$$
\gamma_{\eta}:u \mapsto  u(\tilde{z},1+\eta(\tilde{z})), 
$$
can be extended by continuity to a linear operator from $H^1(\Omega_{\eta})$ to 
$H^s(\Gamma)$ for $0\leq s<{\alpha}/{2}$.
\end{theorem}

By recalling that $v_N=({\bf u}_N)_{|\Gamma}$, this proves the second statement of Lemma~\ref{HsBrzina}.
\qed

Notice that the difficulty associated with bounding the gradient of ${\bf u}_N$ is somewhat artificial, since the gradient
of the fluid velocity ${\bf u}^N$ defined on the physical domain is, in fact, uniformly bounded (by Proposition \ref{velocity_bounds}).
Namely, the difficulty is imposed
by the fact that we decided to work with the problem defined on a fixed domain $\Omega_F$, and not on a family of moving domains. 
This decision, however, simplifies other parts of the main existence proof. The ``expense'' that we had to pay for
this decision is embedded in the proof of Lemma~\ref{HsBrzina}.

\if 1 = 0

The main ingredient in getting the ``integral equicontinuity'' estimate (ii) from Theorem~\ref{SimonLemma} is
Lemma \ref{stabilnost}. Namely, if we 
multiply the third equality of Lemma \ref{stabilnost} by $\Delta t$ we get:
\begin{equation}\label{tauetaN}
\|\tau_{\Delta t}{\bf u}_{N}-{\bf u}_{N}\|^2_{L^2((0,T)\times\Omega_F)}
+\|\tau_{\Delta t} v_{N}-v_{N}\|^2_{L^2((0,T)\times (0,1))}\leq C {\Delta t}.
\end{equation}
This is ``almost'' (ii) except that in this estimate $\varepsilon$ depends on $\Delta t$ (i.e., $N$),
which is not sufficient to show equicontinuity (ii). 
We need to show that estimate \eqref{tauetaN} holds for all the functions  $(v_N)_{N\in\N}$, $({\bf u}_{N})_{N\in\N}$
independently of $N\in\N$. This is why we need to work a little harder to get the following compactness result.

\fi

We are now ready to use Theorem~\ref{SimonLemma} to prove compactness of 
the sequences $(v_N)_{N\in\N}$ and $({\bf u}_{N})_{N\in\N}$.
\begin{theorem}\label{u_v_convergence}
Sequences $(v_N)_{N\in\N}$ and $({\bf u}_{N})_{N\in\N}$ are relatively compact in 
$L^2(0,T;L^2(0,1))$ and $L^2(0,T;L^2(\Omega_F))$, respectively.
\label{LKompaktnost}
\end{theorem}

\proof We use Theorem~\ref{SimonLemma} with $q = 2$, and $X=L^2$. We verify that  both assumptions (i) and (ii) hold.

{Assumption (i):} To show that the sequences $(v_N)_{N\in\N}$ and $({\bf u}_{N})_{N\in\N}$ are relatively compact in 
$L^2(0,1)$ and $L^2(\Omega_F)$, respectively, we use Lemma~\ref{HsBrzina} and the compactness of the embeddings
$H^s(\Omega_F)\hookrightarrow L^2(\Omega_F)$ and
$H^{s/2}(0,1)\hookrightarrow L^2(0,1)$, respectively, for $0<s<1/2$.
Namely, from  Lemma~\ref{HsBrzina} we know that 
sequences $({\bf u}_N)_{N\in\N}$ and $(v_N)_{N\in\N}$ are uniformly bounded 
in $L^2(0,T;H^s(\Omega_F))$ and $L^2(0,T;H^{s/2}(0,1))$, respectively, for $0<s<1/2$.
The compactness of the embeddings $H^s(\Omega_F)\hookrightarrow L^2(\Omega_F)$ and
$H^{s/2}(0,1)\hookrightarrow L^2(0,1)$ verify Assumption (i) of Theorem~\ref{SimonLemma}.

Assumption (ii): We prove that the ``integral equicontinuity'', stated in assumption (ii) of Theorem \ref{SimonLemma},
holds for the sequence $(v_N)_{N\in\N}$. Analogous reasoning can be used for $({\bf u}_{N})_{N\in\N}$.
Thus, we want to show that for each $\varepsilon > 0 $, there exists a $\delta > 0$ such that
\begin{equation}\label{equicontinuity1}
\|\tau_h v_N-v_N\|^2_{L^2(\omega;L^2(0,1))}< \varepsilon, \quad \forall |h| < \delta, \ {\rm independently\ of}\ N\in\N,
\end{equation}
where $\omega$ is an arbitrary compact subset of $\Omega$.
Indeed, we will show that for each $\varepsilon > 0$, the following choice of $\delta$:
$$\delta :=\min\{{\rm dist}(\omega,\partial \Omega) /2,\varepsilon/(2C)\}$$
provides the desired estimate, where 
$C$ is the constant from Lemma \ref{stabilnost} (independent of $N$).

Let  $h$ be an arbitrary real number whose absolute value is less than $\delta$.
We want to show that \eqref{equicontinuity1} holds for all $\Delta t = T/N$.
This will be shown in two steps. First, we will show that \eqref{equicontinuity1} holds for the case when $\Delta t \ge h$ (Case 1), and then for the case when $\Delta t < h$ (Case 2).

A short remark is in order: For a given $\delta > 0$, we will have $\Delta t < \delta$ for infinitely many $N$,  and both cases will apply.
For a finite number of functions $(v_N)$, we will, however, have that $\Delta t \ge \delta$. For those functions \eqref{equicontinuity1}
needs to be proved for all $\Delta t$ such that $|h| < \delta \le \Delta t$, which falls into Case 1 bellow. Thus, Cases 1 and 2 cover all
the possibilities.

\begin{figure}[ht]
\centering{
\includegraphics[scale=0.5]{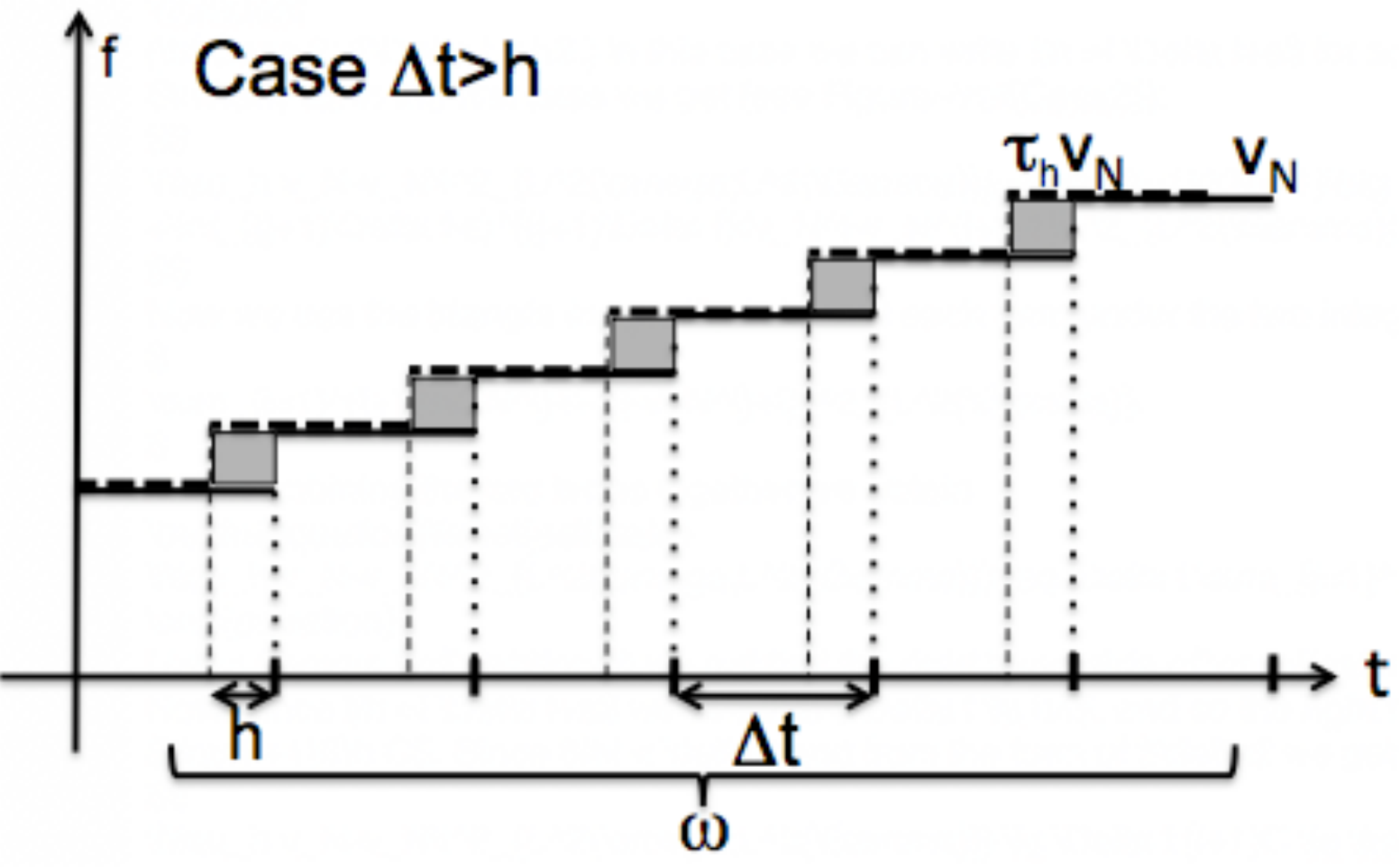}
\includegraphics[scale=0.5]{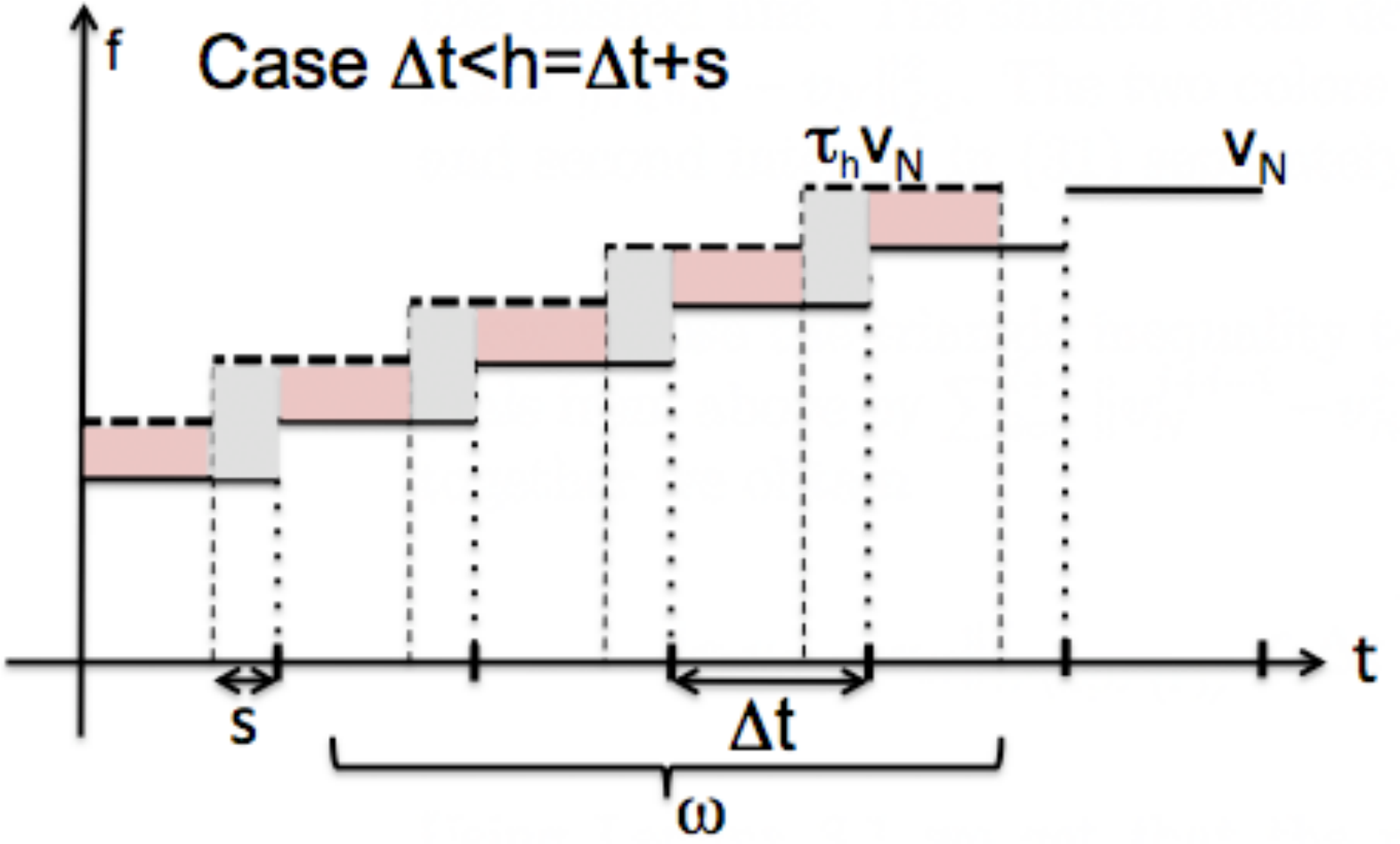}
}
\caption{Left panel--Case 1: $\Delta t \ge h$. The graph of $v_N$ is shown in solid line, while the
graph of the shifted function $\tau_h v_N$ is shown in dashed line.
The shaded area denotes the non-zero contributions to the norm $\|\tau_h v_N-v_N\|^2_{L^2}$.
Right panel--Case2: $\Delta t < h = \Delta t + s, 0<s<\Delta t$. The graph of $v_N$ is shown in solid line, while the
graph of the shifted function $\tau_h v_N$ is shown in the dashed line.
The shaded areas denote non-zero contributions to the norm $\|\tau_h v_N-v_N\|^2_{L^2}$.
The two colors represent the contributions to the first and second integral in \eqref{est} separately.}
\label{fig:case1}
\end{figure}

\vskip 0.1in 
\noindent
{\bf Case 1:  $\Delta t \ge h$.} We calculate the shift by $h$ to obtain (see Figure~\ref{fig:case1}, left):
$$
\|\tau_h v_N-v_N\|^2_{L^2(\omega;L^2(0,1))}\leq\sum_{j=1}^{N-1}\int_{j\Delta t-h}^{j\Delta t}\|v_N^j-v_N^{j+1}\|_{L^2(0,1)}^2=
$$
$$
=h\sum_{j=1}^{N-1}\|v_N^j-v_N^{j+1}\|^2_{L^2(0,1)}\leq hC<\varepsilon/2 < \varepsilon.
$$
The last inequality follows from $|h| < \delta \le \varepsilon/(2C)$.

\vskip 0.2in
\noindent
{\bf Case 2: $\Delta t < h$.} In this case we can write $h =l \Delta t+s$ for some $\l\in\N$, $0<s\leq \Delta t$. 
Similarly, as in the first case, we get (see Figure~\ref{fig:case1}, right):
\begin{equation} \label{est}
\begin{array}{rcl}
\|\tau_h v_N-v_N\|^2_{L^2(\omega;L^2(0,1))}&=&
\displaystyle{\sum_{j=1}^{N-l-1}\big (\int_{j \Delta t}^{(j+1)\Delta t-s}\|v_N^j-v_N^{j+l}\|^2_{L^2(0,1)}}\\
&+&\displaystyle{\int_{(j+1)\Delta t-s}^{(j+1)\Delta t}\|v_N^j-v_N^{j+l+1}\|^2_{L^2(0,1)}\big ).}
\end{array}
\end{equation}
Now we use the triangle inequality to bound each term under the two integrals from above by
$
\sum_{i=1}^{l+1}\|v_N^{j+i-1}-v_N^{j+i}\|^2_{L^2(0,1)}.
$
After combining the two terms together we obtain
\begin{equation}\label{estimate}
\|\tau_h v_N-v_N\|^2_{L^2(\omega;L^2(0,1))}\leq \Delta t \sum_{j=1}^{N-l-1}\sum_{i=1}^{l+1}\|v_N^{j+i-1}-v_N^{j+i}\|^2_{L^2(0,1)}.
\end{equation}
Using Lemma \ref{stabilnost} we get that the right hand-side of \eqref{estimate} is bounded by $\Delta t (l+1)C$.
Now, since $h =l \Delta t+s$ we see that $\Delta t \le h/l$, and so the right hand-side of \eqref{estimate} is bounded by 
$\frac{l+1}{l}h C$. Since $|h| < \delta$ and from the form of $\delta$ we get
$$
\|\tau_h v_N-v_N\|^2_{L^2(\omega;L^2(0,1))} \le \Delta t (l+1)C \le \frac{l+1}{l}h C \le \frac{l+1}{l} \frac{\varepsilon}{2} < \varepsilon.
$$
Thus, if we set $\omega=[\delta/2,T-\delta/2]$ we have shown:
$$
\|\tau_{\delta/2} v_N-v_N\|^2_{L^2(\delta/2,T-\delta/2;L^2(0,1))}<\varepsilon,\quad N\in\N.
$$
To show that condition (ii) from Theorem~\ref{SimonLemma} holds it remains to estimate 
$\|\tau_{\delta/2} v_N-v_N\|^2_{L^2(T-\delta/2,T;L^2(0,1))}$.
From the first inequality in Lemma \ref{stabilnost} 
(boundedness of $v_N^{n+\frac{i}{2}}, i = 1,2$ in $L^2(0,1)$) we have
$$\int_{T-\delta/2}^T\|\tau_{\delta/2} v_N-v_N\|^2_{L^2(0,1)}\leq \frac{\delta}{2}2C<\varepsilon,\quad N\in\N.$$

Thus, we have verified all the assumptions of Theorem~\ref{SimonLemma}, and so the compactness result for $(v_N)_{N\in\N}$ follows
from Theorem~\ref{SimonLemma}.
Similar arguments imply compactness of $(\mathbf{u}_N)_{N\in\N}$.
\qed

To show compactness of $(\eta_N)_{N\in\N}$ we use the approach similar to that in \cite{BorSun},
except that, due to the weaker regularity properties of $\eta_N$, we will have to use
different embedding results (Hilbert interpolation inequalities). 
In the end, compactness of the sequence of lateral boundary approximation will follow due to the Arzel\`a- Ascoli Theorem. 

\begin{figure}[ht]
\centering{
\includegraphics[scale=0.4]{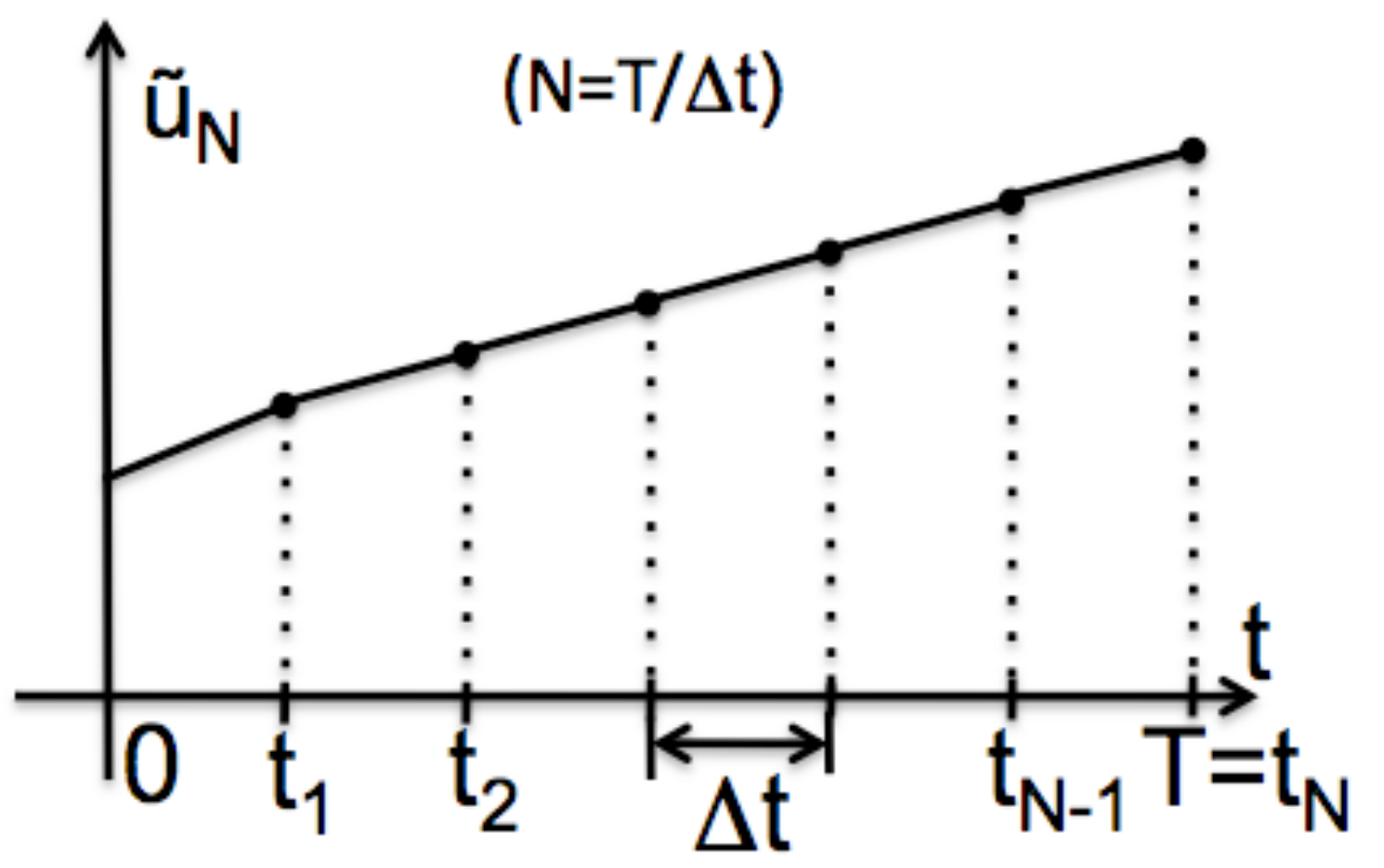}
}
\caption{A sketch of $\tilde{\bf u}_N$.}
\label{fig:u_N_tilde}
\end{figure}

As in \cite{BorSun}, we start by introducing a slightly different set of approximate functions of $\mathbf u$, $v$, $\eta$ and $\bV$.
Namely, for each fixed $\Delta t$ (or $N \in \N$), define $\tilde{\bf u}_N$, $\tilde{\eta}_N$, $\tilde{v}_N$ and $\tilde{\bV}_N$ 
to be continuous, {\sl linear} on
each sub-interval $[(n-1)\Delta t,n\Delta t]$, and such that for $n=0,\dots,N$:
\begin{equation}\label{tilde}
\begin{array}{c}
\tilde{\bf u}_N(n\Delta t,.)={\bf u}_N(n\Delta t,.),\ \tilde{v}_N(n\Delta t,.)={v}_N(n\Delta t,.),\\
 \tilde{\eta}_N(n\Delta t,.)={\eta}_N(n\Delta t,.), 
\;\tilde{\bV}_N(n\Delta t,.)={\bV}_N(n\Delta t,.),
\end{array}
\end{equation}
See Figure~\ref{fig:u_N_tilde}.
A straightforward calculation gives the following inequalities (see \cite{Tem}, p. 328)
\begin{equation}\label{TemLemma}
\begin{array}{c}
\displaystyle{\|v_N-\tilde{v}_N\|^2_{L^2(0,T;L^2(0,1))}\leq{\frac{\Delta t}{3}}\sum_{n=1}^N\|v^{n+1}-v^{n}\|^2_{L^2(0,1)},}
\\ 
\displaystyle{\|{\bf u}_N-\tilde{\bf u}_N\|^2_{L^2(0,T;L^2(\Omega_F))}\leq{\frac{\Delta t}{3}}\sum_{n=1}^N\|{\bf u}^{n+1}-{\bf u}^{n}\|^2_{L^2(\Omega_F)},}
\\ 
\displaystyle{\|\eta_N-\tilde{\eta}_N\|^2_{L^2(0,T;L^2(0,1))}\leq{\frac{\Delta t}{3}}\sum_{n=1}^N\|\eta^{n+1}-\eta^{n}\|^2_{L^2(0,1)},}
\\ 
\displaystyle{\|\bV_N-\tilde{\bV}_N\|^2_{L^2(0,T;L^2(\Omega_S))}\leq{\frac{\Delta t}{3}}\sum_{n=1}^N\|\bV^{n+1}-\bV^{n}\|^2_{L^2(\Omega_S)},}
\end{array}
\end{equation}
We now observe that 
$$\displaystyle{\partial_t\tilde{\eta}_N(t)=\frac{\eta^{n+1}-\eta^n}{\Delta t}= \frac{\eta^{n+1/2}-\eta^n}{\Delta t}=v^{n+\frac 1 2}},\ t\in (n\Delta t,(n+1)\Delta t),$$ 
and so, since $v^*_N$ was defined in  \eqref{aproxNS} as
a piece-wise constant function defined via $v^*_N(t,\cdot)=v^{n+\frac 1 2}$, for $t\in(n\Delta t,(n+1)\Delta t]$,
we see that
\begin{equation}\label{derivativeeta}
\partial_t\tilde{\eta}_N=v^*_N\ a.e.\ {\rm on}\ (0,T).
\end{equation}
%Functions $\tilde{\rm u}_N$ and $\tilde{v}_N$ will be used later in the proof when analysis of time-derivatives is performed. TODO
By using Lemma \ref{stabilnost} (the boundedness of $E_N^{n+\frac{i}{2}}$), we get
$$
(\tilde{\eta}_N)_{N\in\N}\; {\rm is\; bounded\; in}\; L^{\infty}(0,T;H^1_0(0,1))\cap W^{1,\infty}(0,T;L^2(0,1)).
$$
We now use the following  result on continuous embeddings:
\begin{equation}\label{compact_embedding}
L^{\infty}(0,T;H^1_0(0,1))\cap W^{1,\infty}(0,T;L^2(0,1)) \hookrightarrow C^{0,1-\alpha}([0,T]; H^{\alpha}(0,1)),
\end{equation}
for $0<\alpha<1$.
This result follows from the standard Hilbert interpolation inequalities, see \cite{LionsMagenes}.
A slightly different result (assuming higher regularity) was also used in \cite{CG} to deal with a set of mollifying functions approximating
a solution to a moving-boundary problem between a viscous fluid and an elastic plate.
From \eqref{compact_embedding} we see that $(\tilde{\eta}_N)_{N\in\N}$ is also bounded (uniformly in $N$) in $C^{0,1-\alpha}([0,T]; H^{\alpha}(0,1))$.
Now, from the continuous embedding of $H^{\alpha}(0,1)$ into $H^{\alpha-\epsilon}$, and
by applying the Arzel\`{a}-Ascoli Theorem, we conclude that sequence $(\tilde{\eta}_N)_{N\in\N}$ has a convergent subsequence,
which we will again denote by $(\tilde{\eta}_N)_{N\in\N}$, such that
$$\tilde{\eta}_N\rightarrow \tilde{\eta}\ {\rm in} \  C([0,T];H^s(0,1)),\ 0< s<1.$$
Since \eqref{TemLemma} implies that
$(\tilde{\eta}_N)_{N\in\N}$ and $(\eta_N)_{N\in\N}$ have the same limit, we have $\eta=\tilde{\eta}\in C([0,T];H^s(0,1))$,
where $\eta$ is the weak* limit of $(\eta_N)_{N\in\N}$, discussed in \eqref{weakconv}.
Thus, we have
$$\tilde{\eta}_N\rightarrow {\eta}\ {\rm in} \  C([0,T];H^s(0,1)),\ 0< s<1.$$
We can now prove the following Lemma:
\begin{lemma}\label{etaconvergence}
$\eta_N\rightarrow \eta$ in $L^{\infty}(0,T;H^s(0,1))$, $0 < s<1$.
\label{konveta}
\end{lemma}

\proof
The proof is similar to the proof of Lemma 3 in \cite{BorSun}. The result follows from the continuity in time of $\eta$, and from the fact that $\tilde{\eta}_N\rightarrow {\eta}\ {\rm in} \  C([0,T];H^s(0,1))$, for $ 0 < s<1$, applied to the inequality
$$
\|\eta_N(t)-\eta(t)\|_{H^s(0,1)}=
\| \eta_N(t) - \eta(n\Delta t) + \eta(n\Delta t) -\eta(t)\|_{H^s(0,1)}
$$
$$
=\| \eta_N(n\Delta t) - \eta(n\Delta t) + \eta(n\Delta t) -\eta(t)\|_{H^s(0,1)}
$$
$$
\le \| \eta_N(n\Delta t) - \eta(n\Delta t)\| + \| \eta(n\Delta t) -\eta(t)\|_{H^s(0,1)}
$$
$$
= 
\| \tilde\eta_N(n\Delta t) - \eta(n\Delta t)\|_{H^s(0,1)}  + \| \eta(n\Delta t) -\eta(t)\|_{H^s(0,1)}.
$$
\if 1 = 0
Namely, let $\varepsilon>0$. 
Then, from the continuity of $\eta$ in time we have that there exists a $\delta t>0$ such that 
$$
\|\eta(t_1)-\eta(t_2)\|_{H^s(0,1)}<\frac{\varepsilon}{2}, \  {\rm for} \ t_1,t_2\in [0,T], \ {\rm and} \ |t_1-t_2|\leq\delta t.
$$
Furthermore, from the convergence $\tilde{\eta}_N\rightarrow {\eta}\ {\rm in} \  C([0,T];H^s(0,1)),\ 0 < s<1$, we know that 
there exists an $N^* \in \N$ such that 
$$
\|\tilde{\eta}_N-\eta\|_{C([0,T];H^s(0,1))}< \frac{\varepsilon}{2}, \ \forall N\geq N^*.
$$ 
Now, let $N$ be any natural number such that $N>\max\{N^*,T/\delta t\}$. Denote by $\Delta t = T/N$, and let $t \in [0,T]$.
Furthermore, let $n\in\N$ be such that $(n-1)\Delta t<t\leq n\Delta t$.
Recall that $\tilde{\eta}_N(n\Delta t)=\eta_N(n \Delta t)=\eta_N(t)$ from the definition of $\tilde{\eta}_N$ and $\eta_N$.
By using this, and by combining the two estimates above, we get
$$
\|\eta_N(t)-\eta(t)\|_{H^s(0,1)}=
\| \eta_N(t) - \eta(n\Delta t) + \eta(n\Delta t) -\eta(t)\|_{H^s(0,1)}
$$
$$
=\| \eta_N(n\Delta t) - \eta(n\Delta t) + \eta(n\Delta t) -\eta(t)\|_{H^s(0,1)}
$$
$$
\le \| \eta_N(n\Delta t) - \eta(n\Delta t)\| + \| \eta(n\Delta t) -\eta(t)\|_{H^s(0,1)}
$$
$$
= 
\| \tilde\eta_N(n\Delta t) - \eta(n\Delta t)\|_{H^s(0,1)}  + \| \eta(n\Delta t) -\eta(t)\|_{H^s(0,1)}
<\varepsilon.
$$
Here, the first term is bounded by $\varepsilon/2$ due to the convergence $\tilde{\eta}_N \to \eta$, while
the second term is bounded by $\varepsilon/2$ due to the continuity of $\eta$. 
Since the obtained estimate is uniform in $N$ and $t$, the statement of the Lemma is proved.
\fi
\qed
\vskip 0.1in

We summarize the strong convergence results obtained in Theorem~\ref{u_v_convergence} and Lemma~\ref{etaconvergence}.
We have shown that there exist subsequences $(\mathbf u_N)_{N\in\N}$, $(v_ N)_{N\in\N}$ and $(\eta_N)_{N\in\N}$
such that
\begin{equation}\label{strong_convergence}
\begin{array}{rcl}
\mathbf u_N &\to& {\bf u} \  {\rm in}\  L^2(0,T;L^2(\Omega_F)),\\
v_ N &\to& v \  {\rm in}\  L^2(0,T;L^2(0,1)),\\
\tau_{\Delta t} \mathbf  u_N &\to& u \  {\rm in}\  L^2(0,T;L^2(\Omega_F)),\\
\tau_{\Delta t} v_ N &\to& v \  {\rm in}\  L^2(0,T;L^2(0,1)),\\
\eta_N &\rightarrow& \eta  \  {\rm in}\  L^{\infty}(0,T;H^s(0,1)), \ 0\leq s<1.
\end{array}
\end{equation}
Because of the uniqueness of derivatives, we also have $v=\partial_t\eta$ in the sense of distributions.
The statements about the convergence of $(\tau_{\Delta t} \mathbf  u_N)_{N\in\N}$ and $(\tau_{\Delta t} v_ N)_{N\in\N}$ follow directly from 
\begin{equation}\label{tauetaN}
\|\tau_{\Delta t}{\bf u}_{N}-{\bf u}_{N}\|^2_{L^2((0,T)\times\Omega_F)}
+\|\tau_{\Delta t} v_{N}-v_{N}\|^2_{L^2((0,T)\times (0,1))}\leq C {\Delta t},
\end{equation}
which is obtained after multiplying the third equality of Lemma \ref{stabilnost} by $\Delta t$.

Furthermore, one can also show that 
subsequences $(\tilde{v}_N)_N$, $(\tilde{\bf u}_N)_N$ and  $(\tilde{\bV}_N)_N $ also converge to $v$, ${\bf u}$ and $\bV$ respectively. 
More precisely,
\begin{equation}
\begin{array}{rcl}
\tilde{\mathbf u}_N &\to& {\bf u} \  {\rm in}\  L^2(0,T;L^2(\Omega_F)),\\
\tilde v_ N &\to& v \  {\rm in}\  L^2(0,T;L^2(0,1)),\\
\tilde{\bV}_N &\rightharpoonup& {\bV} \  {\rm weakly*}\; {\rm in}\  L^{\infty}(0,T;L^2(\Omega_S))
\end{array}
\end{equation}
This statement follows directly from the inequalities \eqref{TemLemma}
and Lemma \ref{stabilnost}, which provides uniform boundedness of the sums on the right hand-sides of the inequalities.

We conclude this section by showing one last convergence result that will be used in the next section to prove
that the limiting functions satisfy weak formulation of the FSI problem. Namely,
we want to show that
\begin{equation}\label{zeta}
\begin{array}{rcl}
{\eta_N} & \to &\eta\; {\rm in}\; L^{\infty}(0,T;C[0,1]),\\
{\tau_{\Delta t}\eta_N} & \to &\eta\; {\rm in}\; L^{\infty}(0,T;C[0,1]).
\end{array}
\end{equation}
The first statement is a direct consequence of Lemma \ref{konveta}
in which we proved that $\eta_N\to\eta$ in $L^{\infty}(0,T;H^s(0,1))$, $0<s<1$.
For $s>\frac 1 2$ we immediately have 
\begin{equation}
{\eta_N} \to \eta\; {\rm in}\; L^{\infty}(0,T;C[0,1]).
\end{equation}
To show convergence of the shifted displacements $\tau_{\Delta t} \eta_N$
to the same limiting function $\eta$, 
we recall that 
\begin{equation*}
\begin{array}{rcl}
\tilde{\eta}_N & \to &\eta\; {\rm in}\; C([0,T];H^s[0,L]),\; 0<s<1,
\end{array}
\end{equation*}
and that $(\tilde{\eta}_N)_{N\in\N}$ is uniformly bounded in $C^{0,1-\alpha}([0,T]; H^{\alpha}(0,1))$, $0 < \alpha < 1$.
Uniform boundeness of $(\tilde{\eta}_N)_{N\in\N}$ in $C^{0,1-\alpha}([0,T]; H^{\alpha}(0,1))$ implies that 
there exists a constant $C > 0$, independent of $N$, such that
$$
 \|  \tilde{\eta}_N((n-1)\Delta t) - \tilde\eta_N(n \Delta t)\|_{H^{\alpha}(0,1)}  \le C |\Delta t|^{1-\alpha}.
$$
This means that for each $\varepsilon > 0$, there exists an $N_1 > 0$ such that 
$$
 \|  \tilde{\eta}_N((n-1)\Delta t) - \tilde\eta_N(n \Delta t)\|_{H^{\alpha}(0,1)}  \le \frac{\varepsilon}{2}, {\rm for\ all}\ N \ge N_1.
$$ 
Here, $N_1$ is chosen by recalling that $\Delta t = T/N$, and so the right hand-side implies that we want an $N_1$ such that
$$
C\left(\frac{T}{N}\right)^{1-\alpha} <  \frac{\varepsilon}{2} \ {\rm for\ all}\ N\ge N_1.
$$
Now, convergence $\tilde{\eta}_N \to \eta\; {\rm in}\; C([0,T];H^s[0,1]),\; 0<s<1,$ implies that for each $\varepsilon > 0$,
there exists an $N_2 > 0$ such that
$$
 \| \tilde \eta_N (n \Delta t)  - \eta(t) \|_{H^s(0,1)} < \frac{\varepsilon}{2},\ {\rm for \ all}\ N \ge N_2.
$$
We will use this to show that for each $\varepsilon > 0$ there exists an $N^* \ge  {\rm max}\{N_1,N_2\}$, such that
$$
\displaystyle{\| \tau_{\Delta t} \tilde{\eta}_N(t) - \eta(t) \|_{H^s(0,1)} < \varepsilon, \ {\rm for \ all}\ N\ge N^*.}
$$ 
Indeed, let $t\in(0,T)$.  Then there exists an $n$  such that $t\in( (n-1)\Delta t, n \Delta t]$. We calculate
\begin{eqnarray*}
&\displaystyle{\| \tau_{\Delta t} \tilde{\eta}_N(t) - \eta(t) \|_{H^s(0,1)} =
\| \tau_{\Delta t} \tilde{\eta}_N(t) - \tilde\eta_N(n \Delta t) +\tilde \eta_N (n \Delta t)  - \eta(t) \|_{H^s(0,1)}}\\
&\displaystyle{=\|  \tilde{\eta}_N((n-1)\Delta t) - \tilde\eta_N(n \Delta t) +\tilde \eta_N (n \Delta t)  - \eta(t) \|_{H^s(0,1)} } \\
&\displaystyle{ \le \|  \tilde{\eta}_N((n-1)\Delta t) - \tilde\eta_N(n \Delta t)\|_{H^s(0,1)}   + \| \tilde \eta_N (n \Delta t)  - \eta(t) \|_{H^s(0,1)} .}
\end{eqnarray*}
The first term is less than $\varepsilon/2$ by  the uniform boundeness of $(\tilde{\eta}_N)_{N\in\N}$ in $C^{0,1-\alpha}([0,T]; H^{\alpha}(0,1))$,
while the second term is less than $\varepsilon/2$ by the convergence of $\tilde{\eta}_N$ to $\eta$ in 
$C([0,T];H^s[0,1]),\; 0<s<1$.

Now, since $\tau_{\Delta t}\tilde{\eta}_N=\widetilde{({\tau_{\Delta t}\eta}_N)}$, we can use the same argument as
in Lemma~\ref{etaconvergence} to show that sequences
$\widetilde{({\tau_{\Delta t}\eta}_N)}$ and $\tau_{\Delta t}\eta_N$ both converge to the same limit $\eta$ in $L^{\infty}(0,T;H^s(0,1))$,
for $0<s<1$.

\section{The limiting problem and weak solution}\label{sec:limit}

Next we want to show that the limiting functions satisfy the weak form \eqref{VFRef} of the full fluid-structure iteration problem.
In this vein, one of the things that needs to be considered is what happens in the limit as $N \to \infty$, i.e., as $\Delta t \to 0$,
of the weak form of the fluid sub-problem (\ref{D1Prob1}).
Before we pass to the limit we must observe that, unfortunately, the velocity test functions in (\ref{D1Prob1}) depend of $N$!
More precisely, they depend on $\eta^n_N$ because of the requirement that
the transformed divergence-free condition $\nabla^{\eta^n_N}\cdot{\bf q} = 0$ must
be satisfied. 
This is a consequence of the fact that we mapped our fluid sub-problem onto a fixed domain $\Omega_F$.
Therefore, we need to take special care when constructing suitable velocity test functions and passing to the limit in \eqref{D1Prob1}.

\subsection{Construction of the appropriate test functions}
We begin by recalling that test functions  $({\bf q},\psi,\bpsi)$ for the 
limiting problem are defined by the space ${\cal Q}$, given in \eqref{Q}, which depends on $\eta$.
Similarly, the test spaces for the approximate problems depend on $N$ through the dependence on $\eta_N$.
We had to deal with the same difficulty in \cite{BorSun} where a FSI problem with a thin structure modeled by 
the full Koiter shell equations
was studied.
The only difference is that, due to the lower regularity of
the fluid-structure interface in the present paper we also need to additionally show that the sequence of gradients of the fluid velocity 
converges weakly to the gradient of the limiting velocity, and pay special attention when taking the limits in the weak formulation
of the FSI problem.

To deal with the dependence of test functions on $N$, we follow the same ideas as those presented in \cite{CDEM,BorSun}.
We restrict ourselves to a dense subset ${\cal X}$ of all test functions in ${\cal Q}$ 
that is independent of $\eta_N$ even for the approximate problems.
We construct the set  ${\cal X}$ to consist 
of the test functions $({\bf q},\psi,\bpsi) \in{\cal X} = {\cal X}_F \times{\cal X}_W\times {\cal X}_S$, such that
the velocity components ${\bf q}\in {\cal X}_F$ are smooth, independent of $N$, and $\nabla\cdot {\mathbf q} = 0$.
Such functions can be constructed as an algebraic
sum of the functions ${\bf q}_0$ that have compact support in $\Omega_{\eta}\cup\Gamma_{in}\cup\Gamma_{out}\cup\Gamma_b$, plus a function ${\bf q}_1$,
which captures the behavior of the solution at the boundary $\Gamma_\eta$. 
More precisely, let $\Omega_{min}$ and $\Omega_{max}$ denote the fluid domains associated with the radii $R_{min}$ and $R_{max}$,
respectively. 
\begin{enumerate}
\item {\bf Definition of test functions $({\bf q}_0,0,{\bf 0})$ on $(0,T)\times\Omega_{max}\times\Omega_S$}: 
Consider all smooth functions ${\bf q}$ with compact support in $\Omega_{\eta}\cup\Gamma_{in}\cup\Gamma_{out}\cup\Gamma_b$,
and such that $\nabla \cdot {\bf q} = 0$.
Then we can extend ${\bf q}$ by $0$ to a divergence-free vector field on $(0,T)\times\Omega_{{max}}$.  
This defines ${\bf q}_0$. 

Notice that since ${\eta_N}$ converge uniformly to $\eta$,
there exists an $N_q > 0$ such that supp$({\bf q}_0)\subset\Omega_{{\tau_{\Delta t}\eta_N}}$, $\forall N \ge N_q$. 
Therefore, ${\bf q}_0$ is well defined  on infinitely many approximate domains $\Omega_{{\tau_{\Delta t}\eta_N}}$.
%Therefore, we can transform ${\bf q}$ to the reference domain $\Omega$ using $A_{{\tau_{\Delta t}\eta_N}}$.
\item {\bf Definition of test functions $({\bf q}_1,\psi,\bpsi)$ on $(0,T)\times\Omega_{max}\times\Omega_S$}: 
Consider  $\psi\in C_c^1([0,T);H^2_0(0,1))$. 
Define
\begin{equation*}
{\bf q}_1 :=\left\{
\begin{array}{l}
 \left.  
   \begin{array}{l}
       {\rm A\  constant \ extension \  in\  the\  vertical}\\
      {\rm direction\  of}  \ \psi{\bf e}_r\ {\rm on}\  \Gamma_\eta:  {\bf q}_1  := (0,\psi(z))^T;\\
       {\rm  Notice}\ {\rm div} {\bf q}_1 = 0.\\
  \end{array}
 \right\}
    {\rm on}\  \Omega_{{max}}\setminus\Omega_{{min}},\\
    \\
\left.
   \begin{array}{l}
      {\rm A \ divergence-free\  extension\  to } \  \Omega_{min}\\
     {\rm (see,\ e.g.\  \cite{GB2}, \ p.\  127).} 
  \end{array}
\right\}
    \ {\rm on}\ \Omega_{min}.
\end{array}
\right.
\end{equation*}
 From the construction
it is clear that ${\bf q}_1$ is also defined on $\Omega_{{\tau_{\Delta t}\eta_N}}$ for each $N$, and so it can be 
mapped onto the reference domain $\Omega$ by the transformation $A_{{\tau_{\Delta t}\eta_N}}$.
We take $\bpsi\in H^1(\Omega_S)$ such that $\bpsi(t,z,1)=\psi(t,z)$.
\end{enumerate}
For any test function $({\bf q},\psi,\bpsi)\in{\cal Q}$ it is easy to see that
the velocity component ${\bf q}$ can then be written as ${\bf q}={\bf q} - {\bf q}_1 + {\bf q}_1$,
where ${\bf q} - {\bf q}_1$ can be approximated by divergence-free functions ${\bf q}_0$
that  have compact support in $\Omega_{\eta}\cup\Gamma_{in}\cup\Gamma_{out}\cup\Gamma_b$.
Therefore, one can easily see that functions $({\bf q},\psi) = ({\bf q}_0 + {\bf q}_1,\psi)$ in ${\cal X}$ satisfy the following properties:
\begin{itemize}
\item
 ${\cal X}$  is dense in the space  ${\cal Q}$ of all test functions defined on the physical, moving domain $\Omega_\eta$,
defined by \eqref{Q}; furthermore, $\nabla\cdot \mathbf q = 0, \forall \mathbf q \in {\cal X}_F$.
\item For each ${\bf q} \in {\cal X}_F$, define
$$
\tilde{\mathbf q} = \mathbf q \circ A_{\eta}.
$$
The set $\{(\tilde{\mathbf q},\psi,\bpsi)  |  \tilde{\mathbf q} = \mathbf q \circ A_{\eta},  \mathbf q \in {\cal X}_F,\; \psi\in{\cal X}_S,\;
 \bpsi \in {\cal X}_S\}$ is  dense in 
the space ${\cal Q}_\eta$ of all test functions defined on the fixed, reference domain $\Omega_F$,
defined by \eqref{Q_eta}.
\item For each ${\bf q} \in {\cal X}_F$, define
$$
 {\bf q}_N:={\bf q}\circ A_{{\tau_{\Delta t}\eta_N}}.
$$
 Functions $ {\bf q}_N$ are defined on the fixed domain $\Omega_F$, and they satisfy
$\nabla^{{\tau_{\Delta t}\eta_N}}\cdot{\bf q_N}=0$.
\end{itemize}

 Functions $ {\bf q}_N$  will serve as test functions for approximate problems associated with the sequence of  domains 
 $\Omega_{\tau_{\Delta t}\eta_N}$, while functions $\tilde{\mathbf q}$ will serve as test functions associated with the domain $\Omega_\eta$.
 Both sets of test functions are defined on $\Omega_F$.
\begin{lemma}\label{testf}
For every $(\bf q,\psi,\bpsi)\in {\cal X}$ we have 
${\bf q}_N\rightarrow\tilde{\bf q}$
uniformly in $L^\infty(0,T; C(\Omega_F))$.
\end{lemma}
\proof
By the Mean-Value Theorem we get:
\begin{eqnarray*}
|{\bf q}_N(t,z,r)-\tilde{\bf q}(t,z,r)|&=&|{\bf q}(t,z,(1+\tau_{\Delta t}\eta_N)r)-{\bf q}(t,z,(1+\eta)r)|\\
    &=&|\partial_r{\bf q}(t,z,\zeta)r|\ |\eta(t,z)-\eta_N(t-\Delta t,z)|.
\end{eqnarray*}
The uniform convergence of ${\bf q}_N$ follows 
from the uniform convergence of $\eta_N$, since $\mathbf q$ are smooth. 
\qed

We are now ready to identify the weak limit ${\bf G}$ from Lemma~\ref{weak_convergence}.
\begin{proposition}\label{gradeta}
${\bf G}=\nabla^{\eta}{\bf u}$, where ${\bf G}$, ${\bf u}$ and $\eta$ are the weak and weak* limits given by Lemma \ref{weak_convergence}.
\end{proposition}

\proof
As in Lemma~\ref{HsBrzina}, it will be helpful to map the approximate fluid velocities and the limiting fluid velocity onto the physical domains. 
For this purpose, we introduce the following functions
\begin{equation*}
\begin{array}{rlrl}
{\bf u}^N(t,.)&={\bf u}_N(t,.)\circ A^{-1}_{\tau_{\Delta t}\eta_N}(t), \ \ &\tilde{\bf u}(t,.)&={\bf u}(t,.)\circ A^{-1}_{\eta}(t),\\
\\
\chi^N {\bf f}(t,\bf{x}) &=
        \left\{\begin{array}{ll}
                {\bf f}, &\ {\bf x}\in \Omega_{\tau_{\Delta t}\eta_N}(t)\\
                            0,  & \ {\bf x} \notin \Omega_{\tau_{\Delta t}\eta_N}(t)
        \end{array} \right. ,
&   \chi {\bf f}(t,\bf{x}) &=
        \left\{\begin{array}{ll}
                {\bf f}, & {\bf x}\in \Omega_{\eta}(t)\\
                            0,  & {\bf x} \notin \Omega_{\eta}(t)
        \end{array} \right. ,
\end{array}
\end{equation*}
where $A$ is the ALE mapping defined by (\ref{RefTrans}),  $\eta$ is the weak* limit
$
\eta_N \rightharpoonup  \eta \; {\rm in}\; L^{\infty}(0,T;H^1_0(0,1))
$ 
satisfying the uniform convergence property \eqref{zeta},
and ${\bf f}$ is an arbitrary function defined on the physical domain. 
Notice, again, that superscript $N$ is used to denote a function defined on the physical domain, while subscript $N$ 
is used denote a function defined on the fixed domain $\Omega_F$.

The proof consists of three main steps: (1) we will first show that $\chi^N {\bf u}^N \to  \chi \tilde{\bf u}$ 
strongly in $L^2((0,T)\times \Omega_{max})$,
then, by using step (1), we will show (2) $\chi^N \nabla {\bf u}^N \to  \chi \nabla \tilde{\bf u}$ 
weakly in $L^2((0,T)\times \Omega_{max})$, and, finally by using (2) we will show 
(3) $\int_0^T \int_{\Omega_F} {\bf G}:\tilde{\bf q} = \int_0^T \int_{\Omega_F} \nabla^\eta {\bf u} : \tilde{\bf q}$
for every test function $\tilde{\bf q} = {\bf q} \circ A_\eta$.

\vskip 0.1in
{\bf STEP 1.} We will show that 
$
\| \chi^N {\bf u}^N - \chi \tilde{\bf u} \|_{L^2((0,T)\times \Omega_{max})} \to 0.
$
To achieve this goal, we introduce the following auxiliary functions
$$\tilde{\bf u}^N(t,.)={\bf u}_N(t,.)\circ A^{-1}_{\eta}(t),$$
which will be used in the following estimate
$$
\|\chi^N{\bf u}^N-\chi\tilde{\bf u}\|_{L^2((0,T)\times\Omega_{max})}\leq 
\|\chi^N{\bf u}^N-\chi\tilde{\bf u}^N\|_{L^2((0,T)\times\Omega_{max})}+\|\chi\tilde{\bf u}^N-\chi\tilde{\bf u}\|_{L^2((0,T)\times\Omega_{max})}.
$$
The second term on the right-hand side converges to zero because of the strong convergence of ${\bf u}_N$ to ${\bf u}$ on 
the reference domain $\Omega_F$, namely,
$$
\|\chi\tilde{\bf u}^N-\chi\tilde{\bf u}\|^2_{L^2(\Omega_{max})}=\int_{\Omega_F}(1+\eta)|{\bf u}_N-{\bf u}|^2\to 0.
$$

\begin{figure}[ht]
\centering{
\includegraphics[scale=0.35]{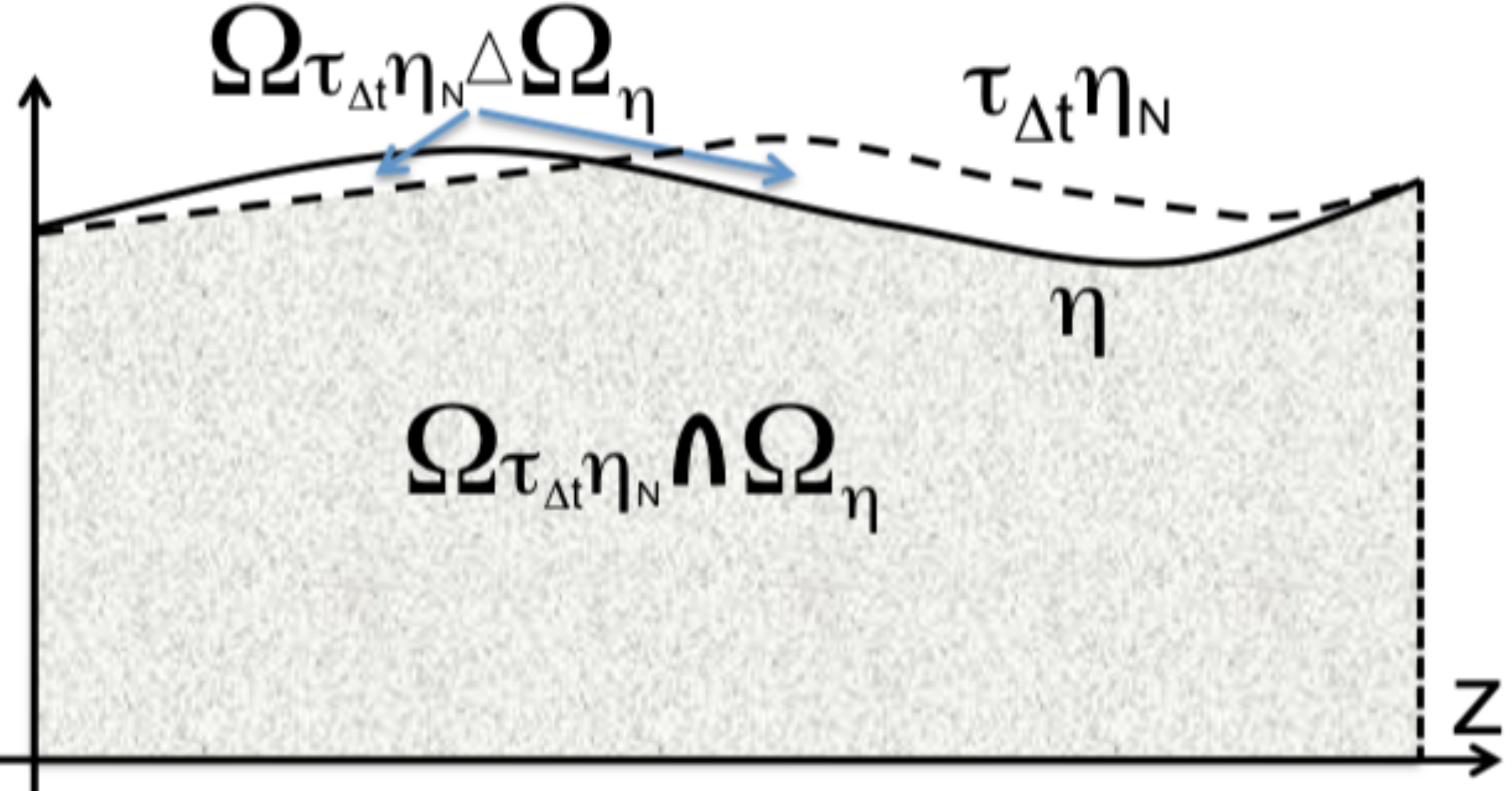}
}
\caption{A sketch of the fluid domains in STEP 1.}
\label{final}
\end{figure}
To show that the first term on the right-hand side converges to zero, first notice that 
$$
\int_0^T\int_{\Omega_{max}}|\chi\tilde{\bf u}^N-\chi^N{\bf u}^N|^2=(\int_0^T\int_{\Omega_{\eta}(t)\triangle\Omega_{\tau_{\Delta t}\eta_N}(t)}
+\int_0^T\int_{\Omega_{\eta}(t)\cap\Omega_{\tau_{\Delta t}\eta_N}(t)})|\chi\tilde{\bf u}^N-\chi^N{\bf u}^N|^2.
$$
Here $A\Delta B := (A\cup B)\setminus (A\cap B)$. See Figure~\ref{final}.
Because of the uniform convergence \eqref{zeta} we can make the measure $|\Omega_{\eta}(t)\triangle\Omega_{\tau_{\Delta t}\eta_N}(t)|$ 
arbitrary small. Furthermore,
by Propostions \ref{eta_bound} and \ref{velocity_bounds} we have that the sequence $(\chi\tilde{\bf u}^N-\chi^N{\bf u}^N)_{N\in\N}$ is
uniformly bounded in $L^2((0,T)\times\Omega_{max}).$ 
Therefore, for every $\varepsilon > 0$, there exists an $N_0\in\N$ such that for every $N\geq N_0$ we have
\begin{equation}\label{AuxEq1}
\int_0^T\int_{\Omega_{\eta}(t)\triangle\Omega_{\tau_{\Delta t}\eta_N}(t)}
|\chi\tilde{\bf u}^N-\chi^N{\bf u}^N|^2<\frac{\varepsilon}{2}.
\end{equation}
To estimate the second term, we need to measure the relative difference between the function ${\bf u}_N$ composed with 
$A_\eta^{-1}(t)$, denoted by $\tilde{\bf u}^N$, and the same function ${\bf u}_N$ composed with $A^{-1}_{\tau_{\Delta t}\eta_N(t)}$,
denoted by ${\bf u}^N$. We will map them both on the same domain and work with one function ${\bf u}_N$, while the convergence
of the $L^2$-integral will be obtained by estimating the difference 
in the ALE mappings. More precisely, we introduce the set
$\omega=A_{\eta}^{-1}(\Omega_{\eta}(t)\cap\Omega_{\tau_{\Delta t}\eta_N}(t))\subset\Omega_F$.
Now, we use the properties of the ALE mapping $A_{\eta}$ and the definitions of $\tilde{\bf u}^N,\;{\bf u}^N$ to get
$$
\int_0^T\int_{\Omega_{\eta}(t)\cap\Omega_{\tau_{\Delta t}\eta_N}(t)}|\chi\tilde{\bf u}^N-\chi^N{\bf u}^N|^2
=
\int_0^T\int_{\omega}\frac{1}{1+\eta}|{\bf u}_N-{\bf u}_N\circ A^{-1}_{\tau_{\Delta t}\eta_N(t)}\circ A_{\eta(t)}|^2
$$
$$
=\int_0^T\int_{\omega}\frac{1}{1+\eta(t,z)}|{\bf u}_N(t,z,r)-{\bf u}_N(t,z,\frac{1+\eta(t,z)}{1+\tau_{\Delta t}\eta_N(t,z)}r)|^2
$$
$$
= \int_0^T\int_{\omega} \left|\partial_r{\bf u}_N(t,z,\zeta)r \left(1-\frac{1+\eta(t,z)}{1+\tau_{\Delta t}\eta_N(t,z)}\right) \right|^2
$$
Now because of the uniform convergence \eqref{zeta} of the sequence $(\tau_{\Delta t}\eta_N)_{N\in\N}$,
and the uniform boundedness of $(\|\partial_r{\bf u}_N\|_{L^2(\Omega_F)})_{N\in\N}$,
which is consequence of Proposition \ref{velocity_bounds}, we can take $N_1\geq N_0$ such that
$$
\int_0^T\int_{\Omega_{\eta}(t)\cap\Omega_{\tau_{\Delta t}\eta_N}(t)}|\chi\tilde{\bf u}^N-\chi^N{\bf u}^N|^2< \frac{\varepsilon}{2},\; N\geq N_1.
$$
This inequality, together with \eqref{AuxEq1} proofs that $\chi^N {\bf u}^N \to  \chi \tilde{\bf u}$ 
strongly in $L^2((0,T)\times \Omega_{max})$.

\vskip 0.1in
{\bf STEP 2.} 
We will now show that 
$
\begin{array}{rcl}
\chi^N\nabla{\bf u}^N \rightharpoonup \chi\nabla\tilde{\bf u} \; {\rm weakly}\; {\rm in}\; L^2((0,T)\times\Omega_{max}). 
\end{array}
$
First notice that from
$$
\nabla{\bf u}^N=\nabla^{\tau_{\Delta t}\eta_N}{\bf u}_N
$$
and from uniform boundedness of $(\nabla^{\tau_{\Delta t}\eta_N}{\bf u}_N)_{N\in\N}$ in $L^2((0,T)\times\Omega_F)$,
established in Proposition~\ref{velocity_bounds}, 
we get that the sequence $(\chi^N\nabla{\bf u}^N)_{N\in\N}$ converges weakly  in $L^2((0,T)\times\Omega_{max})$.
Let us denote the weak limit
of $(\chi^N\nabla{\bf u}^N)_{N\in\N}$ by $\tilde{\bf G}$. 
Therefore,
$$
\int_0^T\int_{\Omega_{max}}\tilde{\bf G}\cdot{\boldsymbol\phi}=\lim_{N\to\infty}\int_0^T\int_{\Omega_{max}}\chi^N\nabla{\bf u}^N\cdot
{\boldsymbol\phi},\quad {\boldsymbol\phi}\in C^{\infty}_c((0,T)\times\Omega_{max}).
$$
We want to show that $\tilde{\bf G} = \chi\nabla\tilde{\bf u}$.

For this purpose, we first consider the set $(\Omega_{max}\setminus\Omega_{\eta}(t))$ and show that $\tilde{\bf G} = 0$ there, and then
the set $\Omega_\eta(t)$ and show that $\tilde{\bf G} = \nabla\tilde{\bf u}$ there.

 Let ${\boldsymbol\phi}$ be a test function such that 
${\rm supp}{\boldsymbol\phi}\subset(0,T)\times\Big (\Omega_{max}\setminus\Omega_{\eta}(t)\Big )$.
Using the uniform convergence of the sequence $\tau_{\Delta t}\eta_N$, obtained in \eqref{zeta}, 
there exists an $N_{\boldsymbol \phi}$ such that $\chi^N({\bf x})=0$, $N\geq N_{\boldsymbol\phi}$, ${\bf x}\in{\rm supp}{\boldsymbol\phi}$. Therefore, we have
$$
\int_0^T\int_{\Omega_{max}}\tilde{\bf G}\cdot{\boldsymbol\phi}=\lim_{N\to\infty}\int_0^T\int_{\Omega_{max}}\chi^N\nabla{\bf u}^N\cdot
{\boldsymbol\phi}=0.
$$
Thus, $\tilde{\bf G}=0$ on $(0,T)\times\Big (\Omega_{max}\setminus\Omega_{\eta}(t)\Big )$. 

Now, let us take a test function ${\boldsymbol\psi}$ such that ${\rm supp}{\boldsymbol\psi}\subset (0,T)\times\Omega_{\eta}(t)$. 
Again using the same argument as before, as well as the uniform convergence of the sequence $\tau_{\Delta t}\eta_N$, obtained in \eqref{zeta},
we conclude that there exists an $N_{\boldsymbol \psi}$ such that 
$\chi^N({\bf x})=1$, $N\geq N_{\boldsymbol\psi}$, ${\bf x}\in{\rm supp}{\boldsymbol\psi}$. Therefore, we have
$$
\int_0^T\int_{\Omega_{max}}\tilde{\bf G}\cdot{\boldsymbol\psi}=\lim_{N\to\infty}\int_0^T\int_{\Omega_{max}}\chi^N\nabla{\bf u}^N\cdot
{\boldsymbol\psi}=\lim_{N\to\infty}\int_0^T\int_{\Omega_{\eta}(t)}\nabla{\bf u}^N\cdot
{\boldsymbol\psi}.
$$
From the strong convergence $\chi^N {\bf u}^N \to  \chi \tilde{\bf u}$  obtained in STEP 1, 
 we have that on the set ${\rm supp}{\boldsymbol\psi}$, ${\bf u}^N \to   \tilde{\bf u}$ in the sense of distributions, and so,
 on the same set  ${\rm supp}{\boldsymbol\psi}$, 
$\nabla{\bf u}^N \to \nabla\tilde{\bf u}$ in the sense of distributions. Therefore we have
$$
\int_0^T\int_{\Omega_{max}}\tilde{\bf G}\cdot{\boldsymbol\psi}=
\lim_{N\to\infty}\int_0^T\int_{\Omega_{\eta}(t)}\nabla{\bf u}^N\cdot
{\boldsymbol\psi}=\int_0^T\int_{\Omega_{\eta}(t)}\nabla\tilde{\bf u}\cdot{\boldsymbol\psi}.
$$
Since this conclusion holds for all the test functions ${\boldsymbol\psi}$ supported in $(0,T)\times\Omega_{\eta(t)}$,
from the uniqueness of the limit,
we conclude $\tilde {\bf G} = \nabla\tilde{\bf u}$ in $(0,T)\times\Omega_{\eta(t)}$.

Therefore, we have shown that
 $$
\begin{array}{rcl}
\chi^N\nabla{\bf u}^N & \rightharpoonup &\chi\nabla\tilde{\bf u}  \; {\rm weakly}\; {\rm in}\; L^2((0,T)\times\Omega_{max}). 
\end{array}
$$

\vskip 0.1in
{\bf STEP 3.} We want to show that $\int_0^T\int_{\Omega_F}{\bf G}:\tilde{\bf q}=\int_0^T\int_{\Omega_F}\nabla^{\eta}{\bf u}:\tilde{\bf q}$
for every test function $\tilde{\bf q} = {\bf q} \circ A_\eta$, ${\bf q} \in {\cal X}_F$.
This will follow from STEP 2, the uniform boundedness and convergence of the gradients $\nabla^{\tau_{\Delta t}\eta_N}\tilde{\bf u}_N$ provided by
Lemma~\ref{weak_convergence}, and from the strong convergence of the test functions ${\bf q}_N \to \tilde{\bf q}$ provided by Lemma \ref{testf}.
More precisely, we have that for every $\tilde{\bf q} = {\bf q} \circ A_\eta$, ${\bf q} \in {\cal X}_F$
$$
\displaystyle{\int_0^T\int_{\Omega_F}{\bf G}:\tilde{\bf q}=\lim_{N\to\infty}\int_0^T\int_{\Omega_F}\nabla^{\tau_{\Delta t}\eta_N}{\bf u}_N:{\bf q}_N}
$$
$$
\displaystyle{=\lim_{N\to\infty}\int_0^T\int_{\Omega_{{max}}}\frac{1}{1+\tau_{\Delta t}\eta_N}\chi^N\nabla{\bf u}^N:{\bf q}
=\int_0^T\int_{\Omega_{\eta}}\frac{1}{1+\eta}\nabla\tilde{\bf u}:{\bf q}=
\int_0^T\int_{\Omega_F}\nabla^{\eta}{\bf u}:\tilde{\bf q}}.
$$
Here, we have used from (\ref{nablaeta}) that 
$\nabla{\bf u}^N=\nabla^{\tau_{\Delta t}\eta_N}{\bf u}_N$, and $\nabla\tilde{\bf u}=\nabla^{\eta}{\bf u}.$
This completes proof.
\qed

\begin{corollary}\label{testf_nabla}
For every $(\bf q,\psi,\bpsi)\in {\cal X}$ we have 
$$
\nabla^{\tau_{\Delta t}\eta_N}{\bf q}_N\rightarrow\nabla^{\eta}\tilde{\bf q},\;{\rm in}\;L^2((0,T)\times\Omega_F).
$$
\end{corollary}

\proof
Since ${\tau_{\Delta t}\eta_N}{\bf q}_N$ and $\tilde{\bf q}$ are the test functions for the velocity fields, the same arguments
as in Proposition~\ref{gradeta} provide weak convergence of $(\nabla^{\tau_{\Delta t}\eta_N}{\bf q}_N)_{N\in\N}$.
To prove strong convergence it is sufficient to prove the convergence of norms 
$\displaystyle{\|\nabla^{\tau_{\Delta t}\eta_N}{\bf q}_N\|_{L^2(\Omega_F)}\to\|\nabla^{\eta}\tilde{\bf q}\|_{L^2(\Omega_F)}}$.
This can be done, by using the uniform convergence of $(\tau_{\Delta t}\eta_N)_{N\in\N}$, in the following way:
$$
\|\nabla^{\tau_{\Delta t}\eta_N}{\bf q}_N\|^2_{L^2(\Omega_F)}=
\int_0^T\int_{\Omega_{max}}\chi^N\frac{1}{1+\tau_{\Delta t}\eta_N}|\nabla{\bf q}|^2\to
\int_0^T\int_{\Omega_{max}}\chi\frac{1}{1+\eta}|\nabla{\bf q}|^2
$$
$$
=\int_0^T\int_{\Omega_{F}}|\nabla^{\eta}\tilde{\bf q}|^2=\|\nabla^{\eta}\tilde{\bf q}\|_{L^2(\Omega_F)}^2.
$$
The notation used here is analogous to that used in the proof of Proposition~\ref{gradeta}.
\qed

Before we can pass to the limit in the weak formulation of the approximate problems, 
there is one more useful observation that we need.
Namely, notice that although ${\bf q}$ are smooth functions both in the spatial variables and in time, 
the functions ${\bf q}_N$ are discontinuous at $n\Delta t$ because ${\tau_{\Delta t}\eta_N}$ is a step function
in time. As we shall see below, it will be useful to approximate each discontinuous function ${\bf q}_N$ in time
by a piece-wise constant function,  $\bar{\bf q}_N$, so that
$$
\bar{\bf q}_N(t,.)={\bf q}(n\Delta t-,.),\quad t\in [(n-1)\Delta t,n\Delta t),\ n=1,\dots,N,
$$
where ${\bf q}_N(n\Delta t-)$ is the limit from the left of ${\bf q}_N$ at  $n\Delta t$, $n=1,\dots,N$.
By using Lemma \ref{testf}, and by applying the same arguments in the proof of Lemma \ref{konveta}, we get
$$
\bar{\bf q}_N\rightarrow \tilde{\bf q}\ {\rm  uniformly \ on}\   [0,T]\times\Omega.
$$

\subsection{Passing to the limit}

To get to the weak formulation of the coupled problem, take the test functions
$(\psi(t),\bpsi(t))\in {\cal X}_W\times{\cal X}_S$ as the test functions in the
weak formulation of the structure sub-problem (\ref{SProb1}) and integrate 
the weak formulation (\ref{SProb1}) with respect to $t$ from $n\Delta t$ to $(n+1)\Delta t$.
Notice that the construction of the test functions is done in such a way that $(\psi(t),\bpsi(t))$ do not depend on $N$, and
are continuous. Then, consider the weak formulation  (\ref{D1Prob1}) of the fluid sub-problem and take the test 
functions $({\bf q}_N(t),\psi(t))$ (where $\mathbf q_N = \mathbf q \circ A_{{\tau_{\Delta t}\eta_N}}$, $\mathbf q \in {\cal X}_F$).
Integrate the fluid sub-problem (\ref{D1Prob1}) with respect to $t$ from $n\Delta t$ to $(n+1)\Delta t$.
Add the two weak formulations together, and take the sum from $n=0,\dots,N-1$ to get the time integrals over $(0,T)$ as follows:
\begin{equation}
\begin{array}{c}
\displaystyle{\int_0^T\int_{\Omega_F}(1+\tau_{\Delta t}\ \eta_N) \Big (\partial_t\tilde{\bf u}_N\cdot {\bf q}_N+
\frac 1 2(\tau_{\Delta t}{\bf u}_N-{\bf w}_N)\cdot\nabla^{{\tau_{\Delta t}\eta_N}}{\bf u}_N\cdot{\bf q}_N}\\ \\
\displaystyle{-\frac 1 2(\tau_{\Delta t}{\bf u}_N-{\bf w}_N)\cdot\nabla^{{\tau_{\Delta t}\eta_N}}{\bf q}_N\cdot{\bf u}_N\Big )+
\frac{1}{2}\int_0^T \int_{\Omega_F}{v^*_N}{\bf u}_N\cdot{\bf q}_N}\\ \\
\displaystyle{+\int_0^T\int_{\Omega_F}(1+\tau_{\Delta t}\eta_N)2{\bf D}^{{\tau_{\Delta t}\eta_N}}({\bf u_N}):{\bf D}^{{\tau_{\Delta t}\eta_N}}({\bf q}_N)+
\int_0^T\int_0^1\partial_{t} \tilde{v}_N\psi}\\ \\
\displaystyle{+\int_0^T\int_0^1\partial_z\eta_N\partial_z\psi
+\int_0^T\int_{\Omega_S}\partial_t\tilde{\bV}_N\cdot{\bpsi}
+\int_0^T\int_{\Omega_S}a_S(\bd_N,\bpsi)}
\\ \\
\displaystyle{ =\int_0^TP_{in}^Ndt\int_0^1q_z(t,0,r)dr-\int_0^TP_{out}^Ndt\int_0^1q_z(t,L,r)dr,}
\end{array}
\label{AproxEqNS}
\end{equation} 
with
\begin{equation}
\begin{array}{c}
\nabla^{\tau_{\Delta t}\eta}\cdot{\bf u}_N=0,\quad v_N=((u_r)_N)_{|\Gamma},\; \eta_N=(\bd_N)_{|\Gamma},
\\[0.3cm]
{\bf u}_N(0,.)={\bf u}_0,\ \eta(0,.)_N=\eta_0,\ v_N(0,.)=v_0.
\end{array}
\label{AproxEqNS}
\end{equation} 
Here $\tilde{\mathbf u}_N$, $\tilde{v}_N$ and $\tilde{\bV}_N$ are the piecewise linear functions defined in \eqref{tilde},
$\tau_{\Delta t}$ is the shift in time by $\Delta t$ to the left, defined in \eqref{shift},
$\nabla^{{\tau_{\Delta t}\eta_N}}$ is the transformed gradient via the ALE mapping $A_{{\tau_{\Delta t}\eta_N}}$,
defined in \eqref{nablaeta}, and $v_N^*$, $\mathbf u_N$, $v_N$, $\eta_N$, $\bd_N$ and $\bV_N$ are defined in \eqref{aproxNS}.

Using the convergence results obtained for the approximate solutions in Section~\ref{sec:convergence}, and 
the convergence results just obtained for the test functions $\mathbf q_N$,
we can pass to the limit directly in all the terms except in the term that contains $\partial_t\tilde{\mathbf u}_N$.
To deal with this term we notice that, 
since ${\bf q}_N$ are smooth on sub-intervals $(j\Delta t,(j+1)\Delta t)$, we can use integration by parts on these sub-intervals to obtain:
$$
\int_0^T\int_{\Omega_F}(1+\tau_{\Delta t}\eta_N)\partial_t\tilde{\bf u}_N\cdot {\bf q}_N=
\sum_{j=0}^{N-1}\int_{j\Delta t}^{(j+1)\Delta t}\int_{\Omega_F}(1+\eta^j_N)\partial_t\tilde{\bf u}_N\cdot {\bf q}_N
$$
$$
=\sum_{j=0}^{N-1}\Big (-\int_{j\Delta t}^{(j+1)\Delta t}\int_{\Omega_F}(1+\tau_{\Delta t}\eta_N)\tilde{\bf u}_N\cdot\partial_t{\bf q}_N
$$
\begin{equation}\label{two_terms}
+\int_{\Omega_F}(1+\eta^{j+1}-\eta^{j+1}+\eta^j){\bf u}_N^{j+1}\cdot{\bf q}_N((j+1)\Delta t-)
-\int_{\Omega_F}(1+\eta^j){\bf u}_N^{j}\cdot{\bf q}_N(j\Delta t+)\Big ).
\end{equation}
Here, we have denoted by ${\bf q}_N((j+1)\Delta t-)$ and ${\bf q}_N(j\Delta t+)$
the limits from the left and right, respectively, of ${\bf q}_N$ at the appropriate points.

The integral involving $\partial_t{\bf q}_N$ can be simplified by recalling that
$\mathbf q_N = \mathbf q \circ A_{\eta_N}$, where $\eta_N$ are constant on each sub-interval $(j\Delta t,(j+1)\Delta t)$.
Thus, by the chain rule, we see that $\partial_t{\bf q}_N= \partial_t \mathbf q$ on $(j\Delta t,(j+1)\Delta t)$.
After summing over all $j = 0,...,N-1$ we obtain
$$
- \sum_{j=0}^{N-1}\int_{j\Delta t}^{(j+1)\Delta t}\int_{\Omega_F}(1+\tau_{\Delta t}\eta_N)\tilde{\bf u}_N\cdot\partial_t{\bf q}_N
=-\int_0^T\int_{\Omega_F}(1+\tau_{\Delta t}\eta_N)\tilde{\bf u}_N\cdot\partial_t{\bf q}.
$$
To deal with the last two terms in \eqref{two_terms} we calculate
$$
\sum_{j=0}^{N-1}\Big(\int_{\Omega_F}(1+\eta^{j+1}_N-\eta^{j+1}_N+\eta^j_N){\bf u}_N^{j+1}\cdot{\bf q}_N((j+1)\Delta t-)
-\int_{\Omega_F}(1+\eta^j_N){\bf u}_N^{j}\cdot{\bf q}_N(j\Delta t+)\Big)
$$
$$
=\sum_{j=0}^{N-1} \int_{\Omega_F}\left((1+\eta^{j+1}_N){\bf u}_N^{j+1}\cdot{\bf q}_N((j+1)\Delta t-)
-(\eta^{j+1}_N-\eta^j_N) {\bf u}_N^{j+1}\cdot{\bf q}_N((j+1)\Delta t-)\right)
$$
$$
- \int_{\Omega}(1+\eta_0) {\bf u}_0\cdot \mathbf q(0)
-\sum_{j=1}^{N-1} \int_{\Omega_F}(1+\eta^j_N){\bf u}_N^{j}\cdot{\bf q}_N(j\Delta t+)\Big)
$$
Now, we can write $(\eta^{j+1}-\eta^j)$ as $\displaystyle{v^{j+\frac 1 2}\Delta t}$, and rewrite the summation indexes
in the first term to obtain that the above expression is equal to
$$
=\sum_{j=1}^{N} \int_{\Omega_F}(1+\eta^{j}_N){\bf u}_N^{j}\cdot{\bf q}_N(j\Delta t-)
-\int_0^T\int_{\Omega_F}v^*_N{\bf u}_N\cdot \bar{\bf q}_N
- \int_{\Omega_F}(1+\eta_0) {\bf u}_0\cdot \mathbf q(0)
-\sum_{j=1}^{N-1}\int_{\Omega_F}(1+\eta_N^j){\bf u}^j_N \cdot {\bf q}_N(j\Delta t+).
$$
Since the test functions have compact support in $[0,T)$, the value of the first term at $j = N$ is zero, and so 
we can combine the two sums to obtain
$$
=\sum_{j=1}^{N} \int_{\Omega_F}(1+\eta^{j}_N){\bf u}_N^{j}\cdot \left( {\bf q}_N(j\Delta t-)-{\bf q}_N(j\Delta t+)\right)
- \int_{\Omega_F}(1+\eta_0) {\bf u}_0\cdot \mathbf q(0)-\int_0^T\int_{\Omega_F}v^*_N{\bf u}_N\cdot \bar{\bf q}_N.
$$
Now we know how to pass to the limit in all the terms expect the first one. 
We continue to rewrite the first expression by using the Mean Value Theorem to obtain:
$$
{\bf q}_N(j\Delta t-,z,r)-{\bf q}_N(j\Delta t+,z,r)={\bf q}(j\Delta t,z,(1+\eta_N^j)r)-{\bf q}(j\Delta t,z,(1+\eta_N^{j+1})r)=
$$
$$
=\partial_r{\bf q}(j\Delta t,z,\zeta)r(\eta_N^j-\eta_N^{j+1})=-\Delta t\partial_r{\bf q}(j\Delta t,z,\zeta)v_N^{j+\frac 1 2}r.
$$ 
Therefore we have:
$$
\sum_{j=1}^{N-1}\int_{\Omega_F}(1+\eta_N^j){\bf u}^j_N\big ({\bf q}(j\Delta t-)-{\bf q}(j\Delta t+))=-
\int_0^{T-\Delta t}\int_{\Omega_F}(1+\eta_N){\bf u}_Nr\tau_{-\Delta t}v^*_N\partial_r\bar{\bf q}.
$$
We can now pass to the limit in this last term to obtain:
$$
\int_0^{T-\Delta t}\int_{\Omega_F}(1+\eta_N){\bf u}_Nr\tau_{-\Delta t}v^*_N\partial_r\bar{\bf q}\rightarrow 
\int_0^T\int_{\Omega_F}(1+\eta){\bf u}r\partial_t\eta\partial_r{\bf q}.
$$
Therefore, by noticing that $\partial_t\tilde{\bf q}=\partial_t {\bf q}+r\partial_t\eta\partial_r{\bf q}$
we have finally obtained
$$
\int_0^T\int_{\Omega_F}(1+\tau_{\Delta t}\eta_N)\partial_t\tilde{\bf u}_N\cdot {\bf q}_N\rightarrow -\int_0^T\int_{\Omega_F}(1+\eta){\bf u}\cdot\partial_t\tilde{\bf q}
-\int_0^T\int_{\Omega_F}\partial_t\eta{\bf u}\cdot\tilde{\bf q}
$$
$$
-\int_{\Omega_F}(1+\eta_0){\bf u}_0\cdot\tilde{\bf q}(0),
$$
where we recall that $\tilde{\mathbf q }= \mathbf q \circ A_\eta$.

Thus, we have shown that the limiting functions $\mathbf u$, $\eta$ and $\bd$ satisfy the weak form of 
problem \eqref{FSIeqRef}-\eqref{IC_ALE} in the sense of Definition \ref{DefWSRef},
for all test functions that belong to a dense subset of ${\cal Q}^\eta$.
By density arguments, we have, therefore, shown the main result of this manuscript:

\begin{theorem}{\rm{\bf{(Main Theorem)}}}
Suppose that the initial data   $v_0\in L^2(0,1)$, ${\bf u}_0\in L^2(\Omega_{\eta_0})$, $\bV_0\in L^2(\Omega_S)$, $\bd_0\in H^1(\Omega_S)$,
and $\eta_0\in H^1_0(0,1)$ are
such that $1+\eta_0(z)>0$, $z\in [0,1]$ and compatibility conditions \eqref{CC} are satisfied. 
Furthermore, let $P_{in}$, $P_{out}\in L^{2}_{loc}(0,\infty)$. 

Then, there exist a $T>0$ and a weak solution $(\mathbf u,\eta,\bd)$ of
problem \eqref{FSIeqRef}-\eqref{IC_ALE} 
(or equivalently problem \eqref{NS}-\eqref{CC}) on $(0,T)$ in the sense of Definition \ref{DefWSRef} (or equivalently Definition \ref{DefWS}),
such that the following energy estimate is satisfied:
\begin{equation}
E(t)+\int_0^tD(\tau)d\tau \leq E_0+C(\|P_{in}\|_{L^{2}(0,t)}^2+\|P_{out}\|_{L^{2}(0,t)}^2),\quad t\in [0,T],
\label{EE1}
\end{equation}
where $C$ depends only on the coefficients in the problem, $E_0$ is the kinetic energy of initial data,
and $E(t)$ and $D(t)$ are given by
\begin{eqnarray*}
E(t) &=& \frac{1}{2}\|{\bf u}\|^2_{L^2(\Omega_F)}+\frac{1}{2}\|\partial_t\eta\|^2_{L^2(0,1)}+\frac 1 2\|\bd\|^2_{L^2(\Omega_S)}
+ \frac{1}{2}\big (\|\partial_z\eta\|^2_{L^2(0,1)}+a_S(\bd,\bd)\big ),\\
D(t)&=&\|{\bf D}({\bf u})\|^2_{L^2(\Omega_{\eta}(t)))}.
\end{eqnarray*}
\vskip 0.1in
 Furthermore, one of the following is true: 
\begin{equation}\label{ExistenceWS}
either \ T=\infty \  \  {or} \  \ {\lim_{t\rightarrow T}\min_{z\in [0,1]}(1+\eta(z))=0}.
\end{equation}
\end{theorem}
\proof
It only remains to prove the last assertion, which states that our result is either global in time, or, in case the walls of the cylinder 
touch each other, our existence result holds until the time of touching. However, 
the proof of this argument follows the same reasoning as the proof of the Main Theorem in \cite{BorSun},
and the proof of the main result in  \cite{CDEM}, p.~397-398. We avoid repeating those arguments here, and
refer the reader to references \cite{BorSun,CDEM}.\qed

\section{Conclusions}

In this manuscript we proved the existence of a weak solution to a FSI problem in which the structure consists of two layers:
a thin layer modeled by the linear wave equation, and a thick layer modeled by the equations of linear elasticity.
The thin layer acts as a fluid-structure interface with mass. 
An interesting new feature of this problem is the fact that the presence of a thin structure with mass
regularizes the solution of this FSI problem. More precisely,
the energy estimates presented in this work
show that the thin structure inertia regularizes the evolution of the thin structure, which affects the solution of the 
entire coupled FSI problem.
Namely, if we were considering a problem in which 
the structure consisted of only one layer, modeled by the equations of linear elasticity, 
from the energy estimates we would not be able to conclude that the fluid-structure interface is even continuous,
since the displacement $\bf d$ of the thick structure would be in $H^{1/2}(\Gamma)$ at the interface.
With the presence of a thin elastic fluid-structure interface with mass (modeled by the wave equation), the energy 
estimates imply that the displacement of the thin interface  is in $H^1(\Gamma)$, which,
  due to the Sobolev embeddings, implies that the interface is H\"{o}lder continuous $C^{0,1/2}(\Gamma)$.

This is reminiscent of the results by Hansen and Zuazua \cite{HansenZuazua} in which the presence of a point mass at the interface between two
linearly elastic strings with solutions in asymmetric spaces (different regularity on
each side) allowed the proof of well-posedness due to the regularizing effects by the
point mass. For a reader with further interest in the area of simplified coupled problems, we also mention
\cite{KochZauzua,RauchZhangZuazua,ZhangZuazuaARMA07}.

Further research by the authors in the direction of simplified coupled problems that  shed light
on the physics of parabolic-hyperbolic coupling with point mass,  is under way \cite{SunBorHYP,BorSun1D}.
 Our preliminary results in \cite{BorSun1D} indicate that the regularizing feature of the interface with mass is not only a consequence 
 of our mathematical methodology, but a physical property of this complex system.
 
\if 1 = 0
This can be proved by using similar argument as in \cite{CDEM} p.~397-398. 
For the sake of completeness we present the arguments here. 

Let $(0,T_1)$, $T_1 > 0$, be the interval on which we have constructed 
our solution $({\bf u},\eta,\bd)$,
and let $m_1=\min_{(0,T_1)\times (0,1)} (1+\eta)$. From Lemma~\ref{eta_bound} we know that $m_1>0$. Furthermore, since
$\eta\in W^{1,\infty}(0,T;L^2(0,1))\cap L^{\infty}(0,T;H^1_0(0,1))$ and ${\bf u}\in L^{\infty}(0,T;L^2(\Omega_F))$,
we can take $T_1$ such that $\eta(T_1)\in H^1_0(0,1)$, $\partial_t\eta(T_1)\in L^2(0,1)$ and ${\bf u}(T_1)\in L^2(\Omega_F)$,
$\bd(T_1)\in H^1(\Omega_S)$ and $\bV(T_1)\in L^2(\Omega_S)$. 
We can now use the first part of Theorem \ref{ExistenceWS} to prolong the solution $({\bf u},\eta)$ to the interval $(0,T_2)$, $T_2>T_1$.
By iteration, we can continue the construction of our solution to the interval $(0,T_{k})$, $k\in\N$, where $(T_k)_{k\in\N}$ is an increasing sequence.
We set $m_k=\min_{(0,T_k)\times (0,1)} (1+\eta)>0$. 

Since $m_k > 0$ we can continue the construction further. Without loss of generality we could choose
a $T_{k+1} > T_k$  so that $m_{k+1}\geq\frac{m_{k}}{2}$. 
From (\ref{EE1}) and (\ref{compact_embedding}), by taking $\alpha=3/4$, we have that 
the displacement $\eta$ is H\"{o}lder continuous in time, namely,
$$
\|\eta\|_{C^{0,1/4}(0,T_{k+1};C[0,1])}\leq C(T_{k+1}).
$$
Therefore, the following estimate holds:
$$
1+\eta(T_{k+1},z)\geq 1+\eta(T_{k},z)-C(T_{k+1})(T_{k+1}-T_{k})^{1/4}\geq m_{k}-C(T_{k+1})(T_{k+1}-T_{k})^{1/4}.
$$
For a $T_{k+1}$ chosen so that $m_{k+1}\geq\frac{m_{k}}{2}$
this estimate implies
\begin{equation}\label{prolong}
\displaystyle{T_{k+1}-T_{k}\geq \frac{m_{k}^4}{16C(T_{k+1})^4}},\; k\in\N.
\end{equation}

Now, let us take $T^*=\sup_{k\in\N}T_k$ and set $m^*=\min_{(0,T^*)\times (0,1)} (1+\eta)$. Obviously, $m_k\geq m^*$, $k\in\N$.
There are two possibilities. Either $m^* = 0$, or $m^* > 0$. 
If $m^*=0$, this means that 
$\displaystyle{\lim_{t\rightarrow T}\min_{z\in [0,L]}(1+\eta(z))=0}$, and the second statement in the theorem is proved.
If $m^* > 0$, we need to show that $T^* = \infty$.
To do that, notice that (\ref{EE1}) gives the form of the constant $C(T)$ which is a non-decreasing function of $T$.
Therefore, we have $C(T_k)\leq C(T^*)$, $\forall k\in\N$. Using this observation and that fact that $m_k\geq m^*$, $k\in\N$,
estimate (\ref{prolong}) implies
$$
T_{k+1}-T_{k}\geq \frac{(m^*)^2}{2C(T^*)^2},\; \forall k\in\N.
$$
Since this holds for all $k \in \N$, we have that $T^*=\infty$.
\fi

\if 1 = 0
\section{Appendix: Trace Theorem for domains which are locally subgraph of H\" older continuous function}
\label{sec:Trace}

In this section we state special cases of results from \cite{BorisTrag} using notation from this paper.

\begin{lemma}[\cite{BorisTrag}, Lemma 3.3.]\label{TraceLemma}
Let $\eta\in C^{0,\alpha}$, $0<\alpha<1$. 
Let $u\in H^1(\Omega_{\eta})$ and let $0<s<\alpha$. Then ${\tilde u}\in W(0,1;s)$, where 
$$
{\tilde u}(\tilde{r},\tilde{z})=u(\tilde{z},(1+\eta(\tilde{z}))\tilde{r}),\quad (\tilde{z},\tilde{r})\in\Omega_F
$$
and
$$
W(0,1;s)=\{f:f\in L^2(0,1;H^s(0,1))),\; \partial_{\tilde{r}} f\in L^2(0,1;L^2(0,1))\}.
$$
\end{lemma}
\begin{theorem}[\cite{BorisTrag}, Theorem 3.1.]\label{TraceTm}
Let $\alpha<1$ and let $\eta$ be such that 
$$
\eta\in C^{0,\alpha}(0,1),\; \eta({\bf x})\geq \eta_{min}>-1,\; {\bf x}\in [0,1],\;  \eta(0)=\eta(1)=1.
$$
Then the operator 
$$
\gamma_{\eta}:C^1(\overline{\Omega_{\eta}})\ni u\mapsto u(\tilde{z},1+\eta(\tilde{z}))\in C(\Gamma)
$$ 
can be extended by continuity to a linear operator from $H^1(\Omega_{\eta})$ to 
$H^s(\omega)$, $0\leq s< \frac{\alpha}{2}$.
\end{theorem}
\fi

\noindent
{\bf Acknowledgements.} The authors would like to thank Prof. Enrique Zuazua for pointing out the references
\cite{HansenZuazua,KochZauzua,RauchZhangZuazua,ZhangZuazuaARMA07}.
Furthermore, the authors would like to thank the following research support:
Muha's research was supported in part by the Texas Higher Education Board 
under grant ARP 003652-0023-2009, and by ESF OPTPDE -
Exchange Grant 4171;  \v{C}ani\'{c}'s research was supported by the National Science Foundation under grants
DMS-1311709, DMS-1109189, DMS-0806941, and by the Texas Higher Education Board 
under grant ARP 003652-0023-2009.

\end{document}